\documentclass[12pt]{article}
\usepackage{amsmath,amsfonts,amssymb,amsthm,enumitem}
\input amssym.def
\usepackage[a4paper,top=2cm,bottom=2cm,left=2.6cm,right=2.6cm]{geometry}
\usepackage[hidelinks]{hyperref}
\allowdisplaybreaks

\numberwithin{equation}{section}
\begin{document}
	\def \Z{\Bbb Z}
	\def \C{\Bbb C}
	\def \R{\Bbb R}
	\def \Q{\Bbb Q}
	\def \N{\Bbb N}
	
	\def \A{{\mathcal{A}}}
	\def \D{{\mathcal{D}}}
	\def \E{{\mathcal{E}}}
	\def \E{{\mathcal{E}}}
	\def \H{\mathcal{H}}
	\def \S{{\mathcal{S}}}
	\def \wt{{\rm wt}}
	\def \tr{{\rm tr}}
	\def \span{{\rm span}}
	\def \Res{{\rm Res}}
	\def \Der{{\rm Der}}
	\def \End{{\rm End}}
	\def \Ind {{\rm Ind}}
	\def \Irr {{\rm Irr}}
	\def \Aut{{\rm Aut}}
	\def \GL{{\rm GL}}
	\def \Hom{{\rm Hom}}
	\def \mod{{\rm mod}}
	\def \ann{{\rm Ann}}
	\def \ad{{\rm ad}}
	\def \rank{{\rm rank}\;}
	\def \<{\langle}
	\def \>{\rangle}
	
	\def \g{{\frak{g}}}
	\def \h{{\hbar}}
	\def \k{{\frak{k}}}
	\def \sl{{\frak{sl}}}
	\def \gl{{\frak{gl}}}
	
	\def \be{\begin{equation}\label}
		\def \ee{\end{equation}}
	\def \bex{\begin{example}\label}
		\def \eex{\end{example}}
	\def \bl{\begin{lem}\label}
		\def \el{\end{lem}}
	\def \bt{\begin{thm}\label}
		\def \et{\end{thm}}
	\def \bp{\begin{prop}\label}
		\def \ep{\end{prop}}
	\def \br{\begin{rem}\label}
		\def \er{\end{rem}}
	\def \bc{\begin{coro}\label}
		\def \ec{\end{coro}}
	\def \bd{\begin{de}\label}
		\def \ed{\end{de}}
	
	\newcommand{\m}{\bf m}
	\newcommand{\n}{\bf n}
	\newcommand{\nno}{\nonumber}
	\newcommand{\nord}{\mbox{\scriptsize ${\circ\atop\circ}$}}
	\newtheorem{thm}{Theorem}[section]
	\newtheorem{prop}[thm]{Proposition}
	\newtheorem{coro}[thm]{Corollary}
	\newtheorem{conj}[thm]{Conjecture}
	\newtheorem{example}[thm]{Example}
	\newtheorem{lem}[thm]{Lemma}
	\newtheorem{rem}[thm]{Remark}
	\newtheorem{de}[thm]{Definition}
	\newtheorem{hy}[thm]{Hypothesis}
	\makeatletter \@addtoreset{equation}{section}
	\def\theequation{\thesection.\arabic{equation}}
	\makeatother \makeatletter
	\newcommand{\corrfootnote}{\footnotemark[0]} 
	\renewcommand{\thefootnote}{\fnsymbol{footnote}} 

	\begin{center}
		{\Large \bf Twisted Verlinde Formula for Finite Cyclic Orbifolds of Vertex Operator Algebras}
	\end{center}
	
	\begin{center}
		{Wenbin Jiang \\
			Department of Mathematics \\
			Shanghai University, Shanghai 200444, China \\
			E-mail address: jiangwenbin8899@163.com}
	\end{center}
	
	\begin{center}
		{Jiancai Sun\footnote[1]{Corresponding author.} \\
			Department of Mathematics and Newtouch Center for Mathematics \\
			Shanghai University, Shanghai 200444, China \\
			E-mail address: jcsun@shu.edu.cn}
	\end{center}

	\begin{abstract}
		Let $V$ be a simple, rational, $C_2$-cofinite and self-dual vertex operator
		algebra of CFT type, and let $G$ be a finite cyclic automorphism group of
		$V$. We establish a twisted Verlinde formula for cyclic groups of arbitrary
		finite order, expressing fusion coefficients among twisted $V$-modules in
		terms of the $S$-matrix in the orbifold theory. As an application, we compute fusion rules
		among twisted modules for the permutation orbifold
		$\left(V^{\otimes k}\right)^{\langle g\rangle}$ associated with a $k$-cycle
		$g=(1\,2\,\cdots\,k)$.
	\end{abstract}

	\section{Introduction}
	
	The relation between modular transformations and fusion rules was first formulated by E. Verlinde in rational conformal field theory \cite{v}.
	In the framework of vertex operator algebras, the modular invariance of trace functions was obtained, which led to a finite-dimensional representation of $SL(2,\mathbb Z)$ on the space spanned by irreducible characters \cite{z}.
	The Verlinde formula for vertex operator algebras was subsequently formulated and proved, revealing that the fusion matrices of irreducible $V$-modules were simultaneously diagonalized by the modular $S$-matrix and that the fusion coefficients were determined by its entries \cite{h1}.
	Under the assumptions of rationality, $C_2$-cofiniteness, self-duality and CFT type, the category of ordinary $V$-modules, together with the tensor product defined through intertwining operators, was further shown to possess the structure of a modular tensor category \cite{h2}.
	As a consequence, the fusion coefficients of irreducible modules can be expressed in terms of the modular $S$-matrix.

	Let $V$ be a vertex operator algebra and $G$ be a finite subgroup of
	$\operatorname{Aut}(V)$. Orbifold theory considers the fixed-point
	subalgebra $V^G$ together with the $g$-twisted $V$-modules for $g\in G$.
	In conformal field theory, the structure and fusion rules of orbifold models
	were studied through modular transformations of one-loop characters
	\cite{dvvv}. For vertex operator algebras, the associative algebra $A_g(V)$
	was introduced in \cite{dlm}, establishing a correspondence between its
	irreducible modules and irreducible admissible $g$-twisted $V$-modules.
	Moreover, for commuting $g,h\in G$, the trace functions associated with
	$h$-stable irreducible $g$-twisted modules span a finite-dimensional
	$SL(2,\mathbb Z)$-invariant space \cite{dlm3}. The matrix
	corresponding to $\tau\mapsto-\frac{1}{\tau}$ is referred to as the $S$-matrix in
	the orbifold theory.

	To relate the orbifold $S$-matrix to the module category of $V^G$, one
	needs a precise description of the irreducible $V^G$-modules arising from
	twisted $V$-modules. Let $M$ be an irreducible $g$-twisted $V$-module and
	let $G_M$ be its stabilizer in $G$. The projective action of $G_M$ on $M$
	gives a decomposition of $M$ into irreducible $V^G$-modules indexed by the
	irreducible projective characters of $G_M$ \cite{dy,mt}. The corresponding
	components are nonzero and irreducible, and their isomorphism classes are
	determined by the character and the $G$-orbit of $M$. The
	$C_2$-cofiniteness of cyclic fixed-point subalgebras was established in
	\cite{m}, while the regularity of $V^G$ for finite solvable automorphism
	groups, under the standard assumptions on $V$, was proved in \cite{cm}.
	Assuming that $V^G$ is regular and that every nontrivial irreducible twisted
	$V$-module has positive conformal weight, it was proved in \cite{drx1} that
	every irreducible $V^G$-module occurs as an irreducible $V^G$-submodule of
	some irreducible $g$-twisted $V$-module, where $g\in G$. Subsequently, the
	restricted $S$-matrix of $V^G$, indexed by the irreducible
	$V^G$-submodules arising from twisted $V$-modules, was expressed in
	\cite{drx2} in terms of the $S$-transformation coefficients of twisted trace
	functions. When these submodules form a complete set of irreducible
	$V^G$-modules, the restricted $S$-matrix coincides with the full
	$S$-matrix of $V^G$.

	The formulation of a twisted Verlinde formula requires a precise definition
	of fusion rules among twisted modules. Twisted intertwining operators were
	introduced in \cite{xxp} and later extended to modules associated with
	finite order automorphisms that need not commute \cite{h3}. Using the theory
	of algebra objects in braided tensor categories \cite{ko,hkl,ckm}, it was
	proved in \cite{dlxy} that, for a regular self-dual vertex operator algebra
	of CFT type and a finite abelian automorphism group $G$, fusion products
	among twisted modules exist and are associative. Their multiplicities agree
	with the dimensions of the corresponding spaces of twisted intertwining
	operators. In permutation orbifold theory, the $S$-matrix for cyclic
	permutation orbifolds was expressed explicitly in terms of the $S$-matrix
	of $V$ in \cite{dxy}, and these results were subsequently used in
	\cite{dxy2} to determine fusion products among twisted modules for cyclic
	permutation groups.

	The categorical formulation of twisted fusion rules suggests a corresponding
	Verlinde formula. A Verlinde formula for $G$-crossed modular fusion categories
	was established in \cite{dem,t}, expressing fusion coefficients in terms of
	crossed $S$-matrices. For vertex operator algebras, the category of twisted
	modules is expected to carry such a $G$-crossed structure \cite{mc}, and its
	crossed $S$-matrix is expected to agree, up to normalization, with the
	$S$-matrix arising from twisted trace functions.

	An early twisted analogue of the Verlinde formula was obtained for
	automorphism groups of order two in conformal field theory \cite{bfs}. In
	the setting of vertex operator algebras, the result was extended in
	\cite{dl} to cyclic automorphism groups of prime order, expressing fusion
	coefficients among twisted modules in terms of the $S$-matrix in the
	orbifold theory. The formula was further applied to affine vertex operator
	algebras at positive integral levels, leading to a twisted analogue of the
	Kac--Walton formula.

	Motivated by the prime order result in \cite{dl}, it is natural to investigate
	whether the twisted Verlinde formula remains valid when $G$ is a finite
	cyclic group whose order is not necessarily prime. In the present paper, we
	establish the formula for cyclic automorphism groups of arbitrary finite
	order. The resulting formula expresses the fusion coefficients among
	twisted modules entirely in terms of the $S$-matrix in the orbifold theory
	and provides a uniform description of twisted fusion rules for finite
	cyclic orbifolds.
	
	In this paper, we assume that $V$ is a simple, rational, $C_2$-cofinite and
	self-dual vertex operator algebra of CFT type, and that
	$\sigma\in\operatorname{Aut}(V)$ has finite order. Set $G=\langle\sigma\rangle.$
	For $g,h\in G$, let $\mathcal{M}(g,h)$ denote the set of isomorphism classes
	of irreducible $h$-stable $g$-twisted $V$-modules. Given
	$g_1,g_2\in G$, set
	\begin{eqnarray*}
		\mathcal{M}(1;g_1,g_2)
		=
		\mathcal{M}(1,g_1)\cap\mathcal{M}(1,g_2).
	\end{eqnarray*}
	Suppose that $g_1,g_2,g_3\in G$ satisfy $g_1\neq1$ or $g_2\neq1$, and $g_3=g_1g_2$.
	For any $M^i\in\mathcal{M}(g_i)$, $i=1,2,3$, we prove that
	\begin{eqnarray*}
		N^{M^3}_{M^1,M^2}
		=
		\sum_{W\in\mathcal{M}(1;g_1,g_2)}
		\frac{
			S_{M^1,W}S_{M^2,W}
			\overline{S_{M^3,W}}
		}{
			S_{V,W}
		}.
	\end{eqnarray*}
	Thus, the fusion coefficients among twisted modules are determined by the
	$S$-matrix arising from the modular transformations of twisted trace
	functions.

	The argument is based on the fixed-point vertex operator algebra $V^G$. We
	first use the classification of irreducible $V^G$-modules in terms of
	$G$-orbits of irreducible twisted $V$-modules and characters of their
	stabilizer subgroups, together with the relation between the $S$-matrix of
	$V^G$ and the $S$-matrix in the orbifold theory. We then establish a
	relation between the fusion coefficients of twisted $V$-modules and those
	of $V^G$-modules. Let
	\begin{eqnarray*}
		\mathcal{F}:\mathcal{C}_{V^G}\longrightarrow\operatorname{Rep}(V)
	\end{eqnarray*}
	be the induction functor associated with the commutative algebra object $V$
	in $\mathcal{C}_{V^G}$.
	For an irreducible twisted $V$-module $M$ and an irreducible $V^G$-submodule
	$M_\chi$ appearing in the decomposition of $M$ as a $V^G$-module, we prove that
	$\mathcal{F}(M_\chi)$ is the direct sum of the twisted modules in the
	$G$-orbit of $M$.
	This yields explicit relations between the fusion
	coefficients of twisted $V$-modules and those of $V^G$-modules.
	
	The main difficulty is the presence of nontrivial proper stabilizer
	subgroups. Consequently, irreducible twisted modules may lie in $G$-orbits
	of different cardinalities, and their decompositions as $V^G$-modules are
	indexed by irreducible characters of different stabilizer subgroups. This
	requires a systematic treatment of stabilizers, orbit representatives, and
	the character sums appearing in the reduction of the Verlinde formula for
	$V^G$ to the twisted Verlinde formula.

	We first prove the main formula when $|G|=p_1p_2,$ where $p_1$ and $p_2$ are primes, not necessarily distinct.
	The proof is divided
	according to the stabilizer subgroups of the three twisted modules. The
	ordinary Verlinde formula for $V^G$, the orthogonality relations for
	characters, the invariance of the orbifold $S$-matrix, and the lengths of
	the relevant $G$-orbits are used in the computation.
	For a general cyclic group, we write
	$|G|=p_1p_2\cdots p_s$,
	where $s$ is a positive integer and the primes $p_1,p_2,\ldots,p_s$
	are not necessarily distinct. We then use induction on the number of
	prime factors.

	As an application,
	we study the permutation orbifolds
	$\left(V^{\otimes k}\right)^{\langle g\rangle}$ associated with the
	$k$-cycle $g=(1\,2\,\cdots\,k)$. Using the twisted Verlinde formula together
	with the $S$-matrix entries of permutation-twisted modules, we derive explicit
	fusion rules between $g^s$-twisted and $g^r$-twisted modules, as well as
	between untwisted and $g^s$-twisted modules, for arbitrary $1\leq r,s\leq k-1$.
	The formula gives a direct expression for arbitrary $r$ and $s$, while
	the non-coprime case in \cite{dxy2} is reduced to fusion products
	associated with shorter cycles.
	By repeatedly applying the twisted Verlinde formula,
	we also obtain an explicit formula for the  fusion product
	$\left(T_g^0\right)^{\boxtimes k}$, which agrees with
	Theorem~5.1 of \cite{dxy2} and provides an alternative proof. Finally, we specialize
	these formulas to the rank-one lattice vertex operator algebra
	$V_{\mathbb{Z}\alpha}$ with $\langle\alpha,\alpha\rangle=2n$. For $k=4$ and
	$n=1$, the resulting fusion rules agree with those in Section~6 of
	\cite{dxy2}.

	The paper is organized as follows. In Section 2, we recall the basic
	definitions and results about twisted modules, modular invariance,
	twisted fusion rules and the ordinary Verlinde formula. In Section 3, we
	review the classification of irreducible $V^G$-modules and the relation
	between the modular $S$-matrix of $V^G$ and the $S$-matrix in the orbifold
	theory. In Section 4, we establish a relationship between the fusion
	rules of $V$ and  $V^G$. In Section 5, we prove the twisted Verlinde formula for
	arbitrary finite cyclic groups. In Section 6, we apply the twisted Verlinde formula to the
	permutation orbifold
	$\left(V^{\otimes k}\right)^{\langle g\rangle}$ associated with a $k$-cycle
	$g=(1\,2\,\cdots\,k)$.

	Throughout this paper, we use the standard notation $\mathbb{C}$,
	$\mathbb{Q}$, and $\mathbb{Z}$ for the sets of complex numbers,
	rational numbers, and integers, respectively. We denote by
	$\mathbb{Z}_{+}$ the set of positive integers and by $\mathbb{N}$ the
	set of nonnegative integers.

	$\mathbf{Acknowledgment}$: The authors would like to thank Professors Haisheng Li and Xiangyu Jiao for their valuable discussions and helpful suggestions. J.-C. Sun would like to thank the support of National Natural Science Foundation of China
	(Nos. 12071276, 11931009, and 12226402) and Science and Technology Commission of Shanghai Municipality (No. 25ZR1401126).

	\section{Preliminaries}
	In this section, we recall some basic notions and results on vertex
	operator algebras and their orbifold theory that will be used
	throughout this paper. We begin with the definitions of twisted
	modules and review several standard properties of vertex operator
	algebras. We then discuss the modular invariance of trace functions
	associated with twisted modules, together with the corresponding
	$S$-matrices. Finally, we recall the fusion rules for twisted modules
	and the Verlinde formula, which will play an essential role in the
	subsequent sections.

	\subsection{Twisted Modules}
	
	In this subsection, we recall the notions of twisted modules for vertex operator algebras. We also review several conditions that will be used throughout this paper, including rationality, self-duality, $C_2$-cofiniteness and CFT type condition.
	
	Let $x$ be a formal variable, and let $x_0,x_1,x_2$ be independent
	formal variables. The formal delta function is defined by
	\begin{eqnarray}
		\delta(x)=\sum_{n\in\mathbb{Z}}x^n.\label{7}
	\end{eqnarray}
	
	For a linear space $V$ , we define several linear spaces of series with coefficients
	in $V$:
	\begin{equation}
		V[x]
		=
		\left\{
		\sum_{n\in\mathbb{N}}v_nx^n
		\mid
		v_n\in V,\ \text{all but finitely many }v_n=0
		\right\},\label{1}
	\end{equation}
	
	\begin{equation}
		V[x,x^{-1}]
		=
		\left\{
		\sum_{n\in\mathbb{Z}}v_nx^n
		\mid
		v_n\in V,\ \text{all but finitely many }v_n=0
		\right\},\label{2}
	\end{equation}
	
	\begin{equation}
		V[[x]]
		=
		\left\{
		\sum_{n\in\mathbb{N}}v_nx^n
		\mid
		v_n\in V
		\right\},\label{3}
	\end{equation}
	
	\begin{equation}
		V[[x,x^{-1}]]
		=
		\left\{
		\sum_{n\in\mathbb{Z}}v_nx^n
		\mid
		v_n\in V
		\right\},\label{4}
	\end{equation}
	
	\begin{equation}
		V((x))
		=
		\left\{
		\sum_{n\in\mathbb{Z}}v_nx^n
		\mid
		v_n\in V,\ v_n=0\text{ for sufficiently negative }n
		\right\}.\label{5}
	\end{equation}
	
	The following is the definition of a vertex operator algebra (see \cite{b}, \cite{flm},  \cite{fhl}, \cite{ll}).
	\bd{definition 1} { \em A $\emph{vertex operator algebra}$ is a quadruple $(V,Y,\mathbf{1},\omega)$,
		where
		\[
		V=\bigoplus_{n\in\mathbb{Z}}V_{(n)} \quad \text{for} \,v\in V_{(n)},\quad n=\operatorname{wt}v,
		\]
		is a $\mathbb{Z}$-graded vector space satisfying
		\begin{eqnarray*}
			&&\dim V_{(n)}<\infty\quad \text{for} \, n\in\mathbb{Z},\\
			&&V_{(n)}=0 \quad \text{for sufficiently negative} \, n,
		\end{eqnarray*}
		there exists a $\emph{vacuum vector}$ $\mathbf{1}\in V_{(0)}$ and a linear map
		\begin{eqnarray*}
			&&Y(\cdot,x):V\to
			(\operatorname{End}V)[[x,x^{-1}]]\\
			&&\hspace{1.6cm}v \mapsto Y(v, x) = \sum_{n \in \mathbb{Z}} v_n x^{-n-1},
		\end{eqnarray*}
		called the $\emph{vertex operator map}$. The following conditions
		are assumed for $u,v\in V$: First, the $\emph{truncation condition}$ holds:
		\[
		u_nv=0\quad \text{for sufficiently large } n,
		\]	
		that is,
		\[
		Y(u,x)v\in V((x));
		\]
		next, the following $\emph{vacuum property}$ is satisfied:
		\[
		Y(\mathbf{1},x)=\operatorname{Id}_V,
		\]
		where $\operatorname{Id}_V$ denotes the identity operator on $V$;
		moreover, the $\emph{creation property}$ holds:
		\begin{eqnarray*}
			Y(v,x)\mathbf{1}\in V[[x]]\quad \text{and}\quad
			\lim_{x\rightarrow0}Y(v,x)\mathbf{1}=v;
		\end{eqnarray*}
		and the $\emph{Jacobi identity}$ is given by
		\begin{eqnarray}
			&&x_0^{-1}\delta\left(\frac{x_1-x_2}{x_0}\right)
			Y(u,x_1)Y(v,x_2)-
			x_0^{-1}\delta\left(\frac{x_2-x_1}{-x_0}\right)
			Y(v,x_2)Y(u,x_1)\notag\\
			&&\hspace{0.5cm}=
			x_2^{-1}\delta\left(\frac{x_1-x_0}{x_2}\right)
			Y(Y(u,x_0)v,x_2).\label{6}
		\end{eqnarray}
		There exists a distinguished homogeneous vector
		$\omega\in V_{(2)}$, called the $\emph{conformal vector}$, such that
		\[
		Y(\omega,x)=\sum_{n\in\mathbb{Z}}L(n)x^{-n-2}.
		\]
		The operators $L(n)$ satisfy the $\emph{Virasoro algebra relations}$
		\[
		[L(m),L(n)]
		=(m-n)L(m+n)
		+\frac{1}{12}(m^{3}-m)\delta_{m+n,0}c_V
		\]
		for all $m,n\in\mathbb{Z}$, where $c_V\in\mathbb{C}$ is called the
		$\emph{central charge}$ of $V$.
		Moreover, we also have:
		\[
		L(0)v=nv=(\operatorname{wt}v)v
		\quad \text{for}\,v\in V_{(n)}, n\in\mathbb{Z},
		\]
		and finally, the $L(-1)$-$\emph{derivative property}$ holds:
		\[
		Y(L(-1)v,x)=\frac{d}{dx}Y(v,x)
		\quad \text{for}\,v\in V.
		\]
	} \ed
	
	We now recall the definition of an automorphism of a vertex operator algebra (cf. \cite{flm}, \cite{fhl}, \cite{ll}).
	\bd{definition 2}{\em
		Let $(V,Y,\mathbf{1},\omega)$ be a vertex operator algebra. An
		$\emph{automorphism}$ of $V$ is a linear bijective map
		\[
		g:V\rightarrow V
		\]
		such that
		\[
		g(Y(u,x)v)=Y(g(u),x)g(v) \quad \text{for}\,u,v\in V
		\]
		and
		\[
		g(\mathbf{1})=\mathbf{1},\quad g(\omega)=\omega .
		\]
		Denote by Aut($V$) the group of all automorphisms of $V$.
	}\ed
	
	Let $g$ be an automorphism of $V$ of finite order $T$.
	The eigenspace decomposition of $V$ with respect to the action of $g$ is
	given by
	\begin{eqnarray*}
		V=\bigoplus_{r\in\mathbb{Z}/T\mathbb{Z}}V^r,
	\end{eqnarray*}
	where $V^r=\{v\in V\mid g v=e^{\frac{-2\pi \sqrt{-1}r}{T}} v\}.$ We use $r$ to denote both an integer in the range $0\leq r\leq T-1$ and its residue class modulo $T$.
	
	The following definitions are taken from \cite{dlm}.
	\bd{definition 3}{\em
		A $\emph{weak g-twisted V-module}$ $M$ is a vector space equipped with a
		linear map
		\begin{eqnarray*}
			&&Y_M(\cdot,x):V\to (\operatorname{End}M)[[x^{\frac{1}{T}},x^{-\frac{1}{T}}]]\\
			&&\hspace{1.8cm}v \mapsto Y_M(v, x) = \sum_{n \in \frac{1}{T}\mathbb{Z}} v_n x^{-n-1}\quad(v_n\in\operatorname{End}M),
		\end{eqnarray*}
		which satisfies the following: for all $0\leq r\leq T-1$, $u\in V^r$, $v\in V$ and $w\in M$,
		\begin{eqnarray}
			&&Y_M(u,x)=\sum_{n\in\frac{r}{T}+\mathbb{Z}}u_nx^{-n-1},\label{10a}\\
			&&u_lw=0\quad \text{for}\, l\gg0,\label{8a}\\
			&&Y_M(\mathbf{1},x)=\operatorname{Id}_M,\label{9a}\\
			&&x_0^{-1}\delta\left(\frac{x_1-x_2}{x_0}\right)
			Y_M(u,x_1)Y_M(v,x_2)-x_0^{-1}\delta\left(\frac{x_2-x_1}{-x_0}\right)
			Y_M(v,x_2)Y_M(u,x_1)\notag\\
			&&\hspace{0.5cm}=
			x_2^{-1}
			\left(\frac{x_1-x_0}{x_2}\right)^{-\frac{r}{T}}
			\delta\left(\frac{x_1-x_0}{x_2}\right)
			Y_M(Y(u,x_0)v,x_2).\label{8}
		\end{eqnarray}	
	}\ed

	\bd{definition 5}{\em
		An $\emph{admissible g-twisted V-module}$ is a weak $g$-twisted $V$-module
		$M$ with a $\frac{1}{T}\mathbb{N}$-grading
		\begin{eqnarray}
			M=\bigoplus_{n\in\frac{1}{T}\mathbb{N}}M(n),\label{10}
		\end{eqnarray}
		which satisfies
		\begin{eqnarray}
			v_m M(n)\subseteq M(n+\operatorname{wt}v-m-1)\label{11}
		\end{eqnarray}
		for homogeneous $v\in V$.
	}\ed

	\bd{definition 4}{\em
		An $\emph{(ordinary) g-twisted V-module}$ is a weak $g$-twisted $V$-module $M$ which
		carries a $\mathbb{C}$-grading induced by the spectrum of $L(0)$, where $L(0)$ is the
		component operator of $Y_M(\omega,x)=\sum_{n\in\mathbb{Z}}L(n)x^{-n-2}.$
		Specifically, we have
		\begin{eqnarray*}
			M=\bigoplus_{\lambda\in\mathbb{C}}M_\lambda,
		\end{eqnarray*}
		where $M_\lambda=\{w\in M\mid L(0)w=\lambda w\}.$
		Moreover, it is required that $\dim M_\lambda<\infty$ for all $\lambda$,
		and for any fixed $\lambda_0$, $M_{\frac{n}{T}+\lambda_0}=0$
		for all sufficiently small integers $n$.
	}\ed

	\br{r1}{\em  Following \cite{dlm3}, we interchange the definitions of $g$-twisted modules
		and $g^{-1}$-twisted modules.
	}\er

	\br{r7}{\em  If $g=1$, these notions reduce to the usual notions of weak,
		admissible and ordinary $V$-modules (cf. \cite{dlm1}).
	}\er
	
	The following definition is due to \cite{dlm}.
	\bd{definition 6}{\em
		Let $g$ be an automorphism of a vertex operator algebra $V$ of finite
		order. We say that $V$ is $\emph{g-rational}$ if the admissible $g$-twisted
		module category is semisimple.
		In particular, a vertex operator algebra $V$ is called $\emph{rational}$ if it is
		$1$-rational.
	}\ed
	
	The following results are well known (see \cite{dlm}).
	
	\bt{thm 8}
	Suppose that $V$ is a $g$-rational vertex operator algebra. Then the following statements hold:

	\begin{enumerate}[
		label=\textup{(\arabic*)},
		leftmargin=*,
		itemsep=0pt,
		topsep=0pt,
		parsep=0pt
		]
		\item Every irreducible admissible $g$-twisted $V$-module is an ordinary
		$g$-twisted $V$-module. Moreover, for any irreducible admissible
		$g$-twisted $V$-module $M$, there exists $\lambda\in\mathbb C$ such that
		\begin{eqnarray*}
			M=\bigoplus_{n\in\frac1T\mathbb N}M_{\lambda+n},
		\end{eqnarray*}
		where $M_{\lambda}\neq0$. The number $\lambda$ is called the conformal
		weight of $M$.
		
		\item There are only finitely many irreducible admissible
		$g$-twisted $V$-modules up to isomorphism.
	\end{enumerate}
	\et

	We next recall the restricted dual. Let $M=\bigoplus_{n\in\frac{1}{T}\mathbb{N}}M(n)$
	be an admissible $g$-twisted $V$-module. Define the $\emph{restricted dual}$ of $M$ as follows:
	\begin{eqnarray*}
		M'=\bigoplus_{n\in\frac{1}{T}\mathbb{N}}M(n)^*,
	\end{eqnarray*}
	where $M(n)^*=\operatorname{Hom}_{\mathbb{C}}(M(n),\mathbb{C})$.
	For $v\in V$, define a vertex operator $Y_{M'}(v,z)$ on $M'$ by
	\begin{eqnarray*}
		\langle Y_{M'}(v,z)f,u\rangle
		=
		\langle f,
		Y_M(e^{zL(1)}(-z^{-2})^{L(0)}v,z^{-1})u
		\rangle ,
	\end{eqnarray*}
	where $\langle f,w\rangle=f(w)$
	denotes the natural pairing between $M'$ and $M$.  It is known from \cite{fhl}, \cite{xxp} and \cite{xf} that
	$(M',Y_{M'})$ naturally carries the structure of an admissible
	$g^{-1}$-twisted $V$-module. More generally, if instead $M=\bigoplus_{\lambda\in\mathbb{C}}M_\lambda$
	is a $g$-twisted $V$-module, we can also define the restricted dual of $M$ in the same way. In this case, $M'$ becomes a $g^{-1}$-twisted $V$-module. Furthermore, $M$ is irreducible if and only if $M'$ is irreducible.
	
	In particular, if $(M,Y_M)$ is a $V$-module, its restricted dual $(M',Y_{M'})$ is also a $V$-module. We call a $V$-module $M$ $self$-$dual$ if there exists a $V$-module isomorphism between $M$ and
	$M'$. Similarly, a vertex operator algebra $V$ is called $self$-$dual$ if
	$V$ and its restricted dual $V'$ are isomorphic as $V$-modules.
	
	We next recall the finiteness condition known as $C_2$-cofiniteness (cf. \cite{z}).
	\bd{definition 9}{\em
		A vertex operator algebra $V$ is said to be $C_2$-$\emph{cofinite}$ if
		\begin{eqnarray*}
			\dim V/C_2(V)<\infty,
		\end{eqnarray*}
		where $C_2(V)=\operatorname{span}\{u_{-2}v\mid u,v\in V\}$.
	}\ed
	
	\bd{definition 10}{\em
		A vertex operator algebra $V$ is said to be of $\emph{CFT type}$ if $V=\bigoplus_{n\geq0}V_{(n)}$ and $V_{(0)}=\mathbb{C}\mathbf{1}.$
	}\ed
	
	Let $V$ be a vertex operator algebra. For any subgroup $G$ of $\operatorname{Aut}(V)$, the
	set of $G$-fixed points
	\begin{eqnarray*}
		V^G=\{v\in V\mid g(v)=v\ \text{for all }g\in G\}
	\end{eqnarray*}
	forms a vertex operator subalgebra. The following result about fixed-point vertex operator algebras was
	proved in \cite{cm} and \cite{m}.
	
	\bt{thm 11}
	Let $V$ be a simple, rational, $C_2$-cofinite and self-dual vertex
	operator algebra of CFT type. Then for any solvable subgroup
	$G\leq\operatorname{Aut}(V)$, the fixed-point subalgebra
	$V^G$ is also a simple, rational, $C_2$-cofinite and self-dual vertex
	operator algebra of CFT type.
	\et
	
	The following result follows from Lemma 4.2 of
	\cite{adjr}.
	
	\bp{pp 12}
	Suppose that $V$ is a simple, rational, $C_2$-cofinite and self-dual
	vertex operator algebra of CFT type. Then $V$ is $g$-rational for any
	finite-order automorphism $g$ of $V$.
	\ep
	
	\subsection{Modular Invariance}
	
	In this subsection, we collect the results on modular invariance  from \cite{z} and \cite{dlm3}.
	Throughout this subsection, we assume that $V$ is a simple, rational, $C_2$-cofinite
	and self-dual vertex operator algebra of CFT type, and let $G$ be a
	finite subgroup of $\operatorname{Aut}(V)$.
	Let $g\in G$ be an automorphism of finite order $T$ and $M$ be a
	$g$-twisted $V$-module. By Theorem \ref{thm 8}, $M$ admits the decomposition
	\begin{eqnarray*}
		M=\bigoplus_{n\in\frac{1}{T}\mathbb{N}}M_{\lambda+n}.
	\end{eqnarray*}
	By Proposition~\ref{pp 12}, $V$ is $g$-rational for $g\in G$.
	
	The automorphism subgroup $G$ naturally acts on the category of weak twisted $V$-modules (cf. \cite{dlm3}). More precisely, let $g,h\in G$ and suppose that
	$(M,Y_M)$ is a weak $g$-twisted $V$-module.
	The vector space $M$ carries a weak
	$h^{-1}gh$-twisted $V$-module structure $(M\circ h,Y_{M\circ h})$, defined by
	\begin{eqnarray*}
		Y_{M\circ h}(v,x)=Y_M(hv,x)\quad \text{for}\,v\in V,
	\end{eqnarray*}
	where $M\circ h \cong M$ as vector space.
	This gives a right action of $G$
	on the set of weak twisted $V$-modules, as well as on their
	isomorphism classes. We use the notation
	\begin{eqnarray*}
		(M,Y_M)\circ h=(M\circ h,Y_{M\circ h})=M\circ h .
	\end{eqnarray*}
	In particular, $h$  acts on the set of
	$g$-twisted $V$-modules if $g$, $h$ commute.
	
	A $g$-twisted $V$-module $M$ is called $\emph{h-stable}$ if $M\circ h\cong M$
	as $g$-twisted $V$-modules. Let $\mathcal{M}(g)$ denote the set of
	isomorphism classes of irreducible $g$-twisted $V$-modules, and define
	\begin{eqnarray*}
		\mathcal{M}(g,h)
		=
		\{M\in\mathcal{M}(g)\mid M\circ h\cong M\}.
	\end{eqnarray*}
	Since $V$ is $g$-rational, Theorem \ref{thm 8} implies that both
	$\mathcal{M}(g)$ and $\mathcal{M}(g,h)$ are finite sets.

	For an irreducible $g$-twisted $V$-module $M$, define $G_M=\{h\in G\mid M\circ h\cong M\}$,
	which is a subgroup of $G$. By Schur's Lemma, there exists a
	projective representation $\varphi$ of $G_M$ on $M$ such that
	\begin{eqnarray}
		\varphi(h)Y_M(u,x)\varphi(h)^{-1}=Y_M(hu,x) \quad \text{for}\, h\in G_M.\label{13}
	\end{eqnarray}
	If $h=1$, we choose $\varphi(1)=1$. Moreover, the element $g$ lies in the
	subgroup $G_M$. In fact, for an admissible $g$-twisted $V$-module $M$,
	the automorphism $g$ acts on each homogeneous subspace $M(n)$ by
	\begin{eqnarray*}
		g|_{M(n)}=e^{2\pi \sqrt{-1} n},
	\end{eqnarray*}
	where $n\in\frac{1}{T}\mathbb{Z}$ (cf. \cite{dxy}). In particular, any admissible $g$-twisted module is $g$-stable.
	
	Set $o(v)=v_{\operatorname{wt}v-1}$ for homogeneous $v\in V$. Let $\mathbb{H}$ denote the complex upper half-plane. We set $q=e^{2\pi\sqrt{-1}\tau}$,
	where $\tau\in\mathbb{H}$. Let $P(G)$ denote the set of ordered commuting pairs
	in $G$. For $(g,h)\in P(G)$ and an irreducible
	$g$-twisted $V$-module $M\in\mathcal{M}(g,h)$, define the trace function by
	\begin{eqnarray}
		Z_M(v,(g,h),\tau)
		=
		\operatorname{tr}_{M}
		o(v)\varphi(h)q^{L(0)-c/24}
		=
		q^{\lambda-c/24}
		\sum_{n\in\frac{1}{T}\mathbb{N}}
		\operatorname{tr}_{M_{\lambda+n}}
		o(v)\varphi(h)q^n,\label{12}
	\end{eqnarray}
	where $v$ is a homogeneous element of $V$. The trace function $Z_M(v,(g,h),\tau)$ is holomorphic on $\mathbb{H}$ (see
	\cite{z}, \cite{dlm3}). For simplicity, we write $Z_M(v,\tau)=Z_M(v,(g,1),\tau)$.

	There is another vertex operator algebra $(V,Y[\, \cdot,\,x \,],\mathbf{1},
	\widetilde{\omega})$ associated to $V$ in \cite{z}. Here $\widetilde{\omega}=\omega-\frac{c}{24}\mathbf{1}$ and
	\begin{eqnarray*}
		Y[v,x]
		=
		Y(v,e^x-1)e^{x\cdot\operatorname{wt}v}
		=
		\sum_{n\in\mathbb{Z}}v[n]x^{-n-1}
	\end{eqnarray*}
	for homogeneous $v\in V$. We also write
	\begin{eqnarray*}
		Y[\widetilde{\omega},x]
		=
		\sum_{n\in\mathbb{Z}}L[n]x^{-n-2}.
	\end{eqnarray*}
	If $v\in V$ is homogeneous with respect to $L[0]$, we denote its
	weight by $\operatorname{wt}[v]$.

	For
	$(g,h)\in P(G)$ and $M\in\mathcal{M}(g,h)$, the trace function
	$Z_M(v,(g,h),\tau)$ is a function on $V\times\mathbb{H}$. Let $W$ be
	the vector space spanned by all such functions, namely,
	\begin{eqnarray*}
		W=\operatorname{span}
		\{Z_M(v,(g,h),\tau)\mid (g,h)\in P(G),\ M\in\mathcal{M}(g,h)\}.
	\end{eqnarray*}
	It follows from \cite{dlm3} that $\dim W=
	\sum_{(g,h)\in P(G)}
	|\mathcal{M}(g,h)|.$

	The modular group $\Gamma=SL(2,\mathbb Z)$ is the group of
	$2\times2$ integral matrices with determinant $1$. We next introduce the action of the modular group $\Gamma=SL(2,\mathbb Z)$ on $W$, such that
	\begin{eqnarray*}
		Z_M|_\gamma(v,(g,h),\tau)
		=
		(c\tau+d)^{-\operatorname{wt}[v]}
		Z_M(v,(g,h),\gamma\tau),
	\end{eqnarray*}
	where
	\begin{eqnarray*}
		\gamma:\tau\mapsto\frac{a\tau+b}{c\tau+d},
		\quad
		\gamma=
		\begin{pmatrix}
			a&b\\
			c&d
		\end{pmatrix}
		\in \Gamma=SL(2,\mathbb Z).
	\end{eqnarray*}
	Moreover, the modular group acts on the right of $P(G)$ by
	\begin{eqnarray*}
		(g,h)\gamma=(g^ah^c,g^bh^d)\quad \text{for}\,\gamma=
		\begin{pmatrix}
			a&b\\
			c&d
		\end{pmatrix}\in \Gamma.
	\end{eqnarray*}
	
	We recall the following result, which was established in \cite{z} and \cite{dlm3}.
	
	\bt{thm 13}
	Let $V$ be a simple, rational, $C_2$-cofinite
	and self-dual vertex operator algebra of CFT type and $G$ a
	finite subgroup of $\operatorname{Aut}(V)$. We assume that $$W=\operatorname{span}
	\{Z_M(v,(g,h),\tau)\mid (g,h)\in P(G),\ M\in\mathcal{M}(g,h)\}.$$ Then there
	exists a representation
	\begin{eqnarray*}
		\rho:\Gamma\longrightarrow GL(W)
	\end{eqnarray*}
	such that for $(g,h)\in P(G)$,
	$\gamma=
	\begin{pmatrix}
		a&b\\
		c&d
	\end{pmatrix}
	\in\Gamma,$
	and  $M\in\mathcal{M}(g,h)$, we have
	\begin{eqnarray*}
		Z_M|_{\gamma}(v,(g,h),\tau)
		=
		\sum_{N\in\mathcal{M}(g^ah^c,g^bh^d)}
		\gamma_{M,N}
		Z_N(v,(g,h)\gamma,\tau),
	\end{eqnarray*}
	where $\rho(\gamma)=(\gamma_{M,N}).$
	That is,
	\begin{eqnarray*}
		Z_M(v,(g,h),\gamma\tau)
		=
		(c\tau+d)^{\operatorname{wt}[v]}
		\sum_{N\in\mathcal{M}(g^ah^c,g^bh^d)}
		\gamma_{M,N}
		Z_N(v,(g,h)\gamma,\tau).
	\end{eqnarray*}
	\et

	Since the modular group $\Gamma=SL(2,\mathbb Z)$ is generated by the two
	elements
	\begin{eqnarray*}
		S=
		\begin{pmatrix}
			0&-1\\
			1&0
		\end{pmatrix}
		\,\text{and}\,\,
		T=
		\begin{pmatrix}
			1&1\\
			0&1
		\end{pmatrix},
	\end{eqnarray*}
	the representation $\rho$ of $\Gamma$ on the space $W$ is
	completely determined by the matrices $\rho(S)$ and $\rho(T)$. The
	matrix $\rho(S)$ is called the $S$-$matrix$ of the orbifold theory.
	There are two special cases of the $S$-transformation:
	\begin{eqnarray*}
		Z_M\left(v,-\frac1\tau\right)
		=
		\tau^{\operatorname{wt}[v]}
		\sum_{N\in\mathcal{M}(1,g^{-1})}
		S_{M,N}Z_N(v,(1,g^{-1}),\tau)
	\end{eqnarray*}
	for $M\in\mathcal{M}(g)=\mathcal{M}(g,1)$ and
	\begin{eqnarray*}
		Z_N\left(v,(1,g),-\frac1\tau\right)
		=
		\tau^{\operatorname{wt}[v]}
		\sum_{M\in\mathcal{M}(g)}
		S_{N,M}Z_M(v,\tau)
	\end{eqnarray*}
	for $N\in\mathcal{M}(1,g)$. Here the matrix $S=(S_{M,N})_{M,N\in\mathcal{M}(1)}$ is called the $S$-$matrix$ of $V$ and is independent of the choice of vector $v$.
	
	\subsection{Fusion Rules for Twisted Modules and Verlinde Formula}
	In this subsection, we recall several results on the fusion rules of
	twisted modules from \cite{dlxy}. We also recall the Verlinde formula, as established in \cite{h1} and \cite{dln}. Throughout this subsection, we assume
	that $V$ is a simple, rational, $C_2$-cofinite, self-dual vertex
	operator algebra of CFT type, and that $G$ is a finite abelian subgroup
	of $\operatorname{Aut}(V)$. It follows from Theorem \ref{thm 11} that
	the fixed-point subalgebra $V^G$ is also
	simple, rational, $C_2$-cofinite, self-dual vertex
	operator algebra of CFT type.

	We first recall some basic notions about modular tensor categories
	from \cite{bk}, \cite{ko} and \cite{egno}.
	\bd{definition 9}{\em Let $\mathcal{C}$ be a $\emph{modular tensor category}$ with
		tensor product functor $\otimes$, unit object $\mathbf{1}_{\mathcal C}$,
		unit isomorphisms $l$ and $r$, associativity isomorphism $a$, and
		braiding isomorphism $c$. An $\emph{algebra}$ in $\mathcal{C}$ consists of an object $A\in\mathcal C$
		together with morphisms
		\begin{eqnarray*}
			\mu:A\otimes A\to A,
			\quad
			\iota_A:\mathbf{1}_{\mathcal C}\hookrightarrow A,
		\end{eqnarray*}
		which satisfy the following associativity and unit conditions:
		\begin{enumerate}[
			label=\textup{(\arabic*)},
			leftmargin=*,
			itemsep=0pt,
			topsep=0pt,
			parsep=0pt
			]
			\item \textit{Associativity.}
			The compositions
			\begin{eqnarray*}
				\mu\circ(\mu\otimes \operatorname{Id}_A)\circ a,
				\quad
				\mu\circ(\operatorname{Id}_A\otimes\mu):
				A\otimes(A\otimes A)\to A
			\end{eqnarray*}
			are equal.
			
			\item \textit{Unit.}
			The composition
			\begin{eqnarray*}
				\mu\circ(\iota_A\otimes\operatorname{Id}_A):
				A=\mathbf{1}_{\mathcal C}\otimes A\to A
			\end{eqnarray*}
			is equal to $l_A$; and the composition
			\begin{eqnarray*}
				\mu\circ(\operatorname{Id}_A\otimes\iota_A):
				A=A\otimes\mathbf{1}_{\mathcal C}\to A
			\end{eqnarray*}
			is equal to $r_A$.
		\end{enumerate}
	}\ed

	We denote the algebra object $(A,\mu,\iota_A)$ simply by
	$A$. The algebra
	object $A$ is called $\emph{commutative}$ if
	\begin{eqnarray*}
		\mu\circ c_{A,A}=\mu,
	\end{eqnarray*}
	where $c_{A,A}:A\otimes A\rightarrow A\otimes A$ is the braiding
	isomorphism.
	
	\bd{definition 9}{\em
		Let $A$ be an algebra object in $\mathcal{C}$. The category
		$\operatorname{Rep}(A)$  is defined as follows. An object
		of $\operatorname{Rep}(A)$ is a pair $(M,\mu_M)$, where
		$M\in\mathcal{C}$ and
		\begin{eqnarray*}
			\mu_M:A\otimes M\to M
		\end{eqnarray*}
		is a morphism in $\mathcal{C}$ satisfying the following conditions:
		\begin{eqnarray*}
			&&\mu_M\circ(\mu\otimes\operatorname{Id}_M)\circ a
			=
			\mu_M\circ(\operatorname{Id}_A\otimes\mu_M):
			A\otimes(A\otimes M)\to M,\\
			&&\mu_M\circ(\iota_A\otimes\operatorname{Id}_M)
			=
			l_M:
			\mathbf{1}_{\mathcal C}\otimes M\to M .
		\end{eqnarray*}
		For two objects $(M_1,\mu_{M_1})$ and $(M_2,\mu_{M_2})$ in
		$\operatorname{Rep}(A)$, the morphisms are given by
		\begin{eqnarray*}
			\operatorname{Hom}_{\operatorname{Rep}(A)}
			((M_1,\mu_{M_1}),(M_2,\mu_{M_2}))
			=
			\{\varphi\in\operatorname{Hom}_{\mathcal C}(M_1,M_2)
			\mid
			\mu_{M_2}\circ(\operatorname{Id}_A\otimes\varphi)
			=
			\varphi\circ\mu_{M_1}\}.
		\end{eqnarray*}
	}\ed
	
	The following result was proved in \cite{ko}.
	
	\bt{thm 17}
	The category $\operatorname{Rep}(A)$ is a tensor category with unit object
	$A$.
	\et
	
	We also need the following results from Theorem~1.6 of \cite{ko}.
	
	\bt{thm 18}
	Define a functor $\mathcal{F}:\mathcal{C}\to\operatorname{Rep}(A)$
	by
	\begin{eqnarray*}
		\mathcal{F}(M)=A\otimes M,
		\quad
		\mu_{\mathcal{F}(M)}=\mu\otimes\operatorname{Id}_M \quad \text{for} \,\, M\in \mathcal{C}.
	\end{eqnarray*}
	Then the following statements hold:
	
	\begin{enumerate}[
		label=\textup{(\arabic*)},
		leftmargin=*,
		itemsep=0pt,
		topsep=0pt,
		parsep=0pt
		]
		\item The functor $\mathcal{F}$ is exact.
		
		\item The functor $\mathcal{F}$ is a tensor functor.
		
		\item Let $G:\operatorname{Rep}(A)\to\mathcal{C}$
		be the forgetful functor defined by $G(M,\mu_M)=M.$
		Then, for every $M \in \mathcal{C}$ and $X\in\operatorname{Rep}(A)$, there exists a natural
		isomorphism
		\begin{eqnarray*}
			\operatorname{Hom}_{\operatorname{Rep}(A)}
			(\mathcal{F}(M),X)
			\cong
			\operatorname{Hom}_{\mathcal C}(M,G(X)).
		\end{eqnarray*}
	\end{enumerate}
	\et

	Let $\mathcal{C}_V$ denote the category of ordinary $V$-modules and
	$\mathcal{C}_{V^G}$ the category of ordinary $V^G$-modules. It follows
	from \cite{h2} that both $\mathcal{C}_V$ and $\mathcal{C}_{V^G}$ are
	modular tensor categories.
	Moreover, by \cite{ko}, \cite{hkl} and \cite{ckm}, the vertex operator algebra $V$ can be
	regarded as a commutative associative algebra object in
	$\mathcal{C}_{V^G}$. More precisely, we have the decomposition $V=\bigoplus_{\chi\in\operatorname{Irr}(G)}V^\chi,$
	where $\operatorname{Irr}(G)$ denotes the set of irreducible characters
	of $G$ and each $V^\chi$ is an irreducible $V^G$-module (cf.
	\cite{djx}). Hence, $V$ is a direct sum of simple objects in the modular
	tensor category $\mathcal{C}_{V^G}$.
	
	Viewing $V$ as an algebra object in $\mathcal{C}_{V^G}$, we define
	the category $\operatorname{Rep}(V)$ as in Definition~\ref{definition 9}. By \cite{ko} (see also \cite{ckm}), there
	exists a categorical tensor product functor
	\begin{eqnarray*}
		\boxtimes_V:\operatorname{Rep}(V)\times\operatorname{Rep}(V)
		\to\operatorname{Rep}(V),
	\end{eqnarray*}
	which is associative. We denote the tensor product of two objects
	$M,N\in\operatorname{Rep}(V)$ by $M\boxtimes_V N.$
	Furthermore, $\operatorname{Rep}(V)$ is a fusion category. In
	particular, it is a semisimple category with finitely many
	inequivalent simple objects.
	
	Before defining the fusion rules among twisted modules, we first recall
	the notion of intertwining operators (cf. \cite{xxp}).
	Let $g_1,g_2,g_3$ be mutually commuting automorphisms of $V$ with
	finite orders $T_1,T_2,T_3$, respectively. The vector space $V$ has a decomposition into common eigenspaces for the actions of $g_1$ and $g_2$:
	\begin{eqnarray*}
		V=\bigoplus_{\substack{0\leq j_1\leq T_1-1,\\0\leq j_2\leq T_2-1}}
		V^{(j_1,j_2)},
	\end{eqnarray*}
	where $V^{(j_1,j_2)}
	=
	\{v\in V\mid g_s(v)=e^{\frac{-2\pi \sqrt{-1}j_s}{T_s}}v,\ s=1,2\},0\leq j_1\leq T_1-1,0\leq j_2\leq T_2-1.$
	
	The following definition is due to \cite{xxp}.
	\bd{definition 19}{\em
		Let $(M_i,Y_{M_i})$ be a weak $g_i$-twisted $V$-module for
		$i=1,2,3$. An $\emph{intertwining operator of type}$
		$\binom{M_3}{M_1\,M_2}$ is a linear map
		\begin{eqnarray*}
			&&\mathcal{Y}(\cdot,x):M_1
			\to
			\operatorname{Hom}(M_2,M_3)\{x\},\\
			&&\hspace{1.6cm}w \mapsto \mathcal{Y}(w,x)
			=
			\sum_{n\in\mathbb C}w_nx^{-n-1},
		\end{eqnarray*}
		which satisfies the following conditions:
		\begin{enumerate}[label=\textup{(\arabic*)}, leftmargin=*]
			
			\item \textit{Lower truncation.}
			For  $w^1\in M_1$, $w^2\in M_2$ and any fixed $c\in\mathbb C$,
			\begin{eqnarray*}
				w^1_{c+n}w^2=0\quad \text{for}\,n\in\mathbb {Q}\,\,\text{sufficiently large};
			\end{eqnarray*}
			
			\item \textit{The Jacobi identity.}
			For any $u\in V^{(j_1,j_2)}$ with $j_1,j_2\in \mathbb{Z}$ and $w^1\in M_1$, the following identity
			holds:
			\begin{eqnarray}
				&&x_0^{-1}\left(\frac{x_1-x_2}{x_0}\right)^{\frac{j_1}{T_1}}
				\delta\left(\frac{x_1-x_2}{x_0}\right)
				Y_{M_3}(u,x_1)\mathcal{Y}(w^1,x_2)\notag
				\\
				&&-
				x_0^{-1}\left(\frac{x_2-x_1}{-x_0}\right)^{\frac{j_1}{T_1}}
				\delta\left(\frac{x_2-x_1}{-x_0}\right)
				\mathcal{Y}(w^1,x_2)Y_{M_2}(u,x_1)\notag
				\\
				&&\hspace{-0.3cm}=
				x_2^{-1}
				\left(\frac{x_1-x_0}{x_2}\right)^{-\frac{j_2}{T_2}}
				\delta\left(\frac{x_1-x_0}{x_2}\right)
				\mathcal{Y}(Y_{M_1}(u,x_0)w^1,x_2);\label{15}
			\end{eqnarray}
			\item \textit{The $L(-1)$-derivative property.}
			For any $w^1\in M_1$,
			\begin{eqnarray}
				\frac{d}{dx}\mathcal{Y}(w^1,x)
				=
				\mathcal{Y}(L(-1)w^1,x).\label{16}
			\end{eqnarray}
			
		\end{enumerate}
		The vector space of all intertwining operators of type
		$\binom{M_3}{M_1\,M_2}$ is denoted by
		\begin{eqnarray*}
			I_V\binom{M_3}{M_1\,M_2}.
		\end{eqnarray*}
	}\ed

	\bd{definition 20}{\em  Let $(M_i,Y_{M_i})$ be a weak $g_i$-twisted $V$-module for
		$i=1,2,3$. The $\emph{fusion rule}$ is given by
		\begin{eqnarray*}
			N_{M_1M_2}^{M_3}
			=
			\dim I_V
			\binom{M_3}{M_1\,M_2}.
		\end{eqnarray*}
	}\ed

	The following formula was given in \cite{xxp}.
	\bp{pp wa} For any $u\in V^{(j_1,j_2)}$ with $j_1,j_2\in \mathbb{Z}$ and $w^1\in M_1, w^2\in M_2$, there exists a non-negative integer $m(u,w^2)$ depending only on $u$ and $w^2$, such that
	\begin{eqnarray*}
		&(x_2+x_0)^{m(u,w^2)+\frac{j_2}{T_2}}
		\mathcal{Y}(Y_{M^1}(u,x_0)w^1,x_2)w^2
		\nonumber\\
		&\hspace{3cm}=
		(x_0+x_2)^{m(u,w^2)+\frac{j_2}{T_2}}
		Y_{M^3}(u,x_0+x_2)\mathcal{Y}(w^1,x_2)w^2 .
	\end{eqnarray*}
	\ep

	The following result follows from Proposition 3.46 of \cite{ckm} (see also Definition 3.1 of \cite{dlxy}).
	\bp{pp 20}
	$\operatorname{Rep}(V)$ is the subcategory of $\mathcal{C}_{V^G}$
	consisting of every $V^G$-module $W$ together with a
	$V^G$-intertwining operator
	\begin{eqnarray*}
		Y_W(\cdot,x):
		V\to(\operatorname{End}W)\{x\}
	\end{eqnarray*}
	of type $\binom{W}{V\,W}$,
	such that the following conditions are satisfied:

	\begin{enumerate}[
		label=\textup{(\arabic*)},
		leftmargin=*,
		itemsep=0pt,
		topsep=0pt,
		parsep=0pt
		]
		
		\item \textit{Associativity.} For any $u,v\in V$, $w\in W$ and $w'\in W'$, the formal series
		\begin{eqnarray*}
			\langle w',
			Y_W(u,x_1)Y_W(v,x_2)w\rangle
		\end{eqnarray*}
		and
		\begin{eqnarray*}
			\langle w',
			Y_W(Y(u,x_1-x_2)v,x_2)w\rangle
		\end{eqnarray*}
		converge on the domains $|x_1|>|x_2|>0$ and $|x_2|>|x_1-x_2|>0,$
		respectively, to multivalued analytic functions which coincide on their
		common domain.
		
		\item \textit{Unit.} $Y_W(\mathbf{1},x)=\operatorname{Id}_W.$
	\end{enumerate}
	\ep
	
	We remark that the category $\operatorname{Rep}(V)$ contains objects
	more general than ordinary $V$-modules. The following result
	can be found in \cite{k1}, \cite{k2} and \cite{dlxy}.

	\bp{pp 21} Let $V,G$ be as before. Then the following
	statements hold:
	\begin{enumerate}[label=\textup{(\arabic*)}, leftmargin=*,
		itemsep=0pt, topsep=0pt]
		
		\item If $W$ is a $g$-twisted $V$-module for some $g\in G$, then $W$ is
		an object of $\operatorname{Rep}(V)$. Moreover, if $W_i$ is a
		$g_i$-twisted $V$-module with $g_i\in G$ for $i=1,2$, then $W_1$ and
		$W_2$ are equivalent as objects in $\operatorname{Rep}(V)$ if and only
		if $g_1=g_2$ and $W_1\cong W_2$ as twisted $V$-modules.
		
		\item If $W$ is a simple object in $\operatorname{Rep}(V)$, then $W$ is
		an irreducible $g$-twisted $V$-module for some $g\in G$.
	\end{enumerate}
	\ep
	
	As a result, the category $\operatorname{Rep}(V)$ decomposes as $\operatorname{Rep}(V)
	=
	\bigoplus_{g\in G}\mathcal{C}_{V,g},$
	where $\mathcal{C}_{V,g}$ denotes the category of $g$-twisted
	$V$-modules (cf. \cite{dxy2}).
	
	We now recall another notion of fusion rules for twisted modules. Suppose that $M^i$ is a $g_i$-twisted $V$-module with
	$g_i\in G$ for $i=1,2,3$. By Proposition~\ref{pp 21}, the modules $M^1$,
	$M^2$ and $M^3$ can be regarded as objects in $\operatorname{Rep}(V)$.
	Therefore, their tensor product
	\begin{eqnarray*}
		M^1\boxtimes_V M^2
	\end{eqnarray*}
	exists in $\operatorname{Rep}(V)$. As $\operatorname{Rep}(V)$ is a fusion category, it is
	semisimple and the decomposition of $M^1\boxtimes_V M^2$ into simple
	objects is unique. We define fusion rules $N_{M^1,M^2}^{M^3}$ to be the
	multiplicity of the simple object $M^3$ in $M^1\boxtimes_V M^2$. By \cite{dlxy}, this definition of fusion rules is equivalent to the one
	given in Definition~\ref{definition 20}.
	
	\br{remark 23}{\em If there exist weak $g_s$-twisted modules $(M_s,Y_{M_s})$ for
		$s=1,2,3$ such that $N_{M_1,M_2}^{M_3}>0,$
		then $g_3=g_1g_2$
		(cf. \cite{xxp}, \cite{dlxy}). Consequently, we will always assume
		\begin{eqnarray*}
			g_3=g_1g_2 .
		\end{eqnarray*}
	}\er
	
	The following property of fusion rules is given in Proposition 3.9 of \cite{dl} (see also Corollary 5.2 of \cite{h3}).
	\bp{pp 24}
	Let $g_1,g_2$ be commuting finite-order automorphisms of a vertex
	operator algebra $V$ and let $M_i$ be a $g_i$-twisted $V$-module for
	$i=1,2,3$ with $g_3=g_1g_2.$
	Then
	\begin{eqnarray}
		N_{M_1,M_2}^{M_3}
		=
		N_{M_2\circ g_1^{-1},M_1}^{M_3}
		=
		N_{M_2,M_1\circ g_2}^{M_3}.\label{14}
	\end{eqnarray}
	\ep
	
	The following results can be found in \cite{h1} and \cite{dln}.
	
	\bt{thm 25}
	Let $V$ be a simple, rational, $C_2$-cofinite and self-dual vertex
	operator algebra of CFT type, and let
	$M^0=V,M^1,\ldots,M^p$ be all inequivalent irreducible $V$-modules,
	where $p$ is a non-negative integer. Then we have
	
	\begin{enumerate}[label=\textup{(\arabic*)}, leftmargin=*,
		itemsep=0pt, topsep=0pt]
		
		\item The  $S$-matrix is symmetric and unitary.
		
		\item For $0\leq i,j,k\leq p$, the Verlinde formula holds:
		\begin{eqnarray*}
			N_{M^i,M^j}^{M^k}
			=
			\sum_{0\leq r\leq p}
			\frac{
				S_{M^i,M^r}
				S_{M^j,M^r}
				\overline{S_{M^k,M^r}}
			}{
				S_{V,M^r}
			}.
		\end{eqnarray*}
		
	\end{enumerate}
	\et

	\section{$\textit{S}$-Matrices}

	In this section, we assume that $V$ is a simple, rational, $C_2$-cofinite
	and self-dual vertex operator algebra of CFT type and that $\sigma$ is
	an automorphism of $V$ with finite order $n$. Let $G=\langle\sigma\rangle$
	be the cyclic subgroup of $\operatorname{Aut}(V)$ generated by
	$\sigma$. By Theorem~\ref{thm 11}, the fixed-point subalgebra $V^G$ is also a
	simple, rational, $C_2$-cofinite and self-dual vertex operator algebra of
	CFT type. Moreover, Proposition~\ref{pp 12} implies that $V$ is $g$-rational
	for any $g\in G$. The purpose of this section is to investigate the relationship between the
	$S$-matrices of $V$ and its fixed-point subalgebra $V^G$.
	
	We first recall some results about irreducible $V^G$-modules from
	\cite{drx1}.
	Let $M$ be an irreducible $g$-twisted $V$-module. Recall that
	\begin{eqnarray*}
		G_M=\{h\in G\mid M\circ h\cong M\}
	\end{eqnarray*}
	is the subgroup of $G$. Since $G$ is a cyclic group of order $n$, the subgroup $G_M$ is also
	cyclic, and its order divides $n$. Recall from Section~2.2 that there exists a projective representation $\varphi$ of $G_M$ on $M$. Following the proof of Proposition~5.3 of \cite{dm}, the corresponding $2$-cocycle is trivialized and henceforth we choose $\varphi$ to be a representation of $G_M$.

	Suppose that $G_M=H=\langle h\rangle$ is a cyclic subgroup of order
	$m$, where $m$ divides $n$. For $0\leq s\leq m-1$, let $\Lambda_s$ be the character of $G_M$ defined by
	\begin{eqnarray*}
		\Lambda_s(h)=e^{\frac{-2\pi\sqrt{-1}s}{m}}.
	\end{eqnarray*}
	Then $\{\Lambda_s\mid 0\leq s\leq m-1\}$
	is precisely the set of all irreducible characters of $G_M$, which is
	denoted by $\operatorname{Irr}(G_M)$. Therefore,
	the representation $\varphi$ induces the following decomposition of $M$:
	\begin{eqnarray}
		M=\bigoplus_{0\leq s\leq m-1}M_{\Lambda_s}, \label{3.2}
	\end{eqnarray}
	where $M_{\Lambda_s}
	=
	\{w\in M\mid \varphi(h)w=\Lambda_s(h)w\}.$

	The following results about  $V^G$-modules were established
	in \cite{dy} and \cite{mt} (see also \cite{drx1}, \cite{drx2}).

	\bt{thm 3.1}
	Let $M$ be an irreducible $g_1$-twisted $V$-module and  $N$  an irreducible $g_2$-twisted $V$-module. Then
	
	\begin{enumerate}[label=\textup{(\arabic*)}, leftmargin=*,
		itemsep=0pt, topsep=0pt]
		
		\item $M_{\chi}$ is nonzero for any $\chi \in \operatorname{Irr}(G_M)$.
		
		\item Each $M_{\chi}$ is an irreducible $V^G$-module.
		
		\item $M_{\chi}$ and $M_{\gamma}$ are equivalent
		$V^G$-modules if and only if $\chi=\gamma$.
		
		\item If $M$ and $N$ are not in the same $G$-orbit, namely, $M$ and $N\circ g$ are inequivalent for all $g\in G$, then
		$M_{\chi}$ and $N_{\gamma}$ are inequivalent
		$V^G$-modules for any $\chi \in \operatorname{Irr}(G_M)$ and $\gamma \in \operatorname{Irr}(G_N)$.
	\end{enumerate}
	\et

	Following \cite{drx1} and \cite{drx2} , the group $G$ acts on the set of all
	twisted modules
	\begin{eqnarray*}
		\mathcal{S}
		=
		\bigcup_{g\in G}\mathcal{M}(g)
	\end{eqnarray*}
	by $M\mapsto M\circ h \quad \text{for}\,h\in G.$
	For any $h\in G$, $M\circ h$ is isomorphic to $M$ as a
	$V^G$-module. Furthermore, the orbit of $M$ has cardinality $|M\circ G|=[G:G_M].$
	
	Note that if $M$ is an irreducible $g$-twisted $V$-module with
	$G_M=\{1\}$, then $g=1$. To describe the stabilizer group $G_M$ of
	$M$, we introduce the following notation:
	\begin{eqnarray*}
		\mathcal{M}_H(g)
		=
		\{M\in \mathcal{M}(g)\mid G_M=H\},
	\end{eqnarray*}
	where $\langle g\rangle\leq H\leq G$.
	Therefore, we have a disjoint union $\mathcal{M}(g)=\bigcup_{g \in H\leq G}\mathcal{M}_H(g)$.
	Clearly, $G_M=G_{M\circ h}$ for any $h\in G$. Thus, $\mathcal{M}_H(g)$ admits a decomposition into $G$-orbits:
	\begin{eqnarray*}
		\mathcal{M}_H(g)
		=
		\bigcup_{j\in J_{g,H}} \mathcal{O}_{g,H,j}.
	\end{eqnarray*}
	For each orbit $\mathcal{O}_{g,H,j}$, we fix an orbit representative $M^{g,H,j}$. Then we have
	\begin{eqnarray}
		\mathcal{O}_{g,H,j}=\{M^{g,H,j}\circ f\mid f\in G \}. \label{3.1}
	\end{eqnarray}
	Here, $\mathcal{O}_{g,H,j}$ denotes the equivalence class under the
	action of $G$. In fact, the elements $f$ in (\ref{3.1}) can
	be chosen as representatives of the cosets in $G/H$.
	
	We decompose $M^{g,H,j}$ according to (\ref{3.2}) as follows:
	\begin{eqnarray}
		M^{g,H,j}
		=
		\bigoplus_{\chi\in\operatorname{Irr}(H)}
		M^{g,H,j}_{\chi},\label{3.3}
	\end{eqnarray}
	where $M^{g,H,j}_{\chi}
	=
	\{w\in M^{g,H,j}\mid \varphi(h)w=\chi(h)w\}$ and $\langle h\rangle=H$.
	
	Now, for each subgroup $H$ of $G$, we define
	\begin{eqnarray*}
		\mathcal{S}_H
		=
		\bigcup_{\substack{\langle g\rangle\leq H}}
		\mathcal{M}_H(g).
	\end{eqnarray*}
	Then the set of all twisted and untwisted modules can be decomposed as
	\begin{eqnarray*}
		\mathcal{S}=\bigcup_{g\in G}\mathcal{M}(g)
		=
		\bigcup_{\substack{H\leq G}}
		\mathcal{S}_H.
	\end{eqnarray*}

	The following result gives a classification of irreducible
	$V^G$-modules (cf. Theorem 3.3 of \cite{drx1}).
	
	\bt{thm 3.2} Let $V$ be a simple, rational, $C_2$-cofinite, and self-dual
	vertex operator algebra of CFT type. Let $\sigma\in\operatorname{Aut}(V)$ be of finite order $n$ and set
	$G=\langle\sigma\rangle$.
	Then the set
	\begin{eqnarray*}
		\{M_{\chi}^{g,H,j}\mid g\in G,\
		\langle g\rangle\leq H\leq G,\
		j\in J_{g,H},\
		\chi\in\operatorname{Irr}(H)\}
	\end{eqnarray*}
	forms a complete list of inequivalent irreducible $V^G$-modules.
	\et

	We next recall some facts about the $S$-matrix of $V^G$
	from \cite{drx2}.
	For an irreducible $g$-twisted $V$-module $M$, set $\mathcal{O}_M=\{M\circ h\mid h\in G\}$,
	which is the orbit of $M$ under the action of $G$.
	Let $g_1,g_2\in G$, and let $M^i$ be an irreducible
	$g_i$-twisted $V$-module for $i=1,2$. Define
	$C_{M^1,M^2}$ to be an least subset of $G$ such that
	\begin{eqnarray*}
		\{M^2\circ k\mid k\in C_{M^1,M^2}\}
		=
		\mathcal{O}_{M^2}
		\cap
		\left(
		\bigcup_{h\in G_{M^1}}
		\mathcal{M}(h,g_1^{-1})
		\right).
	\end{eqnarray*}

	The following lemma follows immediately from the above definition.
	\bl{lemma 3.3}
	Let $M^1$ and $M^2$ be irreducible $g_1$-twisted and
	$g_2$-twisted $V$-modules, respectively. Then
	\begin{eqnarray*}
		C_{M^1,M^2}
		=
		\begin{cases}
			\text{a set of representatives for the cosets of }G/G_{M^2},
			& g_2\in G_{M^1},\ g_1\in G_{M^2},\\
			\emptyset,
			& \text{otherwise}.
		\end{cases}
	\end{eqnarray*}
	\el

	The following result relates the $S$-matrix of $V^G$ to the
	$S$-matrix associated with the orbifold theory of $V$. It follows from
	Lemma~\ref{lemma 3.3} and Corollary 5.4 of \cite{drx2}.
	
	\bt{thm 3.4}
	Let \(V\) be a simple, rational, \(C_2\)-cofinite, self-dual vertex operator algebra of CFT type.
	Let $\sigma\in\operatorname{Aut}(V)$ be of finite order $n$ and set
	$G=\langle\sigma\rangle$. Suppose that $M^i$ is an irreducible
	$g_i$-twisted $V$-module for $i=1,2$. Then we have
	\[
	S_{M^1_{\chi}, M^2_{\gamma}} =
	\begin{cases}
		\dfrac{\overline{\chi(g_2)}\gamma(g_1^{-1})}{|G_{M^1}|}\sum_{k \in C_{M^1,M^2}}S_{M^1,M^2\circ k},
		& g_2\in G_{M^1},\ g_1\in G_{M^2}, \\[1.2ex]
		0, & \text{otherwise},
	\end{cases}
	\]
	where \(\chi \in \operatorname{Irr}(G_{M^1})\) and \(\gamma \in \operatorname{Irr}(G_{M^2})\).
	\et

	\section{Fusion Rules}
	In this section, we assume that $V$ is a simple, rational, $C_2$-cofinite
	and self-dual vertex operator algebra of CFT type and that $\sigma$ is
	an automorphism of $V$ with finite order $n$. Let $G=\langle\sigma\rangle$
	be the cyclic subgroup of $\operatorname{Aut}(V)$ generated by
	$\sigma$. We will establish a relationship between the fusion
	rules of $V$ and  $V^G$.
	
	We first recall the following result from \cite{dl}.	
	\bp{pp 4.1}
	Let $g_1,g_2,g_3\in G$, and $M^i$ be an irreducible
	$g_i$-twisted $V$-module for $i=1,2,3$. Then for any $h\in G$,
	we have
	\begin{eqnarray*}
		N^{M^3\circ h}_{M^1\circ h,M^2\circ h}
		=
		N^{M^3}_{M^1,M^2}.
	\end{eqnarray*}
	\ep

	We will use the following result.
	\bp{pp 4.2}
	Let $\mathcal{F}:\mathcal{C}_{V^G}\to \operatorname{Rep}(V)$ be the functor defined in Theorem $\mathrm{\ref{thm 18}}$. Given $g\in G,\
	\langle g\rangle\leq H\leq G$ and $
	j\in J_{g,H}$, let
	$M^{g,H,j}$ be a representative of the orbit
	$\mathcal{O}_{g,H,j}$.
	For $\chi\in\operatorname{Irr}(H)$, we have
	\begin{eqnarray*}
		\mathcal{F}(M^{g,H,j}_{\chi})
		=
		\bigoplus_{W\in \mathcal{O}_{g,H,j}}W.
	\end{eqnarray*}
	\ep

	\begin{proof}
		For convenience, set $M=M^{g,H,j}$ and $\mathcal{O}_M=\mathcal{O}_{g,H,j}$.
		Let $N$ be a simple object in $\operatorname{Rep}(V)$. By
		Theorem \ref{thm 18} (3), we have
		\begin{eqnarray*}
			\operatorname{Hom}_{\operatorname{Rep}(V)}
			(\mathcal{F}(M_{\chi}),N)
			\cong
			\operatorname{Hom}_{\mathcal{C}_{V^G}}
			(M_{\chi},N).
		\end{eqnarray*}
		If $\operatorname{Hom}_{\mathcal{C}_{V^G}}(M_{\chi},N)\neq0$,
		then Theorem \ref{thm 3.1} implies that $M_{\chi}$ is an irreducible $V^G$-module
		and $N\in\mathcal{O}_M$. Hence
		\begin{eqnarray}
			\dim
			\operatorname{Hom}_{\operatorname{Rep}(V)}
			(\mathcal{F}(M_{\chi}),N)
			=
			\begin{cases}
				1,&N\in\mathcal{O}_M,\\
				0,&N\notin\mathcal{O}_M .
			\end{cases}\label{4.1}
		\end{eqnarray}
		Since $\operatorname{Rep}(V)$ is a semisimple category, we can write
		\begin{eqnarray*}
			\mathcal{F}(M_{\chi})
			=
			\bigoplus_L L^{\bigoplus m_L},
		\end{eqnarray*}
		where $L$ runs over the simple objects of $\operatorname{Rep}(V)$ and
		$m_L$ denotes the multiplicity of $L$. Therefore, we have
		\begin{eqnarray*}
			\operatorname{Hom}_{\operatorname{Rep}(V)}
			(\mathcal{F}(M_{\chi}),N)
			=\bigoplus_L \operatorname{Hom}_{\operatorname{Rep}(V)}(L,N)^{\bigoplus m_L}.
		\end{eqnarray*}
		By Schur's Lemma, we have
		\begin{eqnarray*}
			\operatorname{Hom}_{\operatorname{Rep}(V)}(L,N)
			=
			\begin{cases}
				\mathbb{C},&L\cong N,\\
				0,&L\not\cong N .
			\end{cases}
		\end{eqnarray*}
		Comparing the above decomposition with (\ref{4.1}), we obtain
		\begin{eqnarray*}
			m_L=
			\begin{cases}
				1,&L\in\mathcal{O}_M,\\
				0,&L\notin\mathcal{O}_M .
			\end{cases}
		\end{eqnarray*}
		Therefore,
		\begin{eqnarray*}
			\mathcal{F}(M_{\chi})
			=
			\bigoplus_{W\in\mathcal{O}_M}W .
		\end{eqnarray*}
		
		This completes the proof.\end{proof}

	As an application, the fusion rules of $V$ can be obtained from the
	fusion rules of $V^G$.
	\bt{thm 4.3}
	Let $M^1,M^2,M^3$ be orbit
	representatives of the form $M^{g,H,j}$. For any $\chi_1\in\operatorname{Irr}(G_{M^1}),$ $\chi_2\in\operatorname{Irr}(G_{M^2}),$
	we have
	\begin{eqnarray*}
		\sum_{U\in\mathcal{O}_{M^1}}
		\sum_{W\in\mathcal{O}_{M^2}}
		N^{M^3}_{U,W}
		=
		\sum_{\chi_3\in\operatorname{Irr}(G_{M^3})}
		N^{M^3_{\chi_3}}_{M^1_{\chi_1},M^2_{\chi_2}}.
	\end{eqnarray*}
	\et

	\begin{proof}
		By Proposition \ref{pp 4.2}, we have
		\begin{eqnarray*}
			&&\mathcal{F}(M^1_{\chi_1})
			=
			\bigoplus_{U\in\mathcal{O}_{M^1}}U,\\
			&&\mathcal{F}(M^2_{\chi_2})
			=
			\bigoplus_{W\in\mathcal{O}_{M^2}}W .
		\end{eqnarray*}
		Since $\mathcal{F}$ is a tensor functor by Theorem \ref{thm 18}, we obtain
		\begin{eqnarray}
			\left(
			\bigoplus_{U\in\mathcal{O}_{M^1}}U
			\right)
			\boxtimes_V
			\left(
			\bigoplus_{W\in\mathcal{O}_{M^2}}W
			\right)
			&=&
			\mathcal{F}(M^1_{\chi_1})
			\boxtimes_V
			\mathcal{F}(M^2_{\chi_2})
			\nonumber\notag\\
			&=&
			\mathcal{F}
			\left(
			M^1_{\chi_1}\boxtimes_{V^G}M^2_{\chi_2}
			\right)
			\nonumber\notag\\
			&=&
			\bigoplus_X
			N^{X}_{M^1_{\chi_1},M^2_{\chi_2}}
			\mathcal{F}(X),\label{4.2}
		\end{eqnarray}
		where $X$ runs over the simple objects of $\mathcal{C}_{V^G}$.
		
		Therefore, the multiplicity of $M^3$ on the left-hand side of (\ref{4.2}) is $$\sum_{U\in\mathcal{O}_{M^1}}
		\sum_{W\in\mathcal{O}_{M^2}}
		N^{M^3}_{U,W}.$$
		On the other hand, by Proposition \ref{pp 4.2} and Theorem \ref{thm 3.2}, the multiplicity of $M^3$ on
		the right-hand side (\ref{4.2}) is $$\sum_{\chi_3\in\operatorname{Irr}(G_{M^3})}
		N^{M^3_{\chi_3}}_{M^1_{\chi_1},M^2_{\chi_2}}.$$
		Hence, we have $$\sum_{U\in\mathcal{O}_{M^1}}
		\sum_{W\in\mathcal{O}_{M^2}}
		N^{M^3}_{U,W}
		=
		\sum_{\chi_3\in\operatorname{Irr}(G_{M^3})}
		N^{M^3_{\chi_3}}_{M^1_{\chi_1},M^2_{\chi_2}}.$$
		
		It concludes the proof.\end{proof}

	The following results follow immediately from Theorem~\ref{thm 4.3}.
	\bc{cor 4.4}
	$\mathrm{{(1)}}$ If $M^1,M^2\in \mathcal{S}_G$, then
	\begin{eqnarray*}
		N^{M^3}_{M^1,M^2}
		=
		\sum_{\chi_3\in\operatorname{Irr}(G_{M^3})}
		N^{M^3_{\chi_3}}_{M^1_{\Lambda_0},M^2_{\Lambda_0}},
	\end{eqnarray*}
	where $\Lambda_0$ is the trivial character of $G$.
	\\$\mathrm{{(2)}}$ If $M^2,M^3\in \mathcal{S}_G$, then
	\begin{eqnarray*}
		N^{M^3}_{M^1,M^2}
		=
		\frac{1}{[G:G_{M^1}]}
		\sum_{\chi_3\in\operatorname{Irr}(G)}
		N^{M^3_{\chi_3}}_{M^1_{\chi_1},M^2_{\Lambda_0}},
	\end{eqnarray*}
	where $\chi_1\in\operatorname{Irr}(G_{M^1})$ and $\Lambda_0$ denotes the trivial characters of $G$.
	\\$\mathrm{{(3)}}$ If $M^1,M^3\in \mathcal{S}_G$, then
	\begin{eqnarray*}
		N^{M^3}_{M^1,M^2}
		=
		\frac{1}{[G:G_{M^2}]}
		\sum_{\chi_3\in\operatorname{Irr}(G)}
		N^{M^3_{\chi_3}}_{M^1_{\Lambda_0},M^2_{\chi_2}},
	\end{eqnarray*}
	where $\chi_2\in\operatorname{Irr}(G_{M^2})$ and $\Lambda_0$ denotes the trivial characters of $G$.
	\ec

	\begin{proof}
		The first assertion follows immediately from Theorem~\ref{thm 4.3}.
		
		For the second assertion, by Theorem~\ref{thm 4.3}, we have
		\begin{eqnarray*}
			\sum_{U\in\mathcal{O}_{M^1}}
			N^{M^3}_{U,M^2}
			=
			\sum_{\chi_3\in\operatorname{Irr}(G)}
			N^{M^3_{\chi_3}}_{M^1_{\chi_1},M^2_{\Lambda_0}} .
		\end{eqnarray*}
		On the other hand, by Proposition~\ref{pp 4.1}, for any
		$U\in\mathcal{O}_{M^1}$, we have $N^{M^3}_{U,M^2}
		=
		N^{M^3}_{M^1,M^2}.$
		It follows that $|\mathcal{O}_{M^1}|N^{M^3}_{M^1,M^2}
		=
		\sum_{\chi_3\in\operatorname{Irr}(G)}
		N^{M^3_{\chi_3}}_{M^1_{\chi_1},M^2_{\Lambda_0}}.$
		Since $|\mathcal{O}_{M^1}|=[G:G_{M^1}],$
		we obtain
		\begin{eqnarray*}
			N^{M^3}_{M^1,M^2}
			=
			\frac{1}{[G:G_{M^1}]}
			\sum_{\chi_3\in\operatorname{Irr}(G)}
			N^{M^3_{\chi_3}}_{M^1_{\chi_1},M^2_{\Lambda_0}} .
		\end{eqnarray*}
		
		The third assertion follows similarly.
	\end{proof}

	As an application of Corollary~{\ref{cor 4.4}}, we show that the intertwining
	operators for $V^G$-modules can be realized as restrictions of
	intertwining operators for twisted $V$-modules. We first recall a
	result from Proposition 5.10 of \cite{dl} (see also  Proposition 11.9 of \cite{dle}).
	
	\bp{pp 4.5}
	Let $g_1,g_2,g_3\in G$ and $M^i$ be an irreducible
	$g_i$-twisted module for $i=1,2,3$. If $\mathcal{Y}$ is a nonzero intertwining operator of
	type $\binom{M^3}{M^1\,M^2},$
	then
	\begin{eqnarray*}
		\mathcal{Y}(u,x)v\neq0
	\end{eqnarray*}
	for any nonzero elements $u\in M^1$ and $v\in M^2$.
	\ep

	We next recall a map between intertwining operators for
	$V^G$-modules and those for twisted $V$-modules (cf. \cite{dl}). Let $g_1,g_2,g_3\in G$ and $M^i$ be an irreducible
	$g_i$-twisted module for $i=1,2,3$. Suppose that $N^i$  is an irreducible
	$V^G$-submodule of $M^i$ for $i=1,2$. In other words, $N^i=M^i_{\chi_i}$ for some $\chi_i\in \operatorname{Irr}(G_{M^i})$, $i=1,2$. We define a linear map
	\begin{eqnarray*}
		&&\Psi:
		I_V\binom{M^3}{M^1\,M^2}
		\to
		I_{V^G}\binom{M^3}{N^1\,N^2}\\
		&&\hspace{1.6cm}\mathcal{Y}(\cdot,x)\mapsto \mathcal{Y}^{r}(\cdot,x)
	\end{eqnarray*}
	where $\mathcal{Y}^{r}(\cdot,x)$ denotes the element in
	$I_{V^G}\binom{M^3}{N^1\,N^2}$ satisfying
	\[
	\mathcal{Y}^{r}(u,x)v=\mathcal{Y}(u,x)v
	\]
	for all $u\in N^1$ and $v\in N^2$.
	It follows immediately from Proposition~\ref{pp 4.5} that the above map has the following property (cf. \cite{dl}).
	\bp{pp 4.6}
	The linear map $\Psi$ is injective.
	\ep
	
	Combining Corollary~\ref{cor 4.4} and Proposition~\ref{pp 4.6}, we obtain
	the following result, which generalizes Theorem~5.12 of \cite{dl} to
	cyclic groups of arbitrary finite order.
	
	\bt{thm 4.7}
	Let $g_1,g_2,g_3\in G$ and $M^i$ be an irreducible
	$g_i$-twisted module for $i=1,2,3$. Suppose that $N^i$  is an irreducible
	$V^G$-submodule of $M^i$ for $i=1,2$. Assume that $M^1,M^2\in\mathcal{S}_G$.  Then the map
	\begin{eqnarray*}
		\Psi:
		I_V\binom{M^3}{M^1\,M^2}
		\to
		I_{V^G}\binom{M^3}{N^1\,N^2}
	\end{eqnarray*}
	is an isomorphism.
	\et
	
	\begin{proof}
		Since $M^1,M^2\in\mathcal{S}_G$, we have $G_{M^1}=G_{M^2}=G.$
		Thus, the decompositions of $M^1$ and $M^2$ are given by
		\begin{eqnarray*}
			M^1
			=
			\bigoplus_{0\leq i\leq n-1}M^1_{\Lambda_i},
			\quad
			M^2
			=
			\bigoplus_{0\leq i\leq n-1}M^2_{\Lambda_i}.
		\end{eqnarray*}
		Let $N^1=M^1_{\Lambda_s}, N^2=M^2_{\Lambda_t}$ for any $0\leq s,t\leq n-1$. By the decomposition $M^3
		=
		\bigoplus_{\chi\in\operatorname{Irr}(G_{M^3})}
		M^3_\chi$
		and Theorem~\ref{thm 4.3}, we have
		\begin{eqnarray*}
			N^{M^3}_{M^1,M^2}
			=
			\sum_{\chi\in\operatorname{Irr}(G_{M^3})}
			N^{M^3_\chi}_{M^1_{\Lambda_s},M^2_{\Lambda_t}} .
		\end{eqnarray*}
		Therefore, we obtain
		\begin{eqnarray*}
			\dim I_{V^G}\binom{M^3}{N^1\,N^2}
			&=&
			\dim
			\bigoplus_{\chi\in\operatorname{Irr}(G_{M^3})}
			I_{V^G}\binom{M^3_\chi}{M^1_{\Lambda_s}\,M^2_{\Lambda_t}}
			\\
			&=&
			\sum_{\chi\in\operatorname{Irr}(G_{M^3})}
			N^{M^3_\chi}_{M^1_{\Lambda_s},M^2_{\Lambda_t}}
			\\
			&=&
			N^{M^3}_{M^1,M^2}
			\\
			&=&
			\dim I_V\binom{M^3}{M^1\,M^2}.
		\end{eqnarray*}
		
		By Proposition~\ref{pp 4.6}, $\Psi$ is injective.
		Hence, $\Psi$ is also surjective, and therefore it is an isomorphism.
	\end{proof}

	We next focus on the case that $M^2,M^3\in\mathcal{S}_G$ and
	$M^1\not\in\mathcal{S}_G$. Let
	\begin{eqnarray*}
		G_{M^1}=\{h\in G\mid M^1\circ h\cong M^1\}
	\end{eqnarray*}
	be the stabilizer subgroup of $M^1$. Let $k$ runs over a fixed
	set of representatives of the left cosets of $G_{M^1}$ in $G$. Recall that the irreducible $V^G$-modules obtained from the
	decomposition with respect to the characters are isomorphic for
	elements in the same $G$-orbit. We define a linear map
	$\widetilde{\Psi}$ as follows
	\begin{eqnarray*}
		&&\widetilde{\Psi}:
		\bigoplus_{k\in G/G_{M^1}}
		I_V
		\binom{M^3}{M^1\circ k\,\, M^2}
		\to
		I_{V^G}
		\binom{M^3}{N^1\,N^2}\\
		&&\hspace{1.4cm}(\mathcal{Y}_{M^1\circ k}(\cdot,x))_{k\in G/G_{M^1}}
		\mapsto
		\sum_{k\in G/G_{M^1}}\mathcal{Y}^r_{M^1\circ k}(\cdot,x),
	\end{eqnarray*}
	where
	$\mathcal{Y}^{r}_{M^1\circ k}(\cdot,x)$ is an element of
	$I_{V^G}\binom{M^3}{N^1\,N^2}$ satisfying
	\begin{eqnarray*}
		\mathcal{Y}^{r}_{M^1\circ k}(u,x)v
		=
		\mathcal{Y}_{M^1\circ k}(u,x)v
	\end{eqnarray*}
	for any $u\in N^1$ and $v\in N^2$. Therefore, we obtain the following result about the linear map $\widetilde{\Psi}$.

	The following result generalizes Theorem~5.13 of \cite{dl} to arbitrary cyclic groups.
	\bt{thm 4.8}
	The linear map $\widetilde{\Psi}$ is an isomorphism.
	\et
	\begin{proof}
		Since $M^2,M^3\in\mathcal{S}_G$, we have the decompositions of $M^2$ and $M^3$:
		\begin{eqnarray*}
			M^2
			=
			\bigoplus_{0\leq i\leq n-1}M^2_{\Lambda_i},
			\quad
			M^3
			=
			\bigoplus_{0\leq i\leq n-1}M^3_{\Lambda_i}.
		\end{eqnarray*}
		Moreover, we have $M^1
		=
		\bigoplus_{\chi\in\operatorname{Irr}(G_{M^1})}
		M^1_\chi$.
		Let $N^1=M^1_{\chi}, N^2=M^2_{\Lambda_t}$ for any  $\chi\in \operatorname{Irr}(G_{M^1})$ and $0\leq t\leq n-1$. By an argument similar to that in the proof of Theorem~\ref{thm 4.7},
		we obtain
		\begin{eqnarray*}
			\sum_{k\in G/G_{M^1}}N^{M^3}_{M^1\circ k,M^2}
			&=&
			\sum_{i=0}^{n-1}N^{M^3_{\Lambda_i}}_{M^1_{\chi},M^2_{\Lambda_t}}.
		\end{eqnarray*}
		
		It suffices to show that $\widetilde{\Psi}$ is injective.
		Suppose that $G/G_{M^1}=\{k_1,k_2,\dots,k_s\}$. Then
		\begin{eqnarray*}
			&&\widetilde\Psi(
			\mathcal{Y}_{M^1\circ k_1}(\cdot,x),
			\mathcal{Y}_{M^1\circ k_2}(\cdot,x),\dots,
			\mathcal{Y}_{M^1\circ k_s}(\cdot,x)
			)\\
			&=&\mathcal{Y}_{M^1\circ k_1}^r(\cdot,x)+\mathcal{Y}_{M^1\circ k_2}^r(\cdot,x)+\cdots+\mathcal{Y}_{M^1\circ k_s}^r(\cdot,x)
		\end{eqnarray*}
		Assume that $\widetilde\Psi(
		\mathcal{Y}_{M^1\circ k_1}(\cdot,x),
		\mathcal{Y}_{M^1\circ k_2}(\cdot,x),\dots,
		\mathcal{Y}_{M^1\circ k_s}(\cdot,x)
		)=0.$
		Then, by Proposition~\ref{pp wa} together with the Jacobson density theorem,
		we have
		\begin{eqnarray*}
			0
			&=&
			(x_0+x_2)^{m(u,w^2)+\frac{j_2}{T_2}}
			Y_{M^3}(u,x_0+x_2)
			\left(
			\sum_{i=1}^{s}
			\mathcal{Y}^{r}_{M^1\circ k_i}(m_i,x)w^2
			\right)
			\\
			&=&
			(x_2+x_0)^{m(u,w^2)+\frac{j_2}{T_2}}
			\sum_{i=1}^{s}
			\mathcal{Y}^{r}_{M^1\circ k_i}
			(Y_{M^1}(u,x_0)m_i,x_2)w^2,
		\end{eqnarray*}
		where $u\in V^{(j_1,j_2)}$ with $j_1,j_2\in \mathbb{Z}$, $m_1,m_2,\dots,m_s\in N^1$, $w^2\in N^2$ and $m(u,w^2)$ is a non-negative integer such that
		\begin{eqnarray*}
			&&Y_{M^1}(u,x_0)m_1\neq0,\\
			&&Y_{M^1}(u,x_0)m_i=0\quad \text{for}\,i\neq 1.
		\end{eqnarray*}
		Consequently, we obtain
		\begin{eqnarray*}
			&(x_2+x_0)^{m(u,w^2)+\frac{j_2}{T_2}}
			\mathcal{Y}^{r}_{M^1\circ k_1}
			(Y_{M^1}(u,x_0)m_1,x_2)w^2
			=0.
		\end{eqnarray*}
		It follows that $$\mathcal{Y}^{r}_{M^1\circ k_1}
		(Y_{M^1}(u,x_0)m_1,x_2)w^2
		=0.$$
		By Proposition~\ref{pp 4.6}, we conclude that $\mathcal{Y}_{M^1\circ k_1}(\cdot,x)=0.$
		By a similar argument, we obtain $$\mathcal{Y}_{M^1\circ k_i}(\cdot,x)=0$$ for $i\neq 1$.
		Therefore, $\widetilde\Psi$ is injective.
		
		This completes the proof.
	\end{proof}

	\section{Twisted Verlinde Formula}
	In this section, we will establish a twisted Verlinde formula for
	a cyclic group acting on a vertex operator algebra. Throughout this
	section, we assume that $V$ is a simple, rational, $C_2$-cofinite and
	self-dual vertex operator algebra of CFT type. Let $\sigma$ be an
	automorphism of $V$ of finite order and $G=\langle\sigma\rangle$
	be the cyclic subgroup of $\operatorname{Aut}(V)$ generated by
	$\sigma$. By Theorem~\ref{thm 11}, the fixed-point subalgebra $V^G$
	is a simple, rational, $C_2$-cofinite and self-dual vertex operator
	algebra of CFT type. Moreover, by Proposition~\ref{pp 12}, $V$ is
	$g$-rational for every $g\in G$.
	The aim of this section is to prove the following
	twisted Verlinde formula for vertex operator algebras.

	\bt{thm 5.1}
	Assume that $V$ is a simple, rational, $C_2$-cofinite and self-dual vertex
	operator algebra of CFT type and that $\sigma$ is an automorphism of
	$V$ of finite order. Let $G=\langle\sigma\rangle$ be the cyclic
	subgroup of $\operatorname{Aut}(V)$ generated by $\sigma$. Let
	$g_1,g_2,g_3\in G$ be automorphisms of $V$ such that
	$g_1\neq1$ or $g_2\neq1$, and $g_3=g_1g_2$. Set $\mathcal{M}(1;g_1,g_2)=\mathcal{M}(1,g_1)\cap \mathcal{M}(1,g_2).$
	Then for any
	$M^i\in\mathcal{M}(g_i)$ for $i=1,2,3$, we have
	\begin{eqnarray*}
		N^{M^3}_{M^1,M^2}
		&=&
		\sum_{W\in\mathcal{M}(1; g_1,g_2)}
		\frac{
			S_{M^1,W}S_{M^2,W}\overline{S_{M^3,W}}
		}
		{S_{V,W}} .
	\end{eqnarray*}
	\et

	The proof of Theorem~\ref{thm 5.1} is divided into two steps. Since
	the order of $G$ admits a prime factorization, it suffices to consider
	the case
	\begin{eqnarray*}
		|G|=p_1p_2\cdots p_s,
	\end{eqnarray*}
	with $p_1,p_2,\ldots,p_s$ being primes, not necessarily distinct. We
	first establish the result for the case where $|G|=p_1p_2$ with two
	primes $p_1$ and $p_2$. Then, by mathematical induction, we extend the
	result to the general case.

	\subsection{Proof of Theorem 5.1: the Case $|G|=p_1p_2$}
	In this subsection, we prove Theorem~\ref{thm 5.1} in the case where
	$|G|=p_1p_2$ for primes $p_1$ and $p_2$. The main strategy is to apply
	the Verlinde formula for $V^G$. In Section~3, we establish a relation
	between the $S$-matrix arising from the orbifold theory of $V$ and the
	$S$-matrix of $V^G$. In Section~4, we establish a relation between the
	fusion rules of twisted $V$-modules and those of $V^G$-modules. These
	relations allow us to derive the twisted Verlinde formula from the
	ordinary Verlinde formula for $V^G$.
	
	Throughout the following
	discussion,  we assume that $V$ is a
	simple, rational, $C_2$-cofinite and self-dual vertex operator algebra of
	CFT type, and that $\sigma\in\operatorname{Aut}(V)$ has order $p_1p_2$ for primes $p_1$ and $p_2$.
	Let $G=\langle\sigma\rangle$ be the cyclic subgroup generated by $\sigma$. We first consider the case that $M^1, M^2, M^3\in \mathcal{S}_G$. We use $\mathcal{I}(V^G)$ to denote the set of
	representatives of inequivalent irreducible $V^G$-modules.
	For $g\in G$, let $\mathcal{R}(g)$ be a set of representatives of the
	$G$-orbits on $\mathcal{M}(g)$. For a subgroup $H\leq G$, define
	$\mathcal{R}_{H}(g)$ to be the subset of $\mathcal{R}(g)$ consisting of
	irreducible $g$-twisted $V$-modules $M$ whose stabilizer subgroup is
	$G_M=H$.

	The following theorem gives the special case of Theorem~\ref{thm 5.1}
	when $|G|=p_1p_2$ for primes $p_1$ and $p_2$.
	\bt{thm 5.2}
	Assume that $V$ is a simple, rational, $C_2$-cofinite and self-dual vertex
	operator algebra of CFT type and that $\sigma$ is an automorphism of
	$V$ of finite order $p_1p_2$ for primes $p_1$ and $p_2$. Let $G=\langle\sigma\rangle$ be the cyclic
	subgroup of $\operatorname{Aut}(V)$ generated by $\sigma$. Let
	$g_1,g_2,g_3\in G$ be automorphisms of $V$ such that
	$g_1\neq1$ or $g_2\neq1$, and $g_3=g_1g_2$. Set $\mathcal{M}(1;g_1,g_2)=\mathcal{M}(1,g_1)\cap \mathcal{M}(1,g_2).$
	Then for any
	$M^i\in\mathcal{M}(g_i)$ for $i=1,2,3$, we have
	\begin{eqnarray*}
		N^{M^3}_{M^1,M^2}
		&=&
		\sum_{W\in\mathcal{M}(1; g_1,g_2)}
		\frac{
			S_{M^1,W}S_{M^2,W}\overline{S_{M^3,W}}
		}
		{S_{V,W}} .
	\end{eqnarray*}
	\et

	We recall the following lemma from \cite{drx2}.
	\bl{lem 5.3}
	Let $M\in\mathcal{M}(g,h)$ for $g,h\in G$, and assume that
	$N\in\mathcal{M}(h,g^{-1})$. Then, for any $k\in G$, we have $S_{M,N}=S_{M\circ k,N\circ k}.$
	\el
	
	We first consider the case that $M^1, M^2, M^3\in \mathcal{S}_G$.
	\bp{pp 5.4}
	Let
	$g_1,g_2,g_3\in G$ be automorphisms of $V$ such that
	$g_1\neq1$ or $g_2\neq1$, and $g_3=g_1g_2$. Set $\mathcal{M}(1;g_1,g_2)=\mathcal{M}(1,g_1)\cap \mathcal{M}(1,g_2).$
	Then for any
	$M^i\in\mathcal{M}(g_i)\cap \mathcal{S}_G$ for $i=1,2,3$, we have
	\begin{eqnarray*}
		N^{M^3}_{M^1,M^2}
		&=&
		\sum_{W\in\mathcal{M}(1; g_1,g_2)}
		\frac{
			S_{M^1,W}S_{M^2,W}\overline{S_{M^3,W}}
		}
		{S_{V,W}} .
	\end{eqnarray*}
	\ep
	
	\begin{proof}
		By Corollary~\ref{cor 4.4} (1), we have
		\begin{eqnarray*}
			N^{M^3}_{M^1,M^2}
			&=&
			\sum_{\lambda\in\operatorname{Irr}(G)}
			N^{M^3_{\lambda}}_{M^1_{\Lambda_0},M^2_{\Lambda_0}},
		\end{eqnarray*}
		where $\Lambda_0$ denotes the trivial character of $G$.
		Applying Theorem~\ref{thm 25}, we obtain
		\begin{eqnarray*}
			N^{M^3_{\lambda}}_{M^1_{\Lambda_0},M^2_{\Lambda_0}}
			&=&
			\sum_{N\in\mathcal{I}(V^G)}
			\frac{
				S_{M^1_{\Lambda_0},N}
				S_{M^2_{\Lambda_0},N}
				\overline{S_{M^3_{\lambda},N}}
			}
			{S_{V^G,N}} .
		\end{eqnarray*}

		By Theorem~\ref{thm 3.2}, every irreducible $V^G$-module $N$ is of the
		form $M^{g,H,j}_{\chi}$. For simplicity, we write $N=W_{\chi}$, where
		$W$ is an irreducible $g_W$-twisted $V$-module.
		By Theorem~\ref{thm 3.4}, for $i=1,2$, we have
		\begin{eqnarray*}
			S_{M^i_{\Lambda_0},W_{\chi}}
			=
			\begin{cases}
				\dfrac{
					\overline{\Lambda_0(g_W)}\chi(g_i^{-1})
				}{|G|}
				\displaystyle\sum_{k\in C_{M^i,W}}S_{M^i,W\circ k},
				& g_i\in G_W,\\[1.2ex]
				0,
				& \text{otherwise}.
			\end{cases}
		\end{eqnarray*}
		Since $M^i$ is $G$-stable, we have
		$G_{M^i}=G$. Hence, $M^i\in\mathcal{M}(g_i,g_W^{-1}).$
		If $g_i\in G_W$, then
		$W\in\mathcal{M}(g_W,g_i^{-1})$. Consequently,
		$W\circ k\in\mathcal{M}(g_W,g_i^{-1})$ for
		$k\in C_{M^i,W}$. By Lemma~\ref{lem 5.3}, for any $k\in C_{M^i,W}$, we have $S_{M^i,W\circ k}
		=
		S_{M^i,W}.$ Thus
		\begin{eqnarray*}
			S_{M^i_{\Lambda_0},W_{\chi}}
			=
			\begin{cases}
				\dfrac{
					\chi(g_i^{-1})
				}{|G_W|}
				S_{M^i,W},
				& g_i\in G_W,\\[1.2ex]
				0,
				& \text{otherwise}.
			\end{cases}
		\end{eqnarray*}
		By a similar argument, we have
		\begin{eqnarray*}
			&&S_{M^3_{\lambda},W_{\chi}}
			=
			\begin{cases}
				\dfrac{
					\overline{\lambda(g_W)}\chi(g_3^{-1})
				}{|G_W|}
				S_{M^3,W},
				& g_3\in G_W,\\[1.2ex]
				0,
				& \text{otherwise},
			\end{cases}\\
			&&S_{V^G,W_{\chi}}=\frac{1}{|G_W|}S_{V,W}.
		\end{eqnarray*}
		Since $G$ is a finite abelian group, we have
		\begin{eqnarray*}
			\sum_{\lambda\in\operatorname{Irr}(G)}
			\lambda(g_W)
			=
			\begin{cases}
				|G|, & g_W=1,\\
				0, & g_W\neq1.
			\end{cases}
		\end{eqnarray*}
		Moreover, we have $\sum_{\chi\in\operatorname{Irr}(G_W)}
		{\chi(g_1^{-1})}
		{\chi(g_2^{-1})}
		\overline{\chi(g_3^{-1})}
		=
		|G_W|.$  Thus, we obtain
		\begin{eqnarray*}
			N^{M^3}_{M^1,M^2}
			&=&
			\sum_{\lambda\in\operatorname{Irr}(G)}
			N^{M^3_{\lambda}}_{M^1_{\Lambda_0},M^2_{\Lambda_0}}\\
			&=&\sum_{\lambda\in\operatorname{Irr}(G)}\sum_{N\in\mathcal{I}(V^G)}
			\frac{
				S_{M^1_{\Lambda_0},N}
				S_{M^2_{\Lambda_0},N}
				\overline{S_{M^3_{\lambda},N}}
			}
			{S_{V^G,N}}\\
			&=&\sum_{\lambda\in\operatorname{Irr}(G)}\sum_{g_W\in G}
			\sum_{W\in\mathcal{R}(g_W)}
			\sum_{\chi\in\operatorname{Irr}(G_W)}\frac{
				S_{M^1_{\Lambda_0},W_\chi}
				S_{M^2_{\Lambda_0},W_\chi}
				\overline{S_{M^3_{\lambda},W_\chi}}
			}
			{S_{V^G,W_\chi}}\\
			&=&\sum_{\lambda\in\operatorname{Irr}(G)}\sum_{g_W\in G}
			\sum_{W\in\mathcal{R}(g_W)}
			\sum_{\chi\in\operatorname{Irr}(G_W)}\sum_{\substack{g_1,g_2\in G_W\leq G}}\\
			&&\cdot\frac{
				\chi(g_1^{-1})
				\chi(g_2^{-1})
				\overline{\chi(g_3^{-1})}
				{\lambda(g_W)}
			}
			{|G_W|^2 }\frac{S_{M^1,W}
				S_{M^2,W}
				\overline{S_{M^3,W}}}{S_{V,W}}\\
			&=&\sum_{g_1,g_2\in G_W\leq G}
			\frac{|G|}{|G_W|}
			\sum_{W\in\mathcal{R}_{G_W}(1)}
			\frac{
				S_{M^1,W}
				S_{M^2,W}
				\overline{S_{M^3,W}}
			}
			{S_{V,W}}\\
			&=&
			\sum_{W\in\mathcal{M}(1; g_1,g_2)}
			\frac{
				S_{M^1,W}S_{M^2,W}\overline{S_{M^3,W}}
			}
			{S_{V,W}}.
		\end{eqnarray*}
		The last equality follows from Lemma~\ref{lem 5.3} and
		the length of the $G$-orbit containing $W$.
	\end{proof}
	
	We next consider the case that $M^1, M^2\in \mathcal{S}_G$ and $M^3\not\in \mathcal{S}_G$.
	\bp{pp 5.5}
	Let
	$g_1,g_2,g_3\in G$ be automorphisms of $V$ such that
	$g_1\neq1$ or $g_2\neq1$, and $g_3=g_1g_2$. Set $\mathcal{M}(1;g_1,g_2)=\mathcal{M}(1,g_1)\cap \mathcal{M}(1,g_2).$
	Suppose that
	$M^i\in\mathcal{M}(g_i)\cap\mathcal{S}_G$ for $i=1,2$, and that
	$M^3\notin\mathcal{S}_G$ is an irreducible $g_3$-twisted
	$V$-module. Then we have
	\begin{eqnarray*}
		N^{M^3}_{M^1,M^2}
		&=&
		\sum_{W\in\mathcal{M}(1; g_1,g_2)}
		\frac{
			S_{M^1,W}S_{M^2,W}\overline{S_{M^3,W}}
		}
		{S_{V,W}} .
	\end{eqnarray*}
	\ep
	
	\begin{proof}
		Let $\Lambda_0$ denote the trivial character of $G$. Combining
		Corollary~\ref{cor 4.4} (1), Theorem~\ref{thm 25} and
		Theorem~\ref{thm 3.2}, we obtain
		\begin{eqnarray*}
			N^{M^3}_{M^1,M^2}
			&=&
			\sum_{\lambda\in\operatorname{Irr}(G_{M^3})}
			N^{M^3_{\lambda}}_{M^1_{\Lambda_0},M^2_{\Lambda_0}}
			\\
			&=&
			\sum_{\lambda\in\operatorname{Irr}(G_{M^3})}
			\sum_{N\in\mathcal{I}(V^G)}
			\frac{
				S_{M^1_{\Lambda_0},N}
				S_{M^2_{\Lambda_0},N}
				\overline{S_{M^3_{\lambda},N}}
			}
			{S_{V^G,N}}
			\\
			&=&
			\sum_{\lambda\in\operatorname{Irr}(G_{M^3})}
			\sum_{g_W\in G}
			\sum_{W\in\mathcal{R}(g_W)}
			\sum_{\chi\in\operatorname{Irr}(G_W)}
			\frac{
				S_{M^1_{\Lambda_0},W_\chi}
				S_{M^2_{\Lambda_0},W_\chi}
				\overline{S_{M^3_{\lambda},W_\chi}}
			}
			{S_{V^G,W_\chi}} .
		\end{eqnarray*}
		By Theorem~\ref{thm 3.4} and Lemma~\ref{lem 5.3}, we have
		\begin{eqnarray*}
			&&S_{M^i_{\Lambda_0},W_{\chi}}
			=
			\begin{cases}
				\dfrac{
					\chi(g_i^{-1})
				}{|G_W|}
				S_{M^i,W},
				& g_i\in G_W\\[1.2ex]
				0,
				& \text{otherwise}
			\end{cases}\quad \text{for}\,i=1,2,\\
			&&S_{M^3_{\lambda},W_{\chi}}
			=
			\begin{cases}
				\dfrac{
					\overline{\lambda(g_W)}\chi(g_3^{-1})
				}{|G_{M^3}|}
				\displaystyle\sum_{k\in C_{M^3,W}}S_{M^3,W\circ k},
				& g_W\in G_{M^3},\ g_3\in G_{W},\\[1.2ex]
				0,
				& \text{otherwise},
			\end{cases}\\
			&&S_{V^G,W_{\chi}}=\frac{1}{|G_W|}S_{V,W}.
		\end{eqnarray*}
		Since $G_{M^3}$ is a finite abelian group, we have
		\begin{eqnarray*}
			\sum_{\lambda\in\operatorname{Irr}(G_{M^3})}
			\lambda(g_W)
			=
			\begin{cases}
				|G_{M^3}|, & g_W=1,\\
				0, & g_W\neq1.
			\end{cases}
		\end{eqnarray*}
		Moreover, we have $\sum_{\chi\in\operatorname{Irr}(G_W)}
		{\chi(g_1^{-1})}
		{\chi(g_2^{-1})}
		\overline{\chi(g_3^{-1})}
		=
		|G_W|.$  Thus, we obtain
		\begin{eqnarray*}
			N^{M^3}_{M^1,M^2}
			&=&
			\sum_{\lambda\in\operatorname{Irr}(G_{M^3})}
			\sum_{g_W\in G}
			\sum_{W\in\mathcal{R}(g_W)}
			\sum_{\chi\in\operatorname{Irr}(G_W)}
			\frac{
				S_{M^1_{\Lambda_0},W_\chi}
				S_{M^2_{\Lambda_0},W_\chi}
				\overline{S_{M^3_{\lambda},W_\chi}}
			}
			{S_{V^G,W_\chi}}\\
			&=&\sum_{\lambda\in\operatorname{Irr}(G_{M^3})}\sum_{g_W\in G_{M^3}}
			\sum_{W\in\mathcal{R}(g_W)}
			\sum_{\chi\in\operatorname{Irr}(G_W)}\sum_{\substack{g_1,g_2\in G_W\leq G}}\\
			&&\cdot\frac{
				\chi(g_1^{-1})
				\chi(g_2^{-1})
				\overline{\chi(g_3^{-1})}
				{\lambda(g_W)}
			}
			{|G_W||G_{M^3}| }\frac{S_{M^1,W}
				S_{M^2,W}
				\sum_{k\in C_{M^3,W}}\overline{S_{M^3,W\circ k}}}{S_{V,W}}\\
			&=&
			\sum_{W\in\mathcal{R}(1)}\sum_{g_1,g_2\in G_W\leq G}
			\frac{
				S_{M^1,W}
				S_{M^2,W}
				\sum_{k\in C_{M^3,W}}\overline{S_{M^3,W\circ k}}
			}
			{S_{V,W}}\\
			&=&\sum_{W\in\mathcal{R}(1)}\sum_{g_1,g_2\in G_W\leq G}\sum_{k\in C_{M^3,W}}
			\frac{
				S_{M^1,W\circ k}
				S_{M^2,W\circ k}
				\overline{S_{M^3,W\circ k}}
			}
			{S_{V,W\circ k}}\\
			&=&
			\sum_{W\in\mathcal{M}(1; g_1,g_2)}
			\frac{
				S_{M^1,W}S_{M^2,W}\overline{S_{M^3,W}}
			}
			{S_{V,W}}.
		\end{eqnarray*}
		
		This completes the proof.
	\end{proof}
	
	We proceed to consider the case where $M^2,M^3\in\mathcal{S}_G$ and
	$M^1\notin\mathcal{S}_G$.
	\bp{pp 5.6}
	Let
	$g_1,g_2,g_3\in G$ be automorphisms of $V$ such that
	$g_1\neq1$ or $g_2\neq1$, and $g_3=g_1g_2$. Set $\mathcal{M}(1;g_1,g_2)=\mathcal{M}(1,g_1)\cap \mathcal{M}(1,g_2).$
	Suppose that
	$M^i\in\mathcal{M}(g_i)\cap\mathcal{S}_G$ for $i=2,3$, and that
	$M^1\notin\mathcal{S}_G$ is an irreducible $g_1$-twisted
	$V$-module. Then we have
	\begin{eqnarray*}
		N^{M^3}_{M^1,M^2}
		&=&
		\sum_{W\in\mathcal{M}(1; g_1,g_2)}
		\frac{
			S_{M^1,W}S_{M^2,W}\overline{S_{M^3,W}}
		}
		{S_{V,W}} .
	\end{eqnarray*}
	\ep
	
	\begin{proof}
		Let $\Lambda_0$ denote the trivial character of $G$ and fix $\chi_1\in\operatorname{Irr}(G_{M^1})$. Applying
		Corollary~\ref{cor 4.4} (2), Theorem~\ref{thm 25}, and
		Theorem~\ref{thm 3.2}, we obtain
		\begin{eqnarray*}
			N^{M^3}_{M^1,M^2}
			&=&
			\frac{1}{[G:G_{M^1}]}
			\sum_{\lambda\in\operatorname{Irr}(G)}
			N^{M^3_{\lambda}}_{M^1_{\chi_1},M^2_{\Lambda_0}}
			\\
			&=&
			\frac{1}{[G:G_{M^1}]}\sum_{\lambda\in\operatorname{Irr}(G)}
			\sum_{N\in\mathcal{I}(V^G)}
			\frac{
				S_{M^1_{\chi_1},N}
				S_{M^2_{\Lambda_0},N}
				\overline{S_{M^3_{\lambda},N}}
			}
			{S_{V^G,N}}
			\\
			&=&
			\frac{|G_{M^1}|}{|G|}\sum_{\lambda\in\operatorname{Irr}(G)}
			\sum_{g_W\in G}
			\sum_{W\in\mathcal{R}(g_W)}
			\sum_{\chi\in\operatorname{Irr}(G_W)}
			\frac{
				S_{M^1_{\chi_1},W_\chi}
				S_{M^2_{\Lambda_0},W_\chi}
				\overline{S_{M^3_{\lambda},W_\chi}}
			}
			{S_{V^G,W_\chi}} .
		\end{eqnarray*}
		By Theorem~\ref{thm 3.4} and Lemma~\ref{lem 5.3}, we have
		\begin{eqnarray*}
			N^{M^3}_{M^1,M^2}
			&=&\frac{|G_{M^1}|}{|G|}\sum_{\lambda\in\operatorname{Irr}(G)}
			\sum_{g_W\in G_{M^1}}
			\sum_{W\in\mathcal{R}(g_W)}
			\sum_{\chi\in\operatorname{Irr}(G_W)}\sum_{\substack{g_1,g_2\in G_W\leq G}}\\
			&&\cdot\frac{
				\chi(g_1^{-1})
				\chi(g_2^{-1})
				\overline{\chi(g_3^{-1})}
				{\lambda(g_W)}\overline{\chi_1(g_W)}
			}
			{|G_W||G_{M^1}| }\frac{\sum_{k\in C_{M^1,W}}S_{M^1,W\circ k}
				S_{M^2,W}
				\overline{S_{M^3,W}}}{S_{V,W}}\\
			&=&
			\sum_{W\in\mathcal{R}(1)}\sum_{g_1,g_2\in G_W\leq G}\frac{\sum_{k\in C_{M^1,W}}S_{M^1,W\circ k}
				S_{M^2,W}
				\overline{S_{M^3,W}}}{S_{V,W}}\\
			&=&\sum_{W\in\mathcal{R}(1)}\sum_{g_1,g_2\in G_W\leq G}\sum_{k\in C_{M^1,W}}
			\frac{
				S_{M^1,W\circ k}
				S_{M^2,W\circ k}
				\overline{S_{M^3,W\circ k}}
			}
			{S_{V,W\circ k}}\\
			&=&
			\sum_{W\in\mathcal{M}(1; g_1,g_2)}
			\frac{
				S_{M^1,W}S_{M^2,W}\overline{S_{M^3,W}}
			}
			{S_{V,W}}.
		\end{eqnarray*}
		
		The proof is complete.
	\end{proof}
	
	Combining Proposition~\ref{pp 5.6} with Proposition~\ref{pp 24},
	we can handle the case where $M^1,M^3\in\mathcal{S}_G$ and
	$M^2\notin\mathcal{S}_G$.
	\bc{cor 5.7}
	Let
	$g_1,g_2,g_3\in G$ be automorphisms of $V$ such that
	$g_1\neq1$ or $g_2\neq1$, and $g_3=g_1g_2$. Set $\mathcal{M}(1;g_1,g_2)=\mathcal{M}(1,g_1)\cap \mathcal{M}(1,g_2).$
	Suppose that
	$M^i\in\mathcal{M}(g_i)\cap\mathcal{S}_G$ for $i=1,3$, and that
	$M^2\notin\mathcal{S}_G$ is an irreducible $g_2$-twisted
	$V$-module. Then we have
	\begin{eqnarray*}
		N^{M^3}_{M^1,M^2}
		&=&
		\sum_{W\in\mathcal{M}(1; g_1,g_2)}
		\frac{
			S_{M^1,W}S_{M^2,W}\overline{S_{M^3,W}}
		}
		{S_{V,W}} .
	\end{eqnarray*}
	\ec
	
	We have proved the cases in which all of $M^1,M^2,M^3$ belong to
	$\mathcal{S}_G$ or two of them belong to $\mathcal{S}_G$. We now consider the case where exactly one of $M^1,M^2,M^3$ belongs to
	$\mathcal{S}_G$. Before considering $M^1\in\mathcal{S}_G$, we need the
	following result.
	\bl{lem 5.8}
	Assume that
	$p_1\neq p_2$. Let
	$g_1,g_2,g_3\in G$ be automorphisms of $V$ such that
	$g_1\neq1$ or $g_2\neq1$, and $g_3=g_1g_2$.
	Set $\mathcal{M}(1;g_1,g_2)=\mathcal{M}(1,g_1)\cap\mathcal{M}(1,g_2)$.
	Suppose that
	$M^1\in\mathcal{M}(g_1)\cap\mathcal{S}_G$, and that
	$M^i\notin\mathcal{S}_G$ is an irreducible $g_i$-twisted
	$V$-module for $i=2,3$. Assume that $G_{M^2}G_{M^3}=G$, then we have
	\begin{eqnarray*}
		N^{M^3}_{M^1,M^2}
		&=&
		\sum_{W\in\mathcal{R}(1;g_1,g_2)}
		\sum_{k_1,k_2\in G/G_W}
		\frac{|G_W|}{p_1p_2}
		\frac{
			S_{M^1,W}
			S_{M^2,W\circ k_1}
			\overline{S_{M^3,W\circ k_2}}
		}
		{S_{V,W}},
	\end{eqnarray*}
	where $\mathcal{R}(1;g_1,g_2)$ denotes a set of representatives of the
	$G$-orbits on $\mathcal{M}(1;g_1,g_2)$, and the elements $k_1, k_2$
	are chosen as representatives of the cosets in $G/G_W$.
	In particular, if $g_2\neq 1$ and $g_3\neq 1$, then
	\begin{eqnarray*}
		N^{M^3}_{M^1,M^2}
		&=&
		\sum_{W\in\mathcal{M}(1; g_1,g_2)}
		\frac{
			S_{M^1,W}S_{M^2,W}\overline{S_{M^3,W}}
		}
		{S_{V,W}} .
	\end{eqnarray*}
	\el
	
	\begin{proof}	
		By Theorem~\ref{thm 4.3}, we have
		\begin{eqnarray*}
			\sum_{W\in\mathcal{O}_{M^2}}N^{M^3}_{M^1,W}
			&=&
			\sum_{\lambda\in\operatorname{Irr}(G_{M^3})}
			N^{M^3_{\lambda}}_{M^1_{\Lambda_0},M^2_{\chi_2}},
		\end{eqnarray*}
		where $\Lambda_0$ is the trivial character of $G$ and
		$\chi_2\in\operatorname{Irr}(G_{M^2})$.
		Since $G_{M^2}G_{M^3}=G$, Proposition~\ref{pp 4.1} implies that $N^{M^3}_{M^1,W}
		=
		N^{M^3}_{M^1,M^2}$
		for any $W\in\mathcal{O}_{M^2}$.
		Therefore,
		\begin{eqnarray*}
			N^{M^3}_{M^1,M^2}
			&=&
			\frac{1}{[G:G_{M^2}]}
			\sum_{\lambda\in\operatorname{Irr}(G_{M^3})}
			N^{M^3_{\lambda}}_{M^1_{\Lambda_0},M^2_{\chi_2}} .
		\end{eqnarray*}
		By Theorems~\ref{thm 25}, \ref{thm 3.2}, \ref{thm 3.4} and Lemma \ref{lem 5.3}, we have
		\begin{eqnarray*}
			N^{M^3}_{M^1,M^2}
			&=&
			\frac{1}{[G:G_{M^2}]}
			\sum_{\lambda\in\operatorname{Irr}(G_{M^3})}
			\sum_{g_W\in G}
			\sum_{W\in\mathcal{R}(g_W)}
			\sum_{\chi\in\operatorname{Irr}(G_W)}
			\frac{
				S_{M^1_{\Lambda_0},W_\chi}
				S_{M^2_{\chi_2},W_\chi}
				\overline{S_{M^3_{\lambda},W_\chi}}
			}
			{S_{V^G,W_\chi}}
			\\
			&=&
			\frac{|G_{M^2}|}{p_1p_2}
			\sum_{\lambda\in\operatorname{Irr}(G_{M^3})}
			\sum_{g_W\in G_{M^2}\cap G_{M^3}}
			\sum_{W\in\mathcal{R}(g_W)}
			\sum_{g_1,g_2\in G_W\leq G}
			\sum_{\chi\in\operatorname{Irr}(G_W)}
			\\
			&&\cdot
			\frac{
				\chi(g_1^{-1})
				\overline{\chi_2(g_W)}
				\chi(g_2^{-1})
				\lambda(g_W)
				\overline{\chi(g_3^{-1})}
			}
			{|G_{M^2}||G_{M^3}|}
			\frac{
				S_{M^1,W}
				\sum_{k_1\in C_{M^2,W}}S_{M^2,W\circ k_1}
			}
			{S_{V,W}}
			\\
			&&\cdot\sum_{k_2\in C_{M^3,W}}\overline{S_{M^3,W\circ k_2}}\\
			&=&
			\frac{1}{p_1p_2}
			\sum_{W\in\mathcal{R}(1;g_1,g_2)}
			\frac{
				|G_W|S_{M^1,W}
				\sum_{k_1\in C_{M^2,W}}
				S_{M^2,W\circ k_1}
				\sum_{k_2\in C_{M^3,W}}
				\overline{S_{M^3,W\circ k_2}}
			}
			{S_{V,W}}
			\\
			&=&
			\frac{1}{p_1p_2}
			\sum_{W\in\mathcal{R}(1;g_1,g_2)}
			\sum_{k_1,k_2\in G/G_W}
			|G_W|
			\frac{
				S_{M^1,W}
				S_{M^2,W\circ k_1}
				\overline{S_{M^3,W\circ k_2}}
			}
			{S_{V,W}},
		\end{eqnarray*}
		where $k_1$ and $k_2$ are chosen as representatives of the cosets in
		$G/G_W$.
		
		Since $G_{M^2}$ and $G_{M^3}$ are cyclic groups of prime order,
		if $g_2\neq1$ and $g_3\neq1$, then
		$G_{M^2}= \langle g_2\rangle$ and
		$G_{M^3}= \langle g_3\rangle$.
		Consequently,
		$\langle g_2,g_3\rangle=G$.
		Since $g_3=g_1g_2$, we also have
		$\langle g_1,g_2\rangle=G$.
		Therefore, $G_W=G$. So
		\begin{eqnarray*}
			N^{M^3}_{M^1,M^2}
			&=&
			\sum_{W\in\mathcal{M}(1; g_1,g_2)}
			\frac{
				S_{M^1,W}
				S_{M^2,W}
				\overline{S_{M^3,W}}
			}
			{S_{V,W}} .
		\end{eqnarray*}
		
		This completes the proof.
	\end{proof}
	
	We now prove the result in the case $M^1\in\mathcal{S}_G$ and
	$M^2,M^3\notin\mathcal{S}_G$.
	\bp{pp 5.9}
	Let
	$g_1,g_2,g_3\in G$ be automorphisms of $V$ such that
	$g_1\neq1$ or $g_2\neq1$, and $g_3=g_1g_2$. Set $\mathcal{M}(1;g_1,g_2)=\mathcal{M}(1,g_1)\cap \mathcal{M}(1,g_2).$
	Suppose that
	$M^1\in\mathcal{M}(g_1)\cap\mathcal{S}_G$, and that
	$M^i\notin\mathcal{S}_G$ is an irreducible $g_i$-twisted
	$V$-module for $i=2,3$. Then we have
	\begin{eqnarray*}
		N^{M^3}_{M^1,M^2}
		&=&
		\sum_{W\in\mathcal{M}(1; g_1,g_2)}
		\frac{
			S_{M^1,W}S_{M^2,W}\overline{S_{M^3,W}}
		}
		{S_{V,W}} .
	\end{eqnarray*}
	\ep
	\begin{proof}
		We consider two cases according to whether
		$\langle g_1,g_2\rangle\neq G$ or $\langle g_1,g_2\rangle=G$.
		
		\textbf{\emph{Case 1.}} Suppose that $\langle g_1,g_2\rangle\neq G.$
		Set $H
		=
		\langle g_1,g_2\rangle
		=
		\langle g_2,g_3\rangle$.
		
		If $g_2\neq1$ and $g_3\neq1$. Since $G_{M^2}$ and $G_{M^3}$ are cyclic groups of prime order, we have $G_{M^2}
		=
		\langle g_2\rangle\leq H$ and $G_{M^3}
		=
		\langle g_3\rangle\leq H$.
		Moreover, $H$ is a cyclic group of prime order. Hence $G_{M^2}=G_{M^3}=H.$
		Therefore, we may consider the $H$-orbifold. The stabilizer subgroup of $M^1$ in $H$ is given by $G_{M^1}^{H}
		=
		G_{M^1}\cap H
		=
		H.$
		By Theorem~5.1 of \cite{dl}, we obtain
		\begin{eqnarray*}
			N^{M^3}_{M^1,M^2}
			&=&
			\sum_{W\in\mathcal{M}(1;g_1,g_2)}
			\frac{
				S_{M^1,W}S_{M^2,W}\overline{S_{M^3,W}}
			}
			{S_{V,W}} .
		\end{eqnarray*}
		
		If $g_2=1$ and $g_3\neq1$. Then $H
		=
		\langle g_1\rangle
		=
		\langle g_3\rangle.$
		In this case, $G_{M^2}<G$ and $G_{M^3}=H.$
		Considering the $H$-orbifold again, we have
		\begin{eqnarray*}
			&&G_{M^1}^H=G_{M^1}\cap H
			=
			H,\\
			&&G_{M^2}^H=G_{M^2}\cap H
			=
			\{1\} \,\text{or}\, H.
		\end{eqnarray*}
		Therefore, by Theorem~5.1 of \cite{dl}, we obtain
		\begin{eqnarray*}
			N^{M^3}_{M^1,M^2}
			&=&
			\sum_{W\in\mathcal{M}(1;g_1,g_2)}
			\frac{
				S_{M^1,W}S_{M^2,W}\overline{S_{M^3,W}}
			}
			{S_{V,W}} .
		\end{eqnarray*}
		
		The case $g_2\neq1$ and $g_3=1$ follows from a similar argument to
		that for the case $g_2=1$ and $g_3\neq1$.

		\textbf{\emph{Case 2.}} Suppose that $\langle g_1,g_2\rangle=G.$
		Then $G
		=
		\langle g_1,g_2\rangle
		=
		\langle g_2,g_3\rangle.$

		If $g_2\neq1$ and $g_3\neq1$. Since
		$G_{M^2}$ and $G_{M^3}$ are cyclic groups of prime order, we have $G_{M^2}
		=
		\langle g_2\rangle$ and $G_{M^3}
		=
		\langle g_3\rangle.$
		Hence, $G_{M^2}G_{M^3}
		=
		G$.
		Suppose that $g_1\neq g_2$. By Lemma~\ref{lem 5.8}, we obtain
		\begin{eqnarray*}
			N^{M^3}_{M^1,M^2}
			&=&
			\sum_{W\in\mathcal{M}(1;g_1,g_2)}
			\frac{
				S_{M^1,W}S_{M^2,W}\overline{S_{M^3,W}}
			}
			{S_{V,W}} .
		\end{eqnarray*}
		On the other hand, in the case $g_1=g_2$, then
		$g_3=g_1g_2=g_2^2$, and hence
		$G_{M^2}G_{M^3}$ cannot be equal to $G$, which contradicts the
		assumption above.
		
		It remains to consider the cases where either $g_2=1$ or $g_3=1$.
		If $g_2=1$ and $g_3\neq1$, then $\langle g_3\rangle
		=
		G$,
		which implies that $G_{M^3}=G$.
		This contradicts the assumption
		$M^3\notin\mathcal{S}_G$.
		
		Similarly, if $g_2\neq1$ and $g_3=1$, then $G_{M^2}=G$,
		which contradicts the assumption
		$M^2\notin\mathcal{S}_G$.
	\end{proof}

	We next consider the result in the case $M^2\in\mathcal{S}_G$ and
	$M^1,M^3\notin\mathcal{S}_G$. We need the
	following result.
	\bl{lem 5.10}
	Assume that
	$p_1\neq p_2$. Let
	$g_1,g_2,g_3\in G$ be automorphisms of $V$ such that
	$g_1\neq1$ or $g_2\neq1$, and $g_3=g_1g_2$.
	Set $\mathcal{M}(1;g_1,g_2)=\mathcal{M}(1,g_1)\cap\mathcal{M}(1,g_2)$.
	Suppose that
	$M^2\in\mathcal{M}(g_2)\cap\mathcal{S}_G$, and that
	$M^i\notin\mathcal{S}_G$ is an irreducible $g_i$-twisted
	$V$-module for $i=1,3$. Assume that $G_{M^1}G_{M^3}=G$, then we have
	\begin{eqnarray*}
		N^{M^3}_{M^1,M^2}
		&=&
		\sum_{W\in\mathcal{R}(1;g_1,g_2)}
		\sum_{k_1,k_2\in G/G_W}
		\frac{|G_W|}{p_1p_2}
		\frac{
			S_{M^1,W\circ k_1}
			S_{M^2,W}
			\overline{S_{M^3,W\circ k_2}}
		}
		{S_{V,W}},
	\end{eqnarray*}
	where $\mathcal{R}(1;g_1,g_2)$ denotes a set of representatives of the
	$G$-orbits on $\mathcal{M}(1;g_1,g_2)$, and the elements $k_1, k_2$
	are chosen as representatives of the cosets in $G/G_W$.
	In particular, if $g_1\neq 1$ and $g_3\neq 1$, then
	\begin{eqnarray*}
		N^{M^3}_{M^1,M^2}
		&=&
		\sum_{W\in\mathcal{M}(1; g_1,g_2)}
		\frac{
			S_{M^1,W}S_{M^2,W}\overline{S_{M^3,W}}
		}
		{S_{V,W}} .
	\end{eqnarray*}
	\el
	
	\begin{proof}	
		By Theorem~\ref{thm 4.3}, we have
		\begin{eqnarray*}
			\sum_{W\in\mathcal{O}_{M^1}}N^{M^3}_{W,M^2}
			&=&
			\sum_{\lambda\in\operatorname{Irr}(G_{M^3})}
			N^{M^3_{\lambda}}_{M^1_{\chi_1},M^2_{\Lambda_0}},
		\end{eqnarray*}
		where $\Lambda_0$ is the trivial character of $G$ and
		$\chi_1\in\operatorname{Irr}(G_{M^1})$.
		Since $G_{M^1}G_{M^3}=G$, Proposition~\ref{pp 4.1} implies that $N^{M^3}_{W,M^2}
		=
		N^{M^3}_{M^1,M^2}$
		for any $W\in\mathcal{O}_{M^1}$.
		Therefore,
		\begin{eqnarray*}
			N^{M^3}_{M^1,M^2}
			&=&
			\frac{1}{[G:G_{M^1}]}
			\sum_{\lambda\in\operatorname{Irr}(G_{M^3})}
			N^{M^3_{\lambda}}_{M^1_{\chi_1},M^2_{\Lambda_0}} .
		\end{eqnarray*}
		By Theorems~\ref{thm 25}, \ref{thm 3.2}, \ref{thm 3.4} and Lemma \ref{lem 5.3}, we have
		\begin{eqnarray*}
			N^{M^3}_{M^1,M^2}
			&=&
			\frac{1}{[G:G_{M^1}]}
			\sum_{\lambda\in\operatorname{Irr}(G_{M^3})}
			\sum_{g_W\in G}
			\sum_{W\in\mathcal{R}(g_W)}
			\sum_{\chi\in\operatorname{Irr}(G_W)}
			\frac{
				S_{M^1_{\chi_1},W_\chi}
				S_{M^2_{\Lambda_0},W_\chi}
				\overline{S_{M^3_{\lambda},W_\chi}}
			}
			{S_{V^G,W_\chi}}
			\\
			&=&
			\frac{|G_{M^1}|}{p_1p_2}
			\sum_{\lambda\in\operatorname{Irr}(G_{M^3})}
			\sum_{g_W\in G_{M^1}\cap G_{M^3}}
			\sum_{W\in\mathcal{R}(g_W)}
			\sum_{g_1,g_2\in G_W\leq G}
			\sum_{\chi\in\operatorname{Irr}(G_W)}
			\\
			&&\cdot
			\frac{
				\chi(g_1^{-1})
				\overline{\chi_1(g_W)}
				\chi(g_2^{-1})
				\lambda(g_W)
				\overline{\chi(g_3^{-1})}
			}
			{|G_{M^1}||G_{M^3}|}
			\frac{
				\sum_{k_1\in C_{M^1,W}}S_{M^1,W\circ k_1}
				S_{M^2,W}
			}
			{S_{V,W}}
			\\
			&&\cdot\sum_{k_2\in C_{M^3,W}}\overline{S_{M^3,W\circ k_2}}\\
			&=&
			\frac{1}{p_1p_2}
			\sum_{W\in\mathcal{R}(1;g_1,g_2)}
			\frac{
				|G_W|\sum_{k_1\in C_{M^1,W}}S_{M^1,W\circ k_1}
				S_{M^2,W}
				\sum_{k_2\in C_{M^3,W}}
				\overline{S_{M^3,W\circ k_2}}
			}
			{S_{V,W}}
			\\
			&=&
			\frac{1}{p_1p_2}
			\sum_{W\in\mathcal{R}(1;g_1,g_2)}
			\sum_{k_1,k_2\in G/G_W}
			|G_W|
			\frac{
				S_{M^1,W\circ k_1}
				S_{M^2,W}
				\overline{S_{M^3,W\circ k_2}}
			}
			{S_{V,W}},
		\end{eqnarray*}
		where $k_1$ and $k_2$ are chosen as representatives of the cosets in
		$G/G_W$.
		
		Since $G_{M^1}$ and $G_{M^3}$ are cyclic groups of prime order,
		if $g_1\neq1$ and $g_3\neq1$, then
		$G_{M^1}= \langle g_1\rangle$ and
		$G_{M^3}= \langle g_3\rangle$.
		Consequently,
		$\langle g_1,g_3\rangle=G$.
		Therefore, $G_W=G$. So
		\begin{eqnarray*}
			N^{M^3}_{M^1,M^2}
			&=&
			\sum_{W\in\mathcal{M}(1; g_1,g_2)}
			\frac{
				S_{M^1,W}
				S_{M^2,W}
				\overline{S_{M^3,W}}
			}
			{S_{V,W}} .
		\end{eqnarray*}
		
		The proof is concluded.
	\end{proof}
	
	We proceed to prove the result for the case
	$M^2\in\mathcal{S}_G$ and $M^1,M^3\notin\mathcal{S}_G$.
	\bp{pp 5.11}
	Let
	$g_1,g_2,g_3\in G$ be automorphisms of $V$ such that
	$g_1\neq1$ or $g_2\neq1$, and $g_3=g_1g_2$. Set $\mathcal{M}(1;g_1,g_2)=\mathcal{M}(1,g_1)\cap \mathcal{M}(1,g_2).$
	Suppose that
	$M^2\in\mathcal{M}(g_2)\cap\mathcal{S}_G$, and that
	$M^i\notin\mathcal{S}_G$ is an irreducible $g_i$-twisted
	$V$-module for $i=1,3$. Then we have
	\begin{eqnarray*}
		N^{M^3}_{M^1,M^2}
		&=&
		\sum_{W\in\mathcal{M}(1; g_1,g_2)}
		\frac{
			S_{M^1,W}S_{M^2,W}\overline{S_{M^3,W}}
		}
		{S_{V,W}} .
	\end{eqnarray*}
	\ep
	
	\begin{proof}
		The proof is divided into two cases according to whether
		$\langle g_1,g_2\rangle=G$ or not.
		
		\textbf{\emph{Case 1.}} Suppose that $\langle g_1,g_2\rangle\neq G.$
		Set $H
		=
		\langle g_1,g_2\rangle
		=
		\langle g_1,g_3\rangle$.
		
		If $g_1\neq1$ and $g_3\neq1$. Since $G_{M^1}$, $G_{M^3}$ and $H$ are cyclic groups of prime order, we have $G_{M^1}
		=
		\langle g_1\rangle= H$ and $G_{M^3}
		=
		\langle g_3\rangle= H$.
		Therefore, we may consider the $H$-orbifold. The stabilizer subgroup of $M^2$ in $H$ is given by $G_{M^2}^{H}
		=
		G_{M^2}\cap H
		=
		H.$
		By Theorem~5.1 of \cite{dl}, we obtain
		\begin{eqnarray*}
			N^{M^3}_{M^1,M^2}
			&=&
			\sum_{W\in\mathcal{M}(1;g_1,g_2)}
			\frac{
				S_{M^1,W}S_{M^2,W}\overline{S_{M^3,W}}
			}
			{S_{V,W}} .
		\end{eqnarray*}
		
		If $g_1=1$ and $g_3\neq1$. Then $H
		=
		\langle g_2\rangle
		=
		\langle g_3\rangle.$
		In this case, $G_{M^1}<G$ and $G_{M^3}=H.$
		Considering the $H$-orbifold again, by Theorem~5.1 of \cite{dl}, we obtain
		\begin{eqnarray*}
			N^{M^3}_{M^1,M^2}
			&=&
			\sum_{W\in\mathcal{M}(1;g_1,g_2)}
			\frac{
				S_{M^1,W}S_{M^2,W}\overline{S_{M^3,W}}
			}
			{S_{V,W}} .
		\end{eqnarray*}
		
		The case $g_1\neq1$ and $g_3=1$ follows from a similar argument to
		that for the case $g_1=1$ and $g_3\neq1$.

		\textbf{\emph{Case 2.}} Suppose that $\langle g_1,g_2\rangle=G.$
		Then $G
		=
		\langle g_1,g_2\rangle
		=
		\langle g_1,g_3\rangle.$

		If $g_1\neq1$ and $g_3\neq1$. Since
		$G_{M^1}$ and $G_{M^3}$ are cyclic groups of prime order, we have $G_{M^1}
		=
		\langle g_1\rangle$ and $G_{M^3}
		=
		\langle g_3\rangle.$
		Hence, $G_{M^1}G_{M^3}
		=
		G$.
		Suppose that $g_1\neq g_3$. By Lemma~\ref{lem 5.10}, we obtain
		\begin{eqnarray*}
			N^{M^3}_{M^1,M^2}
			&=&
			\sum_{W\in\mathcal{M}(1;g_1,g_2)}
			\frac{
				S_{M^1,W}S_{M^2,W}\overline{S_{M^3,W}}
			}
			{S_{V,W}} .
		\end{eqnarray*}
		On the other hand, in the case $g_1=g_3$, then
		$G_{M^1}G_{M^3}$ cannot be equal to $G$, which contradicts the
		assumption above.
		
		It remains to consider the cases where either $g_1=1$ or $g_3=1$.
		If $g_1=1$ and $g_3\neq1$, then $\langle g_3\rangle
		=
		G$,
		which implies that $G_{M^3}=G$.
		This contradicts the assumption
		$M^3\notin\mathcal{S}_G$.
		
		Similarly, if $g_1\neq1$ and $g_3=1$, then $G_{M^1}=G$,
		which contradicts the assumption
		$M^1\notin\mathcal{S}_G$.
	\end{proof}
	
	We consider the case $M^3\in\mathcal{S}_G$ and
	$M^1,M^2\notin\mathcal{S}_G$. We first establish the following result.
	\bl{lem 5.12}
	Assume that
	$p_1\neq p_2$. Let
	$g_1,g_2,g_3\in G$ be automorphisms of $V$ such that
	$g_1\neq1$ or $g_2\neq1$, and $g_3=g_1g_2$.
	Set $\mathcal{M}(1;g_1,g_2)=\mathcal{M}(1,g_1)\cap\mathcal{M}(1,g_2)$.
	Suppose that
	$M^3\in\mathcal{M}(g_3)\cap\mathcal{S}_G$, and that
	$M^i\notin\mathcal{S}_G$ is an irreducible $g_i$-twisted
	$V$-module for $i=1,2$. Assume that $G_{M^1}G_{M^2}=G$, then we have
	\begin{eqnarray*}
		N^{M^3}_{M^1,M^2}
		&=&
		\sum_{W\in\mathcal{R}(1;g_1,g_2)}
		\sum_{k_1,k_2\in G/G_W}
		\frac{|G_W|}{p_1p_2}
		\frac{
			S_{M^1,W\circ k_1}
			S_{M^2,W\circ k_2}
			\overline{S_{M^3,W}}
		}
		{S_{V,W}},
	\end{eqnarray*}
	where $\mathcal{R}(1;g_1,g_2)$ denotes a set of representatives of the
	$G$-orbits on $\mathcal{M}(1;g_1,g_2)$, and the elements $k_1, k_2$
	are chosen as representatives of the cosets in $G/G_W$.
	In particular, if $g_1\neq 1$ and $g_2\neq 1$, then
	\begin{eqnarray*}
		N^{M^3}_{M^1,M^2}
		&=&
		\sum_{W\in\mathcal{M}(1; g_1,g_2)}
		\frac{
			S_{M^1,W}S_{M^2,W}\overline{S_{M^3,W}}
		}
		{S_{V,W}} .
	\end{eqnarray*}
	\el
	
	\begin{proof}	
		By Theorem~\ref{thm 4.3}, we have
		\begin{eqnarray*}
			\sum_{U\in\mathcal{O}_{M^1}}\sum_{W\in\mathcal{O}_{M^2}}N^{M^3}_{U,W}
			&=&
			\sum_{\lambda\in\operatorname{Irr}(G)}
			N^{M^3_{\lambda}}_{M^1_{\chi_1},M^2_{\chi_2}},
		\end{eqnarray*}
		where $\chi_1\in\operatorname{Irr}(G_{M^1})$ and $\chi_2\in\operatorname{Irr}(G_{M^2})$.
		Since $G_{M^1}G_{M^2}=G$, Proposition~\ref{pp 4.1} implies that $N^{M^3}_{U,W}
		=
		N^{M^3}_{M^1,M^2}$
		for any $U\in\mathcal{O}_{M^1}, W\in\mathcal{O}_{M^2}$.
		Therefore,
		\begin{eqnarray*}
			N^{M^3}_{M^1,M^2}
			&=&
			\frac{1}{[G:G_{M^1}]}\frac{1}{[G:G_{M^2}]}
			\sum_{\lambda\in\operatorname{Irr}(G)}
			N^{M^3_{\lambda}}_{M^1_{\chi_1},M^2_{\chi_2}} .
		\end{eqnarray*}
		By Theorems~\ref{thm 25}, \ref{thm 3.2}, \ref{thm 3.4} and Lemma \ref{lem 5.3}, we have
		\begin{eqnarray*}
			N^{M^3}_{M^1,M^2}
			&=&
			\frac{|G_{M^1}||G_{M^2}|}{|G|^2}
			\sum_{\lambda\in\operatorname{Irr}(G)}
			\sum_{g_W\in G}
			\sum_{W\in\mathcal{R}(g_W)}
			\sum_{\chi\in\operatorname{Irr}(G_W)}
			\frac{
				S_{M^1_{\chi_1},W_\chi}
				S_{M^2_{\chi_2},W_\chi}
				\overline{S_{M^3_{\lambda},W_\chi}}
			}
			{S_{V^G,W_\chi}}
			\\
			&=&
			\frac{|G_{M^1}||G_{M^2}|}{p_1^2p_2^2}
			\sum_{\lambda\in\operatorname{Irr}(G)}
			\sum_{g_W\in G_{M^1}\cap G_{M^2}}
			\sum_{W\in\mathcal{R}(g_W)}
			\sum_{g_1,g_2\in G_W\leq G}
			\sum_{\chi\in\operatorname{Irr}(G_W)}
			\\
			&&\cdot
			\frac{
				\overline{\chi_1(g_W)}\chi(g_1^{-1})\overline{\chi_2(g_W)}	\chi(g_2^{-1})\overline{\chi(g_3^{-1})}
				\lambda(g_W)
			}
			{|G_{M^1}||G_{M^2}|}
			\frac{
				\sum_{k_1\in C_{M^1,W}}S_{M^1,W\circ k_1}
				\overline {S_{M^3,W}}
			}
			{S_{V,W}}
			\\
			&&\cdot\sum_{k_2\in C_{M^2,W}}{S_{M^2,W\circ k_2}}\\
			&=&
			\frac{1}{p_1p_2}
			\sum_{W\in\mathcal{R}(1;g_1,g_2)}
			\frac{
				|G_W|\sum_{k_1\in C_{M^1,W}}S_{M^1,W\circ k_1}
				\sum_{k_2\in C_{M^2,W}}{S_{M^2,W\circ k_2}}
				\overline{S_{M^3,W}}
			}
			{S_{V,W}}
			\\
			&=&
			\frac{1}{p_1p_2}
			\sum_{W\in\mathcal{R}(1;g_1,g_2)}
			\sum_{k_1,k_2\in G/G_W}
			|G_W|
			\frac{
				S_{M^1,W\circ k_1}
				S_{M^2,W\circ k_2}
				\overline{S_{M^3,W}}
			}
			{S_{V,W}},
		\end{eqnarray*}
		where $k_1$ and $k_2$ are chosen as representatives of the cosets in
		$G/G_W$.
		
		Since $G_{M^1}$ and $G_{M^2}$ are cyclic groups of prime order,
		if $g_1\neq1$ and $g_2\neq1$, then
		$\langle g_1,g_2\rangle=G$.
		Therefore, $G_W=G$. Therefore we have
		\begin{eqnarray*}
			N^{M^3}_{M^1,M^2}
			&=&
			\sum_{W\in\mathcal{M}(1; g_1,g_2)}
			\frac{
				S_{M^1,W}
				S_{M^2,W}
				\overline{S_{M^3,W}}
			}
			{S_{V,W}} .
		\end{eqnarray*}
		
		This completes the proof.
	\end{proof}
	
	We proceed with the case $M^3\in\mathcal{S}_G$ and
	$M^1,M^2\notin\mathcal{S}_G$.
	\bp{pp 5.13}
	Let
	$g_1,g_2,g_3\in G$ be automorphisms of $V$ such that
	$g_1\neq1$ or $g_2\neq1$, and $g_3=g_1g_2$. Set $\mathcal{M}(1;g_1,g_2)=\mathcal{M}(1,g_1)\cap \mathcal{M}(1,g_2).$
	Suppose that
	$M^3\in\mathcal{M}(g_3)\cap\mathcal{S}_G$, and that
	$M^i\notin\mathcal{S}_G$ is an irreducible $g_i$-twisted
	$V$-module for $i=1,2$. Then we have
	\begin{eqnarray*}
		N^{M^3}_{M^1,M^2}
		&=&
		\sum_{W\in\mathcal{M}(1; g_1,g_2)}
		\frac{
			S_{M^1,W}S_{M^2,W}\overline{S_{M^3,W}}
		}
		{S_{V,W}} .
	\end{eqnarray*}
	\ep
	
	\begin{proof}
		We divide the proof into two cases depending on the subgroup generated by
		$g_1$ and $g_2$.
		
		\textbf{\emph{Case 1.}} Suppose that $\langle g_1,g_2\rangle\neq G.$
		Set $H
		=
		\langle g_1,g_2\rangle$.
		
		If $g_1\neq1$ and $g_2\neq1$. Since $G_{M^1}$, $G_{M^2}$ and $H$ are cyclic groups of prime order, we have $G_{M^1}
		=
		\langle g_1\rangle= H$ and $G_{M^2}
		=
		\langle g_2\rangle= H$.
		Therefore, we may consider the $H$-orbifold. The stabilizer subgroup of $M^3$ in $H$ is given by $G_{M^3}^{H}
		=
		G_{M^3}\cap H
		=
		H.$
		By Theorem~5.1 of \cite{dl}, we obtain
		\begin{eqnarray*}
			N^{M^3}_{M^1,M^2}
			&=&
			\sum_{W\in\mathcal{M}(1;g_1,g_2)}
			\frac{
				S_{M^1,W}S_{M^2,W}\overline{S_{M^3,W}}
			}
			{S_{V,W}} .
		\end{eqnarray*}
		
		If $g_1=1$ and $g_2\neq1$. Then $H
		=
		\langle g_2\rangle.$
		In this case, $G_{M^2}=H.$
		Considering the $H$-orbifold again, by Theorem~5.1 of \cite{dl}, we obtain
		\begin{eqnarray*}
			N^{M^3}_{M^1,M^2}
			&=&
			\sum_{W\in\mathcal{M}(1;g_1,g_2)}
			\frac{
				S_{M^1,W}S_{M^2,W}\overline{S_{M^3,W}}
			}
			{S_{V,W}} .
		\end{eqnarray*}
		
		The case $g_1\neq1$ and $g_2=1$ follows from a similar argument to
		that for the case $g_1=1$ and $g_2\neq1$.

		\textbf{\emph{Case 2.}} Suppose that $\langle g_1,g_2\rangle=G.$

		If $g_1\neq1$ and $g_2\neq1$. Since
		$G_{M^1}$ and $G_{M^2}$ are cyclic groups of prime order, we have $G_{M^1}
		=
		\langle g_1\rangle$ and $G_{M^2}
		=
		\langle g_2\rangle.$
		Hence, $G_{M^1}G_{M^2}
		=
		G$.
		Suppose that $g_1\neq g_2$. By Lemma~\ref{lem 5.12}, we obtain
		\begin{eqnarray*}
			N^{M^3}_{M^1,M^2}
			&=&
			\sum_{W\in\mathcal{M}(1;g_1,g_2)}
			\frac{
				S_{M^1,W}S_{M^2,W}\overline{S_{M^3,W}}
			}
			{S_{V,W}} .
		\end{eqnarray*}
		On the other hand, in the case $g_1=g_2$, then
		$G_{M^1}G_{M^2}$ cannot be equal to $G$, which contradicts the
		assumption above.
		
		It remains to consider the cases where either $g_1=1$ or $g_2=1$.
		If $g_1=1$ and $g_2\neq1$, then $\langle g_2\rangle
		=
		G$,
		which implies that $G_{M^2}=G$.
		This contradicts the assumption
		$M^2\notin\mathcal{S}_G$.
		
		Similarly, if $g_1\neq1$ and $g_2=1$, then $G_{M^1}=G$,
		which contradicts the assumption
		$M^1\notin\mathcal{S}_G$.
	\end{proof}
	
	It remains to consider the case
	$M^1,M^2,M^3\notin\mathcal{S}_G$.
	\bp{pp 5.14}
	Let
	$g_1,g_2,g_3\in G$ be automorphisms of $V$ such that
	$g_1\neq1$ or $g_2\neq1$, and $g_3=g_1g_2$. Set $\mathcal{M}(1;g_1,g_2)=\mathcal{M}(1,g_1)\cap \mathcal{M}(1,g_2).$
	Suppose that
	$M^i\notin\mathcal{S}_G$ is an irreducible $g_i$-twisted
	$V$-module for $i=1,2,3$. Then we have
	\begin{eqnarray*}
		N^{M^3}_{M^1,M^2}
		&=&
		\sum_{W\in\mathcal{M}(1; g_1,g_2)}
		\frac{
			S_{M^1,W}S_{M^2,W}\overline{S_{M^3,W}}
		}
		{S_{V,W}} .
	\end{eqnarray*}
	\ep
	\begin{proof}
		The following two cases arise according to the triviality of
		$g_1$, $g_2$, and $g_3$.
		
		\textbf{\emph{Case 1.}} Suppose that $g_i\neq 1$ for $i=1,2,3.$ Since $G_{M^i}$ is a cyclic group of prime order, we have $G_{M^i}
		=
		\langle g_i\rangle$ for $i=1,2,3.$

		If $p_1=p_2$, then $G$ is a cyclic group of order $p_1^2$.
		Since $G_{M^1}$ and $G_{M^2}$ have the same prime order, it follows that $G_{M^1}
		=
		G_{M^2}.$
		
		If $p_1\neq p_2$, we claim that $G_{M^1}
		=
		G_{M^2}.$
		Otherwise, assume that $G_{M^1}
		\neq
		G_{M^2}.$
		Since $G_{M^1}$ and $G_{M^2}$ are distinct subgroups of the cyclic
		group $G$, without loss of generality, we may assume that
		\begin{eqnarray*}
			G_{M^1}
			=
			\langle\sigma^{p_1}\rangle, G_{M^2}
			=
			\langle\sigma^{p_2}\rangle ,
		\end{eqnarray*}
		where $G_{M^1}$ and $G_{M^2}$ have orders $p_2$ and $p_1$,
		respectively. Hence, we have $G_{M^1}G_{M^2}
		=
		G.$
		Since $g_i\in G_{M^i}$ for $i=1,2$, there exist integers $a_1$ and
		$a_2$ satisfying $1\leq a_1\leq p_2-1, 1\leq a_2\leq p_1-1$
		such that
		\begin{eqnarray*}
			g_1
			=
			\sigma^{p_1a_1}, g_2
			=
			\sigma^{p_2a_2}.
		\end{eqnarray*}
		Therefore, $g_3
		=
		g_1g_2
		=
		\sigma^{p_1a_1+p_2a_2}.$
		Since $\gcd(p_1a_1+p_2a_2,p_1p_2)
		=
		1,$
		we obtain
		\begin{eqnarray*}
			G
			=
			\langle\sigma^{p_1a_1+p_2a_2}\rangle
			=
			\langle g_1g_2\rangle
			=
			\langle g_3\rangle
			=
			G_{M^3}.
		\end{eqnarray*}
		This contradicts the assumption that
		$M^3\notin\mathcal{S}_G$.
		Therefore, $G_{M^1}
		=
		G_{M^2}.$
		Thus, we have $G_{M^1}
		=
		G_{M^2}=G_{M^3}.$
		Set $H
		=
		G_{M^1}.$
		Then $H$ is a cyclic group of prime order. Applying
		Theorem~5.1 of \cite{dl} to the $H$-orbifold, we obtain
		\begin{eqnarray*}
			N^{M^3}_{M^1,M^2}
			&=&
			\sum_{W\in\mathcal{M}(1;g_1,g_2)}
			\frac{
				S_{M^1,W}S_{M^2,W}\overline{S_{M^3,W}}
			}
			{S_{V,W}} .
		\end{eqnarray*}

		\textbf{\emph{Case 2.}} Suppose that exactly one of $g_1$, $g_2$, and $g_3$
		is equal to $1$.
		
		If $g_1=1$, then $g_2=g_3$, and hence $	G_{M^2}
		=
		G_{M^3}.$
		Set $H
		=
		G_{M^2}.$
		Then $H$ is a cyclic group of prime order. The stabilizer subgroup of
		$M^1$ in $H$ is given by
		\begin{eqnarray*}
			G_{M^1}^{H}
			=
			G_{M^1}\cap H
			=
			\{1\}\,\text{or}\,H .
		\end{eqnarray*}
		Applying Theorem~5.1 of \cite{dl} to the $H$-orbifold, we obtain
		\begin{eqnarray*}
			N^{M^3}_{M^1,M^2}
			&=&
			\sum_{W\in\mathcal{M}(1;g_1,g_2)}
			\frac{
				S_{M^1,W}S_{M^2,W}\overline{S_{M^3,W}}
			}
			{S_{V,W}} .
		\end{eqnarray*}
		
		The cases $g_2=1$ and $g_3=1$ can be obtained by a similar argument.	
	\end{proof}
	
	Combining Propositions~\ref{pp 5.4}, \ref{pp 5.5}, \ref{pp 5.6},
	\ref{pp 5.9}, \ref{pp 5.11}, \ref{pp 5.13}, \ref{pp 5.14} with
	Corollary~\ref{cor 5.7}, we obtain Theorem~\ref{thm 5.2}.

	\subsection{Proof of Theorem 5.1: the Case $|G|=p_1p_2\cdots p_s$}
	In Section 5.1, we established the result for the case that
	$G$ is a cyclic group of order $p_1p_2$. In this subsection, we extend
	the argument to the general case
	\begin{eqnarray*}
		|G|=p_1p_2\cdots p_s,
	\end{eqnarray*}
	with $p_1,p_2,\ldots,p_s$ being primes, not necessarily distinct.
	
	We use the notation $\mathcal{I}(V^G)$, $\mathcal{R}(g)$, and
	$\mathcal{R}_{H}(g)$ introduced in Section 5.1. We proceed with the
	proof of Theorem~\ref{thm 5.1}.
	\begin{proof}
		By the prime factorization of $|G|$, we have
		\begin{eqnarray*}
			|G|
			=
			p_1p_2\cdots p_s
		\end{eqnarray*}
		with $p_1,p_2,\ldots,p_s$ being primes that are not necessarily
		distinct.
		We argue
		by induction on $s$. The cases $s=1$ and $s=2$ follow from
		\cite[Theorem~5.1]{dl} and Theorem~5.2, respectively. Let $s\geq 3$, and assume that the assertion holds for every
		cyclic group whose order is a product of $t$ primes, where $1\leq t<s$.
		We prove the assertion for the cyclic group $G$ of order
		$p_1p_2\cdots p_s$.
		We next divide the proof into two cases according to whether
		$\langle g_1,g_2\rangle$ equals $G$ or not.

		\textbf{\emph{Case 1.}} Suppose that $\langle g_1,g_2\rangle< G.$
		Set $H
		=
		\langle g_1,g_2\rangle$.
		Since $g_1\neq1$ or
		$g_2\neq1$, we have $H\neq \{1\}$. Moreover, $H$ is a proper subgroup of $G$.
		Hence, the
		order of $H$ has the form
		\begin{eqnarray*}
			|H|
			=
			p_{i_1}p_{i_2}\cdots p_{i_r}
			\quad \text{for}\,1\leq r\leq s-1,
		\end{eqnarray*}
		where $i_1,i_2,\ldots,i_r$ are distinct indices in
		$\{1,2,\ldots,s\}$.
		Consider the $H$-orbifold.
		Note that $g_i\in H$ for $i=1,2,3$, and the stabilizer subgroup of $M^i$ in $H$ is given by
		$G_{M^i}^{H}
		=
		G_{M^i}\cap H$.
		By the induction hypothesis, we obtain
		\begin{eqnarray*}
			N^{M^3}_{M^1,M^2}
			=
			\sum_{W\in\mathcal{M}(1;g_1,g_2)}
			\frac{
				S_{M^1,W}S_{M^2,W}\overline{S_{M^3,W}}
			}
			{S_{V,W}} .
		\end{eqnarray*}
		Since $g_1,g_2,g_3\in H$, this is precisely the formula required
		for these twisted sectors in the $G$-orbifold.

		\textbf{\emph{Case 2.}} Suppose that $\langle g_1,g_2\rangle= G.$
		Then $G
		=
		\langle g_1,g_2\rangle=
		\langle g_2,g_3\rangle=
		\langle g_1,g_3\rangle$.
		Since $g_i\in G_{M^i}$ for $i=1,2,3$, we have $G_{M^1}G_{M^2}
		=
		G_{M^2}G_{M^3}
		=
		G_{M^1}G_{M^3}
		=
		G$.
		By Theorem~\ref{thm 4.3}, we have
		\begin{eqnarray*}
			\sum_{U\in\mathcal{O}_{M^1}}\sum_{W\in\mathcal{O}_{M^2}}N^{M^3}_{U,W}
			&=&
			\sum_{\chi_3\in\operatorname{Irr}(G_{M^3})}
			N^{M^3_{\chi_3}}_{M^1_{\chi_1},M^2_{\chi_2}},
		\end{eqnarray*}
		where $\chi_1\in\operatorname{Irr}(G_{M^1})$ and $\chi_2\in\operatorname{Irr}(G_{M^2})$.
		Since $G_{M^1}G_{M^2}
		=
		G_{M^2}G_{M^3}
		=
		G_{M^1}G_{M^3}
		=
		G$, Proposition~\ref{pp 4.1} implies that $N^{M^3}_{U,W}
		=
		N^{M^3}_{M^1,M^2}$
		for any $U\in\mathcal{O}_{M^1}, W\in\mathcal{O}_{M^2}$.
		Therefore,
		\begin{eqnarray*}
			N^{M^3}_{M^1,M^2}
			&=&
			\frac{1}{[G:G_{M^1}]}\frac{1}{[G:G_{M^2}]}
			\sum_{\chi_3\in\operatorname{Irr}(G_{M^3})}
			N^{M^3_{\chi_3}}_{M^1_{\chi_1},M^2_{\chi_2}}.
		\end{eqnarray*}
		Applying Theorem~\ref{thm 25}, we obtain
		\begin{eqnarray*}
			N^{M^3_{\chi_3}}_{M^1_{\chi_1},M^2_{\chi_2}}
			&=&
			\sum_{N\in\mathcal{I}(V^G)}
			\frac{
				S_{M^1_{\chi_1},N}
				S_{M^2_{\chi_2},N}
				\overline{S_{M^3_{\chi_3},N}}
			}
			{S_{V^G,N}} .
		\end{eqnarray*}
		By Theorem~\ref{thm 3.2}, every irreducible $V^G$-module $N$ is of the
		form $M^{g,H,j}_{\chi}$. For simplicity, we write $N=W_{\chi}$, where
		$W$ is an irreducible $g_W$-twisted $V$-module.
		By Theorem~\ref{thm 3.4} and Lemma~\ref{lem 5.3}, we have
		\begin{eqnarray*}
			&&S_{M^i_{\chi_i},W_{\chi}}
			=
			\begin{cases}
				\dfrac{
					\overline{\chi_i(g_W)}\chi(g_i^{-1})
				}{|G_{M^i}|}
				\displaystyle\sum_{k_i\in C_{M^i,W}}S_{M^i,W\circ k_i},
				& g_W\in G_{M^i},\ g_i\in G_{W}\\[1.2ex]
				0,
				& \text{otherwise}
			\end{cases}
			\quad \text{for}\,i=1,2,3,\\
			&&S_{V^G,W_{\chi}}=\frac{1}{|G_W|}S_{V,W}.
		\end{eqnarray*}
		Since $G_{M^3}$ is a finite abelian group, we have
		\begin{eqnarray*}
			\sum_{\chi_3\in\operatorname{Irr}(G_{M^3})}
			\chi_3(g_W)
			=
			\begin{cases}
				|G_{M^3}|, & g_W=1,\\
				0, & g_W\neq1.
			\end{cases}
		\end{eqnarray*}
		Moreover, we have $\sum_{\chi\in\operatorname{Irr}(G_W)}
		{\chi(g_1^{-1})}
		{\chi(g_2^{-1})}
		\overline{\chi(g_3^{-1})}
		=
		|G_W|.$  Thus, we obtain
		\begin{eqnarray*}
			N^{M^3}_{M^1,M^2}
			&=&
			\frac{|G_{M^1}||G_{M^2}|}{|G|^2}
			\sum_{\chi_3\in\operatorname{Irr}(G_{M^3})}
			\sum_{g_W\in G}
			\sum_{W\in\mathcal{R}(g_W)}
			\sum_{\chi\in\operatorname{Irr}(G_W)}
			\frac{
				S_{M^1_{\chi_1},W_\chi}
				S_{M^2_{\chi_2},W_\chi}
				\overline{S_{M^3_{\chi_3},W_\chi}}
			}
			{S_{V^G,W_\chi}}
			\\
			&=&
			\frac{|G_{M^1}||G_{M^2}|}{|G|^2}
			\sum_{\chi_3\in\operatorname{Irr}(G_{M^3})}
			\sum_{g_W\in G_{M^1}\cap G_{M^2}\cap G_{M^3}}
			\sum_{W\in\mathcal{R}(g_W)}
			\sum_{g_1,g_2\in G_W\leq G}
			\sum_{\chi\in\operatorname{Irr}(G_W)}
			\\
			&&\cdot
			\frac{
				\overline{\chi_1(g_W)}\chi(g_1^{-1})\overline{\chi_2(g_W)}\chi(g_2^{-1})	\chi_3(g_W)\overline{\chi(g_3^{-1})}
			}
			{|G_{M^1}||G_{M^2}||G_{M^3}|}
			\frac{
				|G_W|
			}
			{S_{V,W}}
			\\
			&&\cdot\sum_{k_1\in C_{M^1,W}}S_{M^1,W\circ k_1}\sum_{k_2\in C_{M^2,W}}S_{M^2,W\circ k_2}\sum_{k_3\in C_{M^3,W}}\overline {S_{M^3,W\circ k_3}}\\
			&=&
			\frac{1}{|G|}\sum_{\chi_3\in\operatorname{Irr}(G_{M^3})}\sum_{W\in\mathcal{R}_G(1)}\sum_{\chi\in\operatorname{Irr}(G)}\frac{
				\chi(g_1^{-1})\chi(g_2^{-1})\chi_3{(1)}\overline{\chi(g_3^{-1})}
			}
			{|G_{M^3}|}\frac{
				S_{M^1,W}S_{M^2,W}\overline {S_{M^3,W}}
			}
			{S_{V,W}}
			\\
			&=&
			\sum_{W\in\mathcal{M}(1;g_1,g_2)}
			\frac{
				S_{M^1,W}
				S_{M^2,W}
				\overline{S_{M^3,W}}
			}
			{S_{V,W}}.
		\end{eqnarray*}
		
		The proof is concluded.
	\end{proof}
	
	\section{Example: The Permutation Orbifold
		$\left(V^{\otimes k}\right)^{\langle (1\,2\,\cdots\,k)\rangle}$}
	In this section, we apply the twisted Verlinde formula established in
	Section~5 to permutation orbifolds $\left(V^{\otimes k}\right)^{\langle g\rangle}$ associated with a $k$-cycle
	$g=(1\,2\,\cdots\,k)$. We first recall the construction of permutation twisted
	modules and describe the irreducible twisted and untwisted modules that are
	stable under prescribed powers of $g$. Using the corresponding $S$-matrix
	relations, we derive explicit formulas for the fusion coefficients between
	$g^s$-twisted and $g^r$-twisted modules, as well as those between untwisted
	and $g^s$-twisted modules, for arbitrary $1\leq r,s\leq k-1$. By repeatedly applying the twisted Verlinde formula,
	we also obtain an explicit formula for the  fusion product
	$\left(T_g^0\right)^{\boxtimes k}$. We then
	specialize these formulas to the rank-one lattice vertex operator algebra
	$V_{\mathbb{Z}\alpha}$ with
	$\langle\alpha,\alpha\rangle=2n$, where the fusion coefficients are evaluated
	explicitly by finite Fourier sums. In the case $k=4$ and $n=1$, the resulting fusion
	rules agree with those obtained in Section~6 of \cite{dxy2}.

	We first recall several results on twisted modules for permutation
	orbifolds, including the construction of irreducible twisted and
	untwisted modules, following \cite{bdm} (cf. \cite{dxy}, \cite{dxy2}).

	Let $V=(V,Y,\mathbf{1},\omega)$ be a vertex operator algebra and $k$ be a positive integer. In
	$(\operatorname{End}V)[[x^{\frac{1}{k}},x^{-\frac{1}{k}}]]$, define
	\begin{eqnarray*}
		\Delta_k(x)
		&=&
		\exp\left(
		\sum_{j\in\mathbb{Z}_{+}}
		a_j x^{-\frac{j}{k}}L(j)
		\right)
		k^{-L(0)}
		x^{(\frac{1}{k}-1)L(0)},
	\end{eqnarray*}
	where the coefficients $a_j$ for $j\in\mathbb{Z}_{+}$ are uniquely determined by
	\begin{eqnarray*}
		\exp\left(
		\sum_{j\in\mathbb{Z}_{+}}
		-a_j x^{j+1}\frac{d}{dx}
		\right)x
		=
		\frac{1}{k}(1+x)^k-\frac{1}{k}.
	\end{eqnarray*}

	The following theorem, proved in \cite{bdm}, describes the
	construction of $g$-twisted $V^{\otimes k}$-modules associated with
	the cyclic permutation $g=(1\,2\,\cdots\,k)$ by means of the operator
	$\Delta_k(z)$ (see also \cite{dxy}, \cite{dxy2}).
	\bt{thm 6.1}
	Let $V$ be any vertex operator algebra and let
	$g=(1\,2\,\cdots\,k)\in\operatorname{Aut}(V^{\otimes k})$.
	Then for any (weak, admissible) $V$-module $(W,Y_W)$, we have a
	(weak, admissible) $g$-twisted $V^{\otimes k}$-module
	$(T_g(W),Y_{T_g(W)})$, where $T_g(W)=W$ as a vector space and the
	vertex operator map $Y_{T_g(W)}(\cdot,x)$ is uniquely determined by
	\begin{eqnarray}
		Y_{T_g(W)}(u^1,x)
		=
		Y_W\left(\Delta_k(x)u,x^{\frac{1}{k}}\right)
		\quad \text{for}\,\,u\in V ,
	\end{eqnarray}
	where $u^1
	=
	u\otimes\mathbf{1}^{\otimes(k-1)}
	\in V^{\otimes k}.$
	Furthermore, every (weak, admissible) $g$-twisted
	$V^{\otimes k}$-module is isomorphic to one of this form.
	\et
	
	\br{remark 6.2}{\em
		For $k=1$, we have $\Delta_1(x)=\operatorname{Id}_V$. Hence, by
		Theorem~\ref{thm 6.1}, we have
		\begin{eqnarray*}
			Y_{T_1(W)}(u,x)
			&=&
			Y_W\left(\Delta_1(x)u,x\right)
			=
			Y_W(u,x)
			\quad \text{for}\,\,u\in V.
		\end{eqnarray*}
		Thus, $T_1(W)$ has the same $V$-module structure as $W$, and hence
		$T_1(W)\cong W$.
	}\er
	
	Let $\operatorname{Irr}(V)
	=
	\{V=M^0,M^1,\ldots,M^p\}$
	be the set of irreducible $V$-modules up to isomorphism. For
	$1\leq s<k$, let $d=\gcd(s,k)$. Then the order of $g^s$ is
	$\frac{k}{d}$, and the cycle decomposition of $g^s$ is given by
	\begin{eqnarray*}
		g^s=\tau_1\tau_2\cdots\tau_d,
	\end{eqnarray*}
	where each $\tau_i$ is a $\frac{k}{d}$-cycle for $1\leq i\leq d$.
	By \cite{bdm}, every irreducible $g^s$-twisted
	$V^{\otimes k}$-module can be constructed in the form
	\begin{eqnarray}
		T_{\tau_1,\tau_2,\ldots,\tau_d}^{i_1,i_2,\ldots,i_d}
		=
		T_{\tau_1}(M^{i_1})\otimes T_{\tau_2}(M^{i_2})\otimes\cdots\otimes
		T_{\tau_d}(M^{i_d}),
		\label{eq:twisted-module-decomposition}
	\end{eqnarray}
	where $i_1,i_2,\ldots,i_d\in\{0,1,\ldots,p\}$.
	We shall also denote this module by $T_{g^s}^{i_1,i_2,\ldots,i_d}$.
	
	For later use, we describe the relation between $g^s$-twisted modules
	and $g^r$-stable modules (\cite{dxy}, \cite{dxy2}). Let $1\leq r,s<k$, and set
	\begin{eqnarray*}
		d=\gcd(s,k),\,
		{t}=\gcd(d,r),\,
		l=\frac{k}{d},\,
		c=\frac{d}{t}.
	\end{eqnarray*}
	The set of equivalence classes of irreducible $g^s$-twisted
	$V^{\otimes k}$-modules that are $g^r$-stable is denoted by
	$\mathcal{M}(g^s,g^r)$, and it is given by
	\begin{eqnarray}
		\mathcal{M}(g^s,g^r)
		=
		\left\{
		\left(T_{\alpha_1}(W^1)\right)^{\otimes c}
		\otimes\cdots\otimes
		\left(T_{\alpha_{t}}(W^{t})\right)^{\otimes c}
		\ \middle|\
		W^1,\ldots,W^{t}\in\operatorname{Irr}(V)
		\right\},
		\label{eq:L-gs-gr}
	\end{eqnarray}
	where each $\alpha_i$ is an $l$-cycle for $1\leq i\leq {t}$.
	We further write
	\begin{eqnarray}
		T_{g^s}^{W^1,W^2,\ldots,W^{t};c}
		=
		\left(T_{\alpha_1}(W^1)\right)^{\otimes c}
		\otimes \left(T_{\alpha_2}(W^2)\right)^{\otimes c}
		\otimes\cdots\otimes
		\left(T_{\alpha_{t}}(W^{t})\right)^{\otimes c}.
		\label{eq:notation-Tgs}
	\end{eqnarray}

	When $c=1$, we have $d={t}$. In this case, the
	above module agrees with the notation
	$T_{g^s}^{j_1,j_2,\ldots,j_d}$ introduced previously by taking
	$W^i=M^{j_i}$ and $\alpha_i=\tau_i$.

	We also describe the modules in $\mathcal{M}(g^r,g^s)$.
	Let $d_1=\gcd(r,k)$. Then there exist $m,n\in\mathbb{Z}$ such that
	$rm+kn=d_1$. Moreover, $\gcd(d_1,s)
	=
	\gcd(d,r)
	=
	\gcd(s,k,r)
	=
	{t}.$
	Set $l_1=\frac{k}{d_1},
	$ and $
	a=\frac{d_1}{{t}}$.
	Then every irreducible $g^r$-twisted module in
	$\mathcal{M}(g^r,g^s)$ has the form
	\begin{eqnarray}
		T_{g^r}^{W^1,W^2,\ldots,W^{t};a}
		=
		\left(T_{\beta_1}(W^1)\right)^{\otimes a}
		\otimes \left(T_{\beta_2}(W^2)\right)^{\otimes a}
		\otimes\cdots\otimes
		\left(T_{\beta_{t}}(W^{t})\right)^{\otimes a},
		\label{eq:gr-gs-module}
	\end{eqnarray}
	where $W^i\in\operatorname{Irr}(V)$ and each $\beta_i$ is an
	$l_1$-cycle for $1\leq i\leq {t}$.
	
	In particular, for $a=1$, the above notation reduces to the previous
	notation $T_{g^r}^{j_1,\ldots,j_m}$, where
	$m=\gcd(r,k)$ and the corresponding irreducible $V$-modules are chosen
	appropriately.
	
	We next recall the modules in $\mathcal{M}(1,g^s)$. By the results
	in \cite{dxy} and \cite{dxy2}, every $V^{\otimes k}$-module in $\mathcal{M}(1,g^s)$ is
	of the form
	\begin{eqnarray}
		\left(M^{i_1}\right)^{\otimes l}
		\otimes \left(M^{i_2}\right)^{\otimes l}
		\otimes\cdots\otimes
		\left(M^{i_d}\right)^{\otimes l},
		\label{6.6}
	\end{eqnarray}
	where $d=\gcd(s,k)$, $k=dl$, and $ i_1,i_2, \dots, i_d \in \{0, 1, \dots, p\}$.

	By a direct computation, we have the following result.
	\bl{lem 6}
	For
	$r,s\in\mathbb{Z}$, set $t=\gcd(s,r,k)$ and $l=\frac{k}{t}$.
	Then
	\begin{eqnarray*}
		\mathcal{M}(1;g^s,g^r)
		=
		\left\{
		(M^{a_1})^{\otimes l}
		\otimes      (M^{a_2})^{\otimes l}
		\otimes      \cdots\otimes
		(M^{a_{t}})^{\otimes l}
		\ \middle|\
		0\leq a_1,a_2,\ldots,a_{t}\leq p
		\right\},
	\end{eqnarray*}
	where $\mathcal{M}(1;g^s,g^r)$ denotes the set of isomorphism classes of
	irreducible $V^{\otimes k}$-modules that are simultaneously
	$g^s$-stable and $g^r$-stable.
	\el

	\begin{proof}
		It suffices to show that $\mathcal{M}(1;g^s,g^r)=\mathcal{M}(1,g^{t}).$
		Since ${t}=\gcd(r,s,k)$, we have
		\begin{eqnarray*}
			\langle g^r,g^s\rangle=\langle g^{t}\rangle.
		\end{eqnarray*}
		Thus, for any irreducible $V^{\otimes k}$-module $M$, we obtain
		\begin{eqnarray*}
			g^r,g^s\in G_M
			\quad\Longleftrightarrow\quad
			g^{t}\in G_M,
		\end{eqnarray*}
		where $G_M$ denotes the stabilizer subgroup of $M$ in $\langle g\rangle$.
		The desired equality follows immediately.
	\end{proof}

	The following lemma plays an important role in the explicit computation
	of fusion rules for permutation orbifolds (see the proof of Lemma 4.5 of \cite{dxy}).

	\bl{lem 6.3}
	Let $g=(1\,2\,\cdots\,k)$
	and fix a nonnegative integer $s$. Set $d=\gcd(s,k)$ and
	$l=\frac{k}{d}.$
	For $	i_1,i_2,\ldots,i_d,j_1,j_2,\ldots,j_d
	\in
	\{0,1,\ldots,p\},$
	let
	\begin{eqnarray*}
		X^{i_1,i_2,\ldots,i_d}
		=
		(M^{i_1})^{\otimes l}
		\otimes (M^{i_2})^{\otimes l}
		\otimes\cdots\otimes
		(M^{i_d})^{\otimes l} \in \mathcal{M}(1,g^{-s})
	\end{eqnarray*}
	and $T_{g^s}^{{j_1},{j_2},\ldots,{j_d}}$ be a $g^s$-twisted
	$V^{\otimes k}$-module as above.
	Then
	\begin{eqnarray*}
		S_{
			T_{g^s}^{{j_1},{j_2},\ldots,{j_d}},
			X^{i_1,i_2,\ldots,i_d}
		}
		=
		S_{M^{j_1},M^{i_1}}S_{M^{j_2},M^{i_2}}
		\cdots
		S_{M^{j_d},M^{i_d}}.
	\end{eqnarray*}
	\el

	We now apply the twisted Verlinde formula to compute the fusion rules
	for twisted modules of
	$\left(V^{\otimes k}\right)^{\langle g\rangle}$
	with $g=(1\,2\,\cdots\,k)$.

	\bt{thm 6.5}
	Assume that $V$ is a simple, rational, $C_2$-cofinite and self-dual vertex
	operator algebra of CFT type.  Let $g=(1\,2\,\cdots\,k)$ be a $k$-cycle, and let $1\leq s,r\leq k-1$. Suppose that
	\begin{eqnarray*}
		g^s=\tau_1\tau_2\cdots\tau_d,\quad
		g^r=\sigma_1\sigma_2\cdots\sigma_m,\quad
		g^{s+r}=\gamma_1\gamma_2\cdots\gamma_q
	\end{eqnarray*}
	are products of disjoint cycles, where $d=\gcd(s,k),
	m=\gcd(r,k)$ and
	$q=\gcd(s+r,k).$
	Write $S_{i,j}
	:=
	S_{M^i,M^j}$ for $0\leq i,j\leq p$  and ${t}=\gcd(s,r,k)$.
	Then, for any $i_1,i_2,\ldots,i_d,\,
	j_1,j_2,\ldots,j_m,\,
	t_1,t_2,\ldots,t_q
	\in\{0,1,\ldots,p\},$ we have
	\begin{eqnarray}
		N^{\,T_{g^{s+r}}^{t_1,\ldots,t_q}}_{
			T_{g^s}^{i_1,\ldots,i_d},
			T_{g^r}^{j_1,\ldots,j_m}}
		=
		\prod_{\nu=1}^{{t}}
		\left(
		\sum_{a=0}^{p}
		\frac{
			\displaystyle
			\prod_{u=1}^{\frac{d}{t}}
			S_{i_{(\nu-1)\frac{d}{t}+u},a}
			\prod_{v=1}^{\frac{m}{t}}
			S_{j_{(\nu-1)\frac{m}{t}+v},a}
			\prod_{w=1}^{\frac{q}{t}}
			\overline{S_{t_{(\nu-1)\frac{q}{t}+w},a}}
		}{
			\displaystyle
			S_{0,a}^{\frac{k}{t}}
		}
		\right). \label{6.51}
	\end{eqnarray}
	\et
	\begin{proof}
		For any $i_1,i_2,\ldots,i_d,\,
		j_1,j_2,\ldots,j_m,\,
		t_1,t_2,\ldots,t_q
		\in\{0,1,\ldots,p\},$ Theorem \ref{thm 5.1} gives
		\begin{eqnarray}
			N^{\,T_{g^{s+r}}^{t_1,\ldots,t_q}}_{
				T_{g^s}^{i_1,\ldots,i_d},
				T_{g^r}^{j_1,\ldots,j_m}}
			=
			\sum_{W\in\mathcal{M}(1;g^s,g^r)}
			\frac{
				S_{T_{g^s}^{i_1,\ldots,i_d},W}
				S_{T_{g^r}^{j_1,\ldots,j_m},W}
				\overline{
					S_{T_{g^{s+r}}^{t_1,\ldots,t_q},W}
				}
			}{
				S_{V^{\otimes k},W}
			}.\label{6.5}
		\end{eqnarray}
		Set ${t}=\gcd(s,r,k)$,
		by Lemma \ref{lem 6}, we have
		\begin{eqnarray*}
			\mathcal{M}(1;g^s,g^r)
			=
			\left\{
			(M^{a_1})^{\otimes \frac{k}{t}}
			\otimes (M^{a_2})^{\otimes \frac{k}{t}}
			\otimes\cdots\otimes
			(M^{a_{t}})^{\otimes \frac{k}{t}}
			\ \middle|\
			0\leq a_1,a_2,\ldots,a_{t}\leq p
			\right\}.
		\end{eqnarray*}
		Since ${t}$ divides $d$, for each $1\leq\nu\leq {t}$ we have
		\begin{eqnarray*}
			(M^{a_\nu})^{\otimes \frac{k}{t}}
			=
			\bigotimes_{u=1}^{\frac{d}{t}}
			(M^{a_\nu})^{\otimes \frac{k}{d}}.
		\end{eqnarray*}
		Consequently, every
		$W\in\mathcal{M}(1;g^s,g^r)$ may be written in the form
		\begin{eqnarray*}
			W
			=
			\bigotimes_{\nu=1}^{{t}}
			\bigotimes_{u=1}^{\frac{d}{t}}
			(M^{a_\nu})^{\otimes \frac{k}{d}}.
		\end{eqnarray*}
		It follows from Lemma \ref{lem 6.3} and Remark~\ref{remark 6.2} that
		\begin{eqnarray*}
			S_{T_{g^s}^{i_1,\ldots,i_d},
				\,(M^{a_1})^{\otimes \frac{k}{t}}\otimes\cdots
				\otimes(M^{a_{t}})^{\otimes \frac{k}{t}}}
			&=&
			\prod_{\nu=1}^{{t}}
			\prod_{u=1}^{\frac{d}{t}}
			S_{i_{(\nu-1)\frac{d}{t}+u},a_\nu},
			\\[6pt]
			S_{T_{g^r}^{j_1,\ldots,j_m},
				\,(M^{a_1})^{\otimes \frac{k}{t}}\otimes\cdots
				\otimes(M^{a_{t}})^{\otimes \frac{k}{t}}}
			&=&
			\prod_{\nu=1}^{{t}}
			\prod_{v=1}^{\frac{m}{t}}
			S_{j_{(\nu-1)\frac{m}{t}+v},a_\nu},
			\\[6pt]
			\overline{
				S_{T_{g^{s+r}}^{t_1,\ldots,t_q},
					\,(M^{a_1})^{\otimes \frac{k}{t}}\otimes\cdots
					\otimes(M^{a_{t}})^{\otimes \frac{k}{t}}}
			}
			&=&
			\prod_{\nu=1}^{{t}}
			\prod_{w=1}^{\frac{q}{t}}
			\overline{
				S_{t_{(\nu-1)\frac{q}{t}+w},a_\nu}
			},
			\\[6pt]
			S_{V^{\otimes k},
				\,(M^{a_1})^{\otimes \frac{k}{t}}\otimes\cdots
				\otimes(M^{a_{t}})^{\otimes \frac{k}{t}}}
			&=&
			\prod_{\nu=1}^{{t}}S_{0,a_\nu}^{\,\frac{k}{t}}.
		\end{eqnarray*}
		Substituting these identities into (\ref{6.5}) gives
		\begin{eqnarray*}
			N^{\,T_{g^{s+r}}^{t_1,\ldots,t_q}}_{
				T_{g^s}^{i_1,\ldots,i_d},
				T_{g^r}^{j_1,\ldots,j_m}}
			&=&
			\sum_{a_1,\ldots,a_{t}=0}^{p}
			\prod_{\nu=1}^{{t}}
			\frac{
				\displaystyle
				\prod_{u=1}^{\frac{d}{t}}
				S_{i_{(\nu-1)\frac{d}{t}+u},a_\nu}
				\prod_{v=1}^{\frac{m}{t}}
				S_{j_{(\nu-1)\frac{m}{t}+v},a_\nu}
				\prod_{w=1}^{\frac{q}{t}}
				\overline{S_{t_{(\nu-1)\frac{q}{t}+w},a_\nu}}
			}{
				\displaystyle
				S_{0,a_\nu}^{\,\frac{k}{t}}
			}
			\\[6pt]
			&=&
			\prod_{\nu=1}^{{t}}
			\left(
			\sum_{a=0}^{p}
			\frac{
				\displaystyle
				\prod_{u=1}^{\frac{d}{t}}
				S_{i_{(\nu-1)\frac{d}{t}+u},a}
				\prod_{v=1}^{\frac{m}{t}}
				S_{j_{(\nu-1)\frac{m}{t}+v},a}
				\prod_{w=1}^{\frac{q}{t}}
				\overline{S_{t_{(\nu-1)\frac{q}{t}+w},a}}
			}{
				\displaystyle
				S_{0,a}^{\,\frac{k}{t}}
			}
			\right).
		\end{eqnarray*}
		
		This completes the proof.
	\end{proof}

	\br{remark 6.6}{\em
		(1) If $\gcd(s,r)=1$, then ${t}=\gcd(s,r,k)=1.$
		Therefore,  formula (\ref{6.51}) reduces to
		\begin{eqnarray*}
			N^{\,T_{g^{s+r}}^{t_1,\ldots,t_q}}_{
				T_{g^s}^{i_1,\ldots,i_d},
				T_{g^r}^{j_1,\ldots,j_m}}
			=
			\sum_{a=0}^{p}
			\frac{
				\displaystyle
				S_{i_1,a}\cdots S_{i_d,a}
				S_{j_1,a}\cdots S_{j_m,a}
				\overline{S_{t_1,a}}\cdots
				\overline{S_{t_q,a}}
			}{
				\displaystyle
				S_{0,a}^{\,k}
			},
		\end{eqnarray*}
		which coincides with the formula in Theorem~4.7 of \cite{dxy2}.
		
		\noindent(2) For $\gcd(r,s)\neq 1$, the fusion product in \cite{dxy2}
		is reduced to fusion products associated with the disjoint cycles of
		$g^{\gcd(r,s)}$. The formula in Theorem \ref{thm 6.5} gives a direct and uniform expression for
		arbitrary $r$ and $s$ with respect to the fixed $k$-cycle $g$.
	}\er

	As a consequence of Theorem~\ref{thm 6.5}, the following result follows from repeated applications of the twisted Verlinde formula for cyclic groups. It agrees with Theorem 5.1 of \cite{dxy2} and provides an alternative proof.
	\bc{cor 6}
	Let $g=(1\,2\,\cdots\,k)$ be a $k$-cycle.
	For convenience, write
	\begin{eqnarray*}
		M^{i_1,i_2,\ldots,i_k}
		=
		M^{i_1}\otimes M^{i_2}\otimes\cdots\otimes M^{i_k},
	\end{eqnarray*}
	where $i_1,i_2,\ldots,i_k\in\{0,1,\ldots,p\}$.
	Then
	\begin{eqnarray*}
		\left(T_g^0\right)^{\boxtimes k}
		&=&
		\sum_{i_1,i_2,\ldots,i_k=0}^{p}
		\left(
		\sum_{a=0}^{p}
		S_{0,a}^{\,2-(k-1)(k-2)}
		\prod_{\ell=1}^{k}
		\frac{\overline{S_{i_{\ell},a}}}{S_{0,a}}
		\right)
		M^{i_1,i_2,\ldots,i_k}.
	\end{eqnarray*}
	\ec
	\begin{proof}
		The  fusion product is written
		as
		\begin{eqnarray*}
			\left(T_g^0\right)^{\boxtimes k}
			=
			\left(\cdots\left(
			\left(T_g^0\boxtimes T_g^0\right)
			\boxtimes T_g^0
			\right)\cdots\right)\boxtimes T_g^0.
		\end{eqnarray*}
		The intermediate fusion products belong successively to the
		$g^2$-, $g^3$-, $\ldots$, $g^{k-1}$-twisted sectors, and the final
		fusion product belongs to the untwisted sector.
		
		We first compute the fusion product $T_g^0\boxtimes T_g^0.$
		Let $T_{g^2}^{j_1,j_2,\ldots,j_{d_2}}$
		be an irreducible $g^2$-twisted module, where $d_2=\gcd(2,k)$ and $j_1,j_2,\ldots,j_{d_2}\in\{0,1,\ldots,p\}.$ By Theorem~\ref{thm 6.5}, we have
		\begin{eqnarray*}
			N_{T_g^0,T_g^0}^{T_{g^2}^{j_1,j_2,\ldots,j_{d_2}}}
			=
			\sum_{a=0}^{p}
			\frac{
				S_{0,a}S_{0,a}
				\displaystyle\prod_{u=1}^{d_2}
				\overline{S_{j_u,a}}
			}{
				S_{0,a}^{k}
			}
			=
			\sum_{a=0}^{p}
			S_{0,a}^{\,2-k}
			\prod_{u=1}^{d_2}
			\overline{S_{j_u,a}}.
		\end{eqnarray*}
		Therefore, the multiplicity of
		$T_{g^2}^{j_1,j_2,\ldots,j_{d_2}}$ in
		$T_g^0\boxtimes T_g^0$ is $\sum_{a=0}^{p}
		S_{0,a}^{\,2-k}
		\prod_{u=1}^{d_2}
		\overline{S_{j_u,a}}.$
		
		We next  compute $T_{g^2}^{j_1,j_2,\ldots,j_{d_2}}\boxtimes T_g^0.$
		Let $T_{g^3}^{t_1,t_2,\ldots,t_{d_3}}$
		be an irreducible $g^3$-twisted module, where $d_3=\gcd(3,k)$ and $t_1,t_2,\ldots,t_{d_3}\in\{0,1,\ldots,p\}.$ Applying Theorem~\ref{thm 6.5} again, we
		obtain
		\begin{eqnarray*}
			N_{T_{g^2}^{j_1,j_2,\ldots,j_{d_2}},T_g^0}
			^{T_{g^3}^{t_1,t_2,\ldots,t_{d_3}}}
			&=&
			\sum_{b=0}^{p}
			\frac{
				\left(
				\displaystyle\prod_{u=1}^{d_2}S_{j_u,b}
				\right)
				S_{0,b}
				\left(
				\displaystyle\prod_{v=1}^{d_3}
				\overline{S_{t_v,b}}
				\right)
			}{
				S_{0,b}^{k}
			}
			\\
			&=&
			\sum_{b=0}^{p}
			S_{0,b}^{\,1-k}
			\left(
			\prod_{u=1}^{d_2}S_{j_u,b}
			\right)
			\left(
			\prod_{v=1}^{d_3}
			\overline{S_{t_v,b}}
			\right).
		\end{eqnarray*}
		It follows that the multiplicity of
		$T_{g^3}^{t_1,t_2,\ldots,t_{d_3}}$ in
		$\left(T_g^0\right)^{\boxtimes 3}$ is
		\begin{eqnarray}
			&&
			\sum_{j_1,j_2,\ldots,j_{d_2}=0}^{p}
			\left(
			\sum_{a=0}^{p}
			S_{0,a}^{\,2-k}
			\prod_{u=1}^{d_2}
			\overline{S_{j_u,a}}
			\right)
			\left(
			\sum_{b=0}^{p}
			S_{0,b}^{\,1-k}
			\prod_{u=1}^{d_2}
			S_{j_u,b}
			\prod_{v=1}^{d_3}
			\overline{S_{t_v,b}}
			\right).\label{1a}
		\end{eqnarray}
		Since the
		$S$-matrix is unitary,
		we obtain
		\begin{eqnarray*}
			\sum_{j_1,j_2,\ldots,j_{d_2}=0}^{p}
			\prod_{u=1}^{d_2}
			\overline{S_{j_u,a}}S_{j_u,b}
			=
			\prod_{u=1}^{d_2}
			\left(
			\sum_{j_u=0}^{p}
			\overline{S_{j_u,a}}S_{j_u,b}
			\right)
			=
			\delta_{a,b},
		\end{eqnarray*}
		where $\delta_{a,b}$ denotes the Kronecker delta.
		Consequently, the above
		multiplicity (\ref{1a}) reduces to
		\begin{eqnarray*}
			\sum_{a=0}^{p}
			S_{0,a}^{\,2-k}
			S_{0,a}^{\,1-k}
			\prod_{v=1}^{d_3}
			\overline{S_{t_v,a}}
			&=&
			\sum_{a=0}^{p}
			S_{0,a}^{\,3-2k}
			\prod_{v=1}^{d_3}
			\overline{S_{t_v,a}}.
		\end{eqnarray*}
		Repeating the above argument, let
		$T_{g^{k-1}}^{q_1,\ldots,q_{d_{k-1}}}$ be an irreducible
		$g^{k-1}$-twisted module, where $d_{k-1}=\gcd(k-1,k)=1$
		and
		$q_1,q_2,\ldots,q_{d_{k-1}}\in\{0,1,\ldots,p\}$.
		The multiplicity of
		$T_{g^{k-1}}^{q_1,q_2,\ldots,q_{d_{k-1}}}$ in
		$\left(T_g^0\right)^{\boxtimes(k-1)}$ is
		\begin{eqnarray*}
			\sum_{c=0}^{p}
			S_{0,c}^{\,k-1-(k-2)k}
			\prod_{u=1}^{d_{k-1}}
			\overline{S_{q_u,c}}.
		\end{eqnarray*}
		For $i_1,i_2,\ldots,i_k\in\{0,1,\ldots,p\}$, set $M^{i_1,\ldots,i_k}
		=
		M^{i_1}\otimes\cdots\otimes M^{i_k}.$
		By Theorem~\ref{thm 6.5}, we have
		\begin{eqnarray*}
			N_{T_{g^{k-1}}^{q_1,q_2,\ldots,q_{d_{k-1}}},T_g^0}
			^{M^{i_1,i_2,\ldots,i_k}}
			&=&
			\sum_{b=0}^{p}
			S_{0,b}^{\,1-k}
			\left(
			\prod_{u=1}^{d_{k-1}}S_{q_u,b}
			\right)
			\left(
			\prod_{\ell=1}^{k}
			\overline{S_{i_\ell,b}}
			\right).
		\end{eqnarray*}
		Combining these two formulas, summing over
		$q_1,q_2,\ldots,q_{d_{k-1}}$, and using the unitarity of the
		$S$-matrix, we find that the multiplicity of
		$M^{i_1,i_2,\ldots,i_k}$ in $\left(T_g^0\right)^{\boxtimes k}$ is
		\begin{eqnarray*}
			\sum_{a=0}^{p}
			S_{0,a}^{\,k-1-(k-2)k}
			S_{0,a}^{\,1-k}
			\prod_{\ell=1}^{k}
			\overline{S_{i_\ell,a}}
			=
			\sum_{a=0}^{p}
			S_{0,a}^{\,2-(k-1)(k-2)}
			\prod_{\ell=1}^{k}
			\frac{\overline{S_{i_\ell,a}}}{S_{0,a}}.
		\end{eqnarray*}
		Therefore, we obtain
		\begin{eqnarray*}
			\left(T_g^0\right)^{\boxtimes k}
			&=&
			\sum_{i_1,i_2,\ldots,i_k=0}^{p}
			\left(
			\sum_{a=0}^{p}
			S_{0,a}^{\,2-(k-1)(k-2)}
			\prod_{\ell=1}^{k}
			\frac{\overline{S_{i_\ell,a}}}{S_{0,a}}
			\right)
			M^{i_1,i_2,\ldots,i_k}.
		\end{eqnarray*}
	\end{proof}

	We next consider the fusion rules between untwisted and twisted modules for
	$\left(V^{\otimes k}\right)^{\langle g\rangle}$
	with $g=(1\,2\,\cdots\,k)$.

	\bt{thm 6.7}
	Assume that $V$ is a simple, rational, $C_2$-cofinite and self-dual vertex
	operator algebra of CFT type.  Let $g=(1\,2\,\cdots\,k)$ be a $k$-cycle, and let $1\leq s\leq k-1$. Suppose that $g^s=\tau_1\tau_2\cdots\tau_d,$
	is products of disjoint cycles, where $d=\gcd(s,k)$ and $k=dl.$
	Write $S_{i,j}
	:=
	S_{M^i,M^j}$ for $0\leq i,j\leq p$.
	Then, for any $i_1,\ldots,i_k,\,
	j_1,\ldots,j_d,\,
	t_1,\ldots,t_d
	\in\{0,1,\ldots,p\},$ we have
	\begin{eqnarray*}
		N^{\,T_{g^s}^{t_1,t_2,\ldots,t_d}}_{
			M^{i_1}\otimes M^{i_2}\otimes\cdots\otimes M^{i_{k}},
			T_{g^s}^{j_1,j_2,\ldots,j_d}}
		&=&
		\prod_{\nu=1}^{d}
		\left(
		\sum_{q=0}^{p}
		\frac{
			\displaystyle
			S_{i_{(\nu-1)l+1},q}
			\cdots
			S_{i_{\nu l},q}
			S_{j_\nu,q}
			\overline{S_{t_\nu,q}}
		}{
			\displaystyle
			S_{0,q}^{\,l}
		}
		\right).
	\end{eqnarray*}
	\et
	\begin{proof}
		It follows from (6.6) that
		\begin{eqnarray*}
			\mathcal{M}(1,g^s)
			=
			\left\{
			(M^{a_1})^{\otimes l}
			\otimes (M^{a_2})^{\otimes l}\cdots
			\otimes\otimes
			(M^{a_d})^{\otimes l}
			\ \middle|\
			0\leq a_1,\ldots,a_d\leq p
			\right\}.
		\end{eqnarray*}
		For convenience, set $W(a_1,a_2,\ldots,a_d)
		=
		(M^{a_1})^{\otimes l}
		\otimes (M^{a_2})^{\otimes l}
		\otimes\cdots\otimes
		(M^{a_d})^{\otimes l}.$
		For any $i_1,\ldots,i_k,j_1,\ldots,j_d,$ $t_1,\ldots,t_d
		\in\{0,1,\ldots,p\},$
		by Lemma \ref{lem 6.3} and Remark~\ref{remark 6.2}, we have
		\begin{eqnarray*}
			S_{M^{i_1}\otimes M^{i_2}\otimes\cdots\otimes M^{i_k},\,W(a_1,a_2,\ldots,a_d)}
			&=&
			\prod_{\nu=1}^{d}
			\prod_{u=1}^{l}
			S_{i_{(\nu-1)l+u},a_\nu},
			\\[6pt]
			S_{T_{g^s}^{j_1,j_2,\ldots,j_d},\,W(a_1,a_2,\ldots,a_d)}
			&=&
			\prod_{\nu=1}^{d}
			S_{j_\nu,a_\nu},
			\\[6pt]
			\overline{
				S_{T_{g^s}^{t_1,t_2,\ldots,t_d},\,W(a_1,a_2,\ldots,a_d)}
			}
			&=&
			\prod_{\nu=1}^{d}
			\overline{S_{t_\nu,a_\nu}},
			\\[6pt]
			S_{V^{\otimes k},\,W(a_1,a_2,\ldots,a_d)}
			&=&
			\prod_{\nu=1}^{d}
			S_{0,a_\nu}^{\,l}.
		\end{eqnarray*}
		Applying Theorem \ref{thm 5.1}, we obtain
		\begin{eqnarray*}
			N^{\,T_{g^s}^{t_1,\ldots,t_d}}_{
				M^{i_1}\otimes\cdots\otimes M^{i_k},
				T_{g^s}^{j_1,\ldots,j_d}}
			&=&
			\sum_{W\in\mathcal{M}(1,g^s)}
			\frac{
				S_{M^{i_1}\otimes\cdots\otimes M^{i_k},W}
				S_{T_{g^s}^{j_1,\ldots,j_d},W}
				\overline{
					S_{T_{g^s}^{t_1,\ldots,t_d},W}
				}
			}{
				S_{V^{\otimes k},W}
			}
			\\[6pt]
			&=&
			\sum_{a_1,\ldots,a_d=0}^{p}
			\prod_{\nu=1}^{d}
			\frac{
				\displaystyle
				\prod_{u=1}^{l}
				S_{i_{(\nu-1)l+u},a_\nu}
				S_{j_\nu,a_\nu}
				\overline{S_{t_\nu,a_\nu}}
			}{
				\displaystyle
				S_{0,a_\nu}^{\,l}
			}
			\\[6pt]
			&=&
			\prod_{\nu=1}^{d}
			\left(
			\sum_{q=0}^{p}
			\frac{
				\displaystyle
				S_{i_{(\nu-1)l+1},q}
				\cdots
				S_{i_{\nu l},q}
				S_{j_\nu,q}
				\overline{S_{t_\nu,q}}
			}{
				\displaystyle
				S_{0,q}^{\,l}
			}
			\right).
		\end{eqnarray*}

		This completes the proof.
	\end{proof}

	\br{remark 6.8}{\em
		Related fusion rules between untwisted and twisted modules were obtained in
		\cite{dlxy} and further discussed in \cite{dxy2}. Theorem \ref{thm 6.7}
		provides a uniform $S$-matrix expression for these fusion coefficients and
		derives them directly from the twisted Verlinde formula.
	}\er

	We now compute a concrete example. Let $V_{\mathbb{Z}\alpha}$ be the
	rank-one lattice vertex operator algebra associated with the lattice
	$\mathbb{Z}\alpha$, where
	$\langle\alpha,\alpha\rangle=2n$ for some positive integer $n$. Let
	\begin{eqnarray*}
		g=(1\,2\,\cdots\,k)
	\end{eqnarray*}
	be a $k$-cycle. We apply the twisted Verlinde formula to determine the
	fusion rules involving the untwisted and twisted sectors associated with
	the permutation orbifold
	\begin{eqnarray*}
		\left(V_{\mathbb{Z}\alpha}^{\otimes k}\right)^{\langle g\rangle}.
	\end{eqnarray*}

	It is known from \cite{d} (see also \cite{ll}) that the irreducible
	$V_{\mathbb{Z}\alpha}$-modules are precisely
	$V_{\mathbb{Z}\alpha+\lambda_\ell}$, where $\lambda_\ell=\frac{\ell}{2n}\alpha$ and
	$0\leq \ell\leq 2n-1$.
	For convenience, we denote $M^\ell=V_{\mathbb{Z}\alpha+\lambda_\ell}.$
	Therefore,
	\begin{eqnarray*}
		\operatorname{Irr}(V_{\mathbb{Z}\alpha})
		=
		\{M^0,M^1,\ldots,M^{2n-1}\}.
	\end{eqnarray*}
	By \cite{dxy2} (see also \cite{s}), the modular transformation matrix
	of $V_{\mathbb{Z}\alpha}$ is determined by
	\begin{eqnarray}
		S_{\ell,j}
		&:=&
		S_{M^\ell,M^j}
		=
		\frac{1}{\sqrt{2n}}
		e^{-2\pi \sqrt{-1}(\lambda_\ell,\lambda_j)}
		\nonumber\\
		&=&
		\frac{1}{\sqrt{2n}}
		e^{-2\pi \sqrt{-1}
			\left(
			\frac{\ell}{2n}\alpha,
			\frac{j}{2n}\alpha
			\right)}\notag\\
		&=&
		\frac{1}{\sqrt{2n}}
		e^{-\frac{\pi \ell j \sqrt{-1}}{n}} \label{6.7}
	\end{eqnarray}
	for $\ell,j=0,1,\ldots,2n-1 .$

	We first determine the fusion rules involving the untwisted and
	$g^s$-twisted sectors associated with the permutation orbifold.
	Let $1\leq s\leq k-1$. Suppose that $g^s=\tau_1\cdots\tau_d,$
	is products of disjoint cycles, where $d=\gcd(s,k)$ and $k=dl.$
	For any $i_1,\ldots,i_k,\,
	j_1,\ldots,j_d,\,$
	$t_1,\ldots,t_d
	\in\{0,1,\ldots,2n-1\},$ Theorem \ref{thm 6.7} gives
	\begin{eqnarray}
		N^{\,T_{g^s}^{t_1,\ldots,t_d}}_{
			M^{i_1}\otimes\cdots\otimes M^{i_{k}},
			T_{g^s}^{j_1,\ldots,j_d}}
		&=&
		\prod_{\nu=1}^{d}
		\left(
		\sum_{q=0}^{2n-1}
		\frac{
			\displaystyle
			S_{i_{(\nu-1)l+1},q}
			\cdots
			S_{i_{\nu l},q}
			S_{j_\nu,q}
			\overline{S_{t_\nu,q}}
		}{
			\displaystyle
			S_{0,q}^{\,l}
		}
		\right).\label{6.10}
	\end{eqnarray}
	For each fixed $1\leq \nu\leq d$, substituting (\ref{6.7}) into the
	corresponding factor in (\ref{6.10}), we obtain
	\begin{eqnarray}
		&&
		\sum_{q=0}^{2n-1}
		\frac{
			S_{i_{(\nu-1)l+1},q}S_{i_{(\nu-1)l+2},q}\cdots
			S_{i_{\nu l},q}
			S_{j_\nu,q}
			\overline{S_{t_\nu,q}}
		}{
			S_{0,q}^{\,l}
		}
		\notag\\
		&=&
		\frac{1}{2n}
		\sum_{q=0}^{2n-1}
		\exp\left(
		-\frac{2\pi \sqrt{-1}q}{2n}
		\left(
		i_{(\nu-1)l+1}+i_{(\nu-1)l+2}+\cdots+i_{\nu l}+j_\nu-t_\nu
		\right)
		\right).\label{6.11}
	\end{eqnarray}
	By the orthogonality relation for finite Fourier sums, the last
	expression in (\ref{6.11}) is equal to
	\begin{eqnarray*}
		\left\{
		\begin{array}{ll}
			1,
			&
			t_\nu\equiv
			i_{(\nu-1)l+1}+i_{(\nu-1)l+2}+\cdots+i_{\nu l}+j_\nu
			\pmod{2n},
			\\[4pt]
			0,
			&
			\text{otherwise}.
		\end{array}
		\right.
	\end{eqnarray*}
	It follows that
	\begin{eqnarray*}
		N^{\,T_{g^s}^{t_1,\ldots,t_d}}_{
			M^{i_1}\otimes\cdots\otimes M^{i_{k}},
			T_{g^s}^{j_1,\ldots,j_d}}
		=
		\left\{
		\begin{array}{ll}
			1,
			&
			t_\nu\equiv
			i_{(\nu-1)l+1}+\cdots+i_{\nu l}+j_\nu
			\pmod{2n}
			\quad \text{ for }1\leq\nu\leq d,
			\\[6pt]
			0,
			&
			\text{otherwise}.
		\end{array}
		\right.
	\end{eqnarray*}
	
	We next determine the fusion rules involving the $g^s$-twisted  and
	$g^r$-twisted sectors associated with the permutation orbifold.
	Let $1\leq s,r\leq k-1$. Suppose that
	\begin{eqnarray*}
		g^s=\tau_1\tau_2\cdots\tau_d,\quad
		g^r=\sigma_1\sigma_2\cdots\sigma_m,\quad
		g^{s+r}=\gamma_1\gamma_2\cdots\gamma_q
	\end{eqnarray*}
	are products of disjoint cycles, where $d=\gcd(s,k),
	m=\gcd(r,k)$ and
	$q=\gcd(s+r,k).$ Set ${t}=\gcd(s,r,k)$.
	For any $i_1,\ldots,i_d,\,
	j_1,\ldots,j_m,\,
	t_1,\ldots,t_q
	\in\{0,1,\ldots,2n-1\},$ by Theorem \ref{thm 6.5} we have
	\begin{eqnarray*}
		N^{\,T_{g^{s+r}}^{t_1,\ldots,t_q}}_{
			T_{g^s}^{i_1,\ldots,i_d},
			T_{g^r}^{j_1,\ldots,j_m}}
		=
		\prod_{\nu=1}^{{t}}
		\left(
		\sum_{a=0}^{2n-1}
		\frac{
			\displaystyle
			\prod_{u=1}^{\frac{d}{t}}
			S_{i_{(\nu-1)\frac{d}{t}+u},a}
			\prod_{v=1}^{\frac{m}{t}}
			S_{j_{(\nu-1)\frac{m}{t}+v},a}
			\prod_{w=1}^{\frac{q}{t}}
			\overline{S_{t_{(\nu-1)\frac{q}{t}+w},a}}
		}{
			\displaystyle
			S_{0,a}^{\frac{k}{t}}
		}
		\right).
	\end{eqnarray*}
	For each $1\leq \nu\leq {t}$, set
	\begin{eqnarray*}
		P_\nu
		=
		\sum_{u=1}^{\frac{d}{t}}
		i_{(\nu-1)\frac{d}{t}+u}
		+
		\sum_{v=1}^{\frac{m}{t}}
		j_{(\nu-1)\frac{m}{t}+v}
		-
		\sum_{w=1}^{\frac{q}{t}}
		t_{(\nu-1)\frac{q}{t}+w}.
	\end{eqnarray*}
	Using the $S$-matrix formula (\ref{6.7}), the $\nu$-th factor becomes
	\begin{eqnarray*}
		&&
		\sum_{a=0}^{2n-1}
		\frac{
			\displaystyle
			\prod_{u=1}^{\frac{d}{t}}
			S_{i_{(\nu-1)\frac{d}{t}+u},a}
			\prod_{v=1}^{\frac{m}{t}}
			S_{j_{(\nu-1)\frac{m}{t}+v},a}
			\prod_{w=1}^{\frac{q}{t}}
			\overline{
				S_{t_{(\nu-1)\frac{q}{t}+w},a}
			}
		}{
			\displaystyle
			S_{0,a}^{\,\frac{k}{t}}
		}
		\\[6pt]
		&=&
		(2n)^{\frac{k-d-m-q}{2{t}}}
		\sum_{a=0}^{2n-1}
		\exp\left(
		-\frac{2\pi\sqrt{-1}aP_\nu}{2n}
		\right)
		\\[6pt]
		&=&
		\left\{
		\begin{array}{ll}
			(2n)^{1+\frac{k-d-m-q}{2{t}}},
			&
			P_\nu\equiv 0\pmod{2n},
			\\[4pt]
			0,
			&
			P_\nu\not\equiv 0\pmod{2n}.
		\end{array}
		\right.
	\end{eqnarray*}
	Set
	\begin{eqnarray*}
		\eta
		=
		\frac{k+2{t}-d-m-q}{2}.
	\end{eqnarray*}
	Consequently, we have
	\begin{eqnarray*}
		N^{\,T_{g^{s+r}}^{t_1,t_2,\ldots,t_q}}_{
			T_{g^s}^{i_1,i_2,\ldots,i_d},
			T_{g^r}^{j_1,j_2,\ldots,j_m}}
		=
		\left\{
		\begin{array}{ll}
			(2n)^{\eta},
			&
			P_\nu\equiv 0\pmod{2n}
			\quad\text{for }1\leq\nu\leq {t},
			\\[6pt]
			0,
			&
			\text{otherwise},
		\end{array}
		\right.
	\end{eqnarray*}
	where $\eta
	=
	\frac{k+2{t}-d-m-q}{2}.$

	In particular, let $g=(1\,2\,3\,4).$
	For
	\begin{eqnarray*}
		\ell,j,p_1,p_2,q_1,q_2,i_1,i_2,i_3,i_4
		\in\{0,1,\ldots,2n-1\},
	\end{eqnarray*}
	the relevant fusion coefficients are given by
	\begin{eqnarray*}
		N^{\,T_{g^2}^{q_1,q_2}}_{T_g^\ell,T_g^j}
		&=&
		\left\{
		\begin{array}{ll}
			2n,
			&
			q_1+q_2\equiv \ell+j\pmod{2n},
			\\[4pt]
			0,
			&
			q_1+q_2\not\equiv \ell+j\pmod{2n},
		\end{array}
		\right.
		\\[8pt]
		N^{\,T_{g^3}^{j}}_{T_g^\ell,T_{g^2}^{q_1,q_2}}
		&=&
		\left\{
		\begin{array}{ll}
			2n,
			&
			j\equiv \ell+q_1+q_2\pmod{2n},
			\\[4pt]
			0,
			&
			j\not\equiv \ell+q_1+q_2\pmod{2n},
		\end{array}
		\right.
		\\[8pt]
		N^{\,M^{i_1,i_2,i_3,i_4}}_{T_g^\ell,T_{g^3}^j}
		&=&
		\left\{
		\begin{array}{ll}
			1,
			&
			i_1+i_2+i_3+i_4\equiv \ell+j\pmod{2n},
			\\[4pt]
			0,
			&
			i_1+i_2+i_3+i_4\not\equiv \ell+j\pmod{2n},
		\end{array}
		\right.
		\\[8pt]
		N_{T_{g^2}^{p_1,p_2},\,T_{g^2}^{q_1,q_2}}
		^{M^{i_1,i_2,i_3,i_4}}
		&=&
		\left\{
		\begin{array}{ll}
			1,
			&
			i_1+i_2\equiv p_1+q_1 \pmod{2n}
			\ \text{and}\
			i_3+i_4\equiv p_2+q_2 \pmod{2n},
			\\[1mm]
			0,
			&
			\text{otherwise}.
		\end{array}
		\right.
		\\[8pt]
		N^{\,T_{g^r}^j}_{M^{i_1,i_2,i_3,i_4},T_{g^r}^\ell}
		&=&
		\left\{
		\begin{array}{ll}
			1,
			&
			j\equiv \ell+i_1+i_2+i_3+i_4\pmod{2n}
			\\[4pt]
			0,
			&
			j\not\equiv \ell+i_1+i_2+i_3+i_4\pmod{2n}
		\end{array}
		\right.
		\quad \text{for} \,\,r=1,3,
		\\[8pt]
		N^{\,T_{g^2}^{q_1,q_2}}_{
			M^{i_1,i_2,i_3,i_4},T_{g^2}^{\ell,j}}
		&=&
		\left\{
		\begin{array}{ll}
			1,
			&
			q_1\equiv \ell+i_1+i_2\pmod{2n}
			\ \text{and}\
			q_2\equiv j+i_3+i_4\pmod{2n},
			\\[4pt]
			0,
			&
			\text{otherwise}.
		\end{array}
		\right.
	\end{eqnarray*}
	When $n=1$, these fusion rules coincide with those obtained in
	Section~6 of \cite{dxy2}.

\end{document}